# SEIBERG-WITTEN MONOPOLES ON SEIFERT FIBERED SPACES

TOMASZ MROWKA, PETER OZSVÁTH, AND BAOZHEN YU

Abstract. In this paper, we investigate the Seiberg-Witten gauge theory for Seifert fibered spaces. The monopoles over these three-manifolds, for a particular choice of metric and perturbation, are completely described. Gradient flow lines between monopoles are identified with holomorphic data on an associated ruled surface, and a dimension formula for such flows is calculated.

## 1. Introduction

This paper is concerned with Seiberg-Witten gauge theory over Seifert fibered spaces. A Seifert fibered space $\pi\colon Y \to \Sigma$ is the unit circle bundle of an orbifold line bundle over $\Sigma$ whose total space is a smooth three-manifold. The case where $Y$ is diffeomorphic to the product of a circle with a Riemann surface has been studied extensively in [28]; in this paper, we will restrict attention to those Seifert fibered spaces which are not diffeomorphic to a product (indeed, those of non-zero *degree*, in a sense which will be made precise).

Let $i\eta$ denote the connection form of the circle bundle, and $g_\Sigma$ be an orbifold metric over $\Sigma$ with constant curvature, then we can endow $Y$ with the metric

$$g_Y = \eta^2 + \pi^*(g_\Sigma),$$

for which the tangent bundle $TY$ has an orthogonal splitting

$$TY \cong \mathbb{R} \oplus \pi^*(T\Sigma).$$

The Levi-Civita connection on $\Sigma$ then canonically induces a reducible connection ${}^\circ\nabla$ which respects this splitting. We study the solution space to the Seiberg-Witten equations over $Y$, using the above metric and connection on $TY$ (rather than the Levi-Civita connection, which is usually used in the definition of these equations).

Let $\Sigma$ be a two-dimensional orbifold of genus $g$ with marked points $x_1, ..., x_n$ with associated multiplicities $\alpha_1, ..., \alpha_n$. Viewing the marked points as "fractional points,"

The first author was partially supported by an NSF/NYI fellowship (grant number DMS-9357641) and a Sloan Foundation Fellowship.

The second author was partially supported by an NSF Postdoctoral Fellowship; some of the work was done at MSRI, supported in part by NSF grant DMS-9022150.

one can define a natural Euler characteristic by the formula

$$\chi(\Sigma) = 2 - 2g + \sum_{i=1}^{n} (\frac{1}{\alpha_i} - 1), \tag{1}$$

where $g$ is the genus of the underlying smooth curve $|\Sigma|$.

An orbifold line bundle $E$ over $\Sigma$ (hence, in particular, a Seifert fibered space) is specified up to isomorphism by its Seifert data (see Section 2),

$$\mathbf{b}(E) = (b(E), \beta_1(E), ..., \beta_n(E)),$$

a vector of integers with

$$0 \leq \beta_i(E) < \alpha_i.$$

It has an orbifold degree defined by the formula

$$\deg(E) = b(E) + \sum \frac{\beta_i(E)}{\alpha_i}.$$

**Theorem 1.** *Let $Y$ be a Seifert fibered space of non-zero degree over the orbifold $\Sigma$. The moduli space of solutions to the Seiberg-Witten Equations on $Y$ with metric $g_Y$ and connection ${}^\circ\nabla$ is naturally identified with the moduli space of flat $S^1$-connections over $Y$ and two copies of the space of effective orbifold divisors over $\Sigma$ with orbifold degree no greater than $-\frac{\chi(\Sigma)}{2}$.*

**Remark 1.0.1.** *In the above statement, the solutions have not been partitioned according to the* $\mathrm{Spin}_c(3)$ *structures in which they live. See Theorem 5.9.1 for a restatement of Theorem 1 which makes this partitioning explicit.*

In other words, the moduli space of solutions to the Seiberg-Witten Equations consists of a collection of components corresponding to the $S^1$-representation variety of $Y$ (the reducibles): $Hom(\pi_1(Y), S^1)$, which is a certain number of copies of the Jacobian $H^1(\Sigma, \mathbb{R})/H^1(\Sigma, \mathbb{Z})$, and a collection of components $\mathcal{C}^+(\mathbf{e})$ and $\mathcal{C}^-(\mathbf{e})$, labeled by all vectors $\mathbf{e}$ of non-negative integers $(e, \epsilon_1, ..., \epsilon_n)$, with

$$0 \leq \epsilon_i < \alpha_i,$$

and

$$e + \sum_{i=1}^{n} \frac{\epsilon_i}{\alpha_i} \leq -\frac{\chi(\Sigma)}{2}.$$

The spaces $\mathcal{C}^+(\mathbf{e})$ and $\mathcal{C}^-(\mathbf{e})$ both correspond to the space of effective orbifold divisors on $\Sigma$ with multiplicities $\epsilon_i$ over the singular points $x_i$, and a background degree $e$. As such, these spaces are diffeomorphic to the $e^{th}$ symmetric product of the smooth curve $|\Sigma|$, $\mathrm{Sym}^e(|\Sigma|)$.

The orbifold $Y$ has a naturally associated orbifold ruled surface $R$. We have a correspondence between flow-lines between a pair of components in the Seiberg-Witten moduli space of $Y$ and divisors on $R$. To state the correspondence, we must make the following definition.

**Definition 1.0.2.** *Let $\Sigma_-, \Sigma_+ \subset R$ denote the two "sections at infinity" in the ruled surface. Given a pair of orbifold line bundles $E_1$ and $E_2$ a* divisor on $R$ interpolating between $E_1$ and $E_2$ *is an effective divisor $D \subset R$ which does not contain $\Sigma_-$ or $\Sigma_+$ and whose restrictions to these curves are divisors in $E_1$ and $E_2$ respectively.*

**Theorem 2.** *Let $Y$ as in Theorem 1, and fix components $\mathcal{C}^+(\mathbf{e}_1)$, $\mathcal{C}^+(\mathbf{e}_2)$ in the moduli space $\mathcal{M}_{sw}(Y)$. Let $E_1$, $E_2$ be line bundles over $\Sigma$ with Seifert data $\mathbf{e}_1$ and $\mathbf{e}_2$ respectively. Then, there is a natural identification of the moduli space of flows between $\mathcal{C}^+(\mathbf{e}_1)$ and $\mathcal{C}^+(\mathbf{e}_2)$ with the space of divisors in $R$ interpolating between $E_1$ and $E_2$. Similarly for the space of flows between $\mathcal{C}^-(\mathbf{e}_1)$ and $\mathcal{C}^-(\mathbf{e}_2)$. Moreover, there are no flows between $\mathcal{C}^+(\mathbf{e}_1)$ and $\mathcal{C}^-(\mathbf{e}_2)$.*

Similar statements hold for flows which involve the reducible locus (see Theorem 10.0.15).

In general, when $\Sigma$ has singular points, the ruled surface associated to $Y$, $R$, is a singular complex surface, with isolated quotient singularities. These singularities can be resolved, to obtain a smooth surface $\widehat{R}$. The theory of divisors on $R$ is naturally identified with the theory of certain divisors in $\widehat{R}$.

**Theorem 3.** *There is a natural identification of the moduli space of divisors on $R$ interpolating between $E_1$ and $E_2$ with the moduli space of divisors on $\widehat{R}$ in the homology class determined by $E_1$ and $E_2$, via the procedure described in Section 3.*

These identifications enable us to write down a formula for the dimensions of the moduli spaces of flow-lines. To describe the formula, we must introduce some notation. Suppose that the Seifert fibered space $Y$ has Seifert invariants $(b, \beta_1, ..., \beta_n)$. Let $d^j_k$ be the $k^{th}$ denominator in the Hirzebruch-Jung continued fraction expansion of $\alpha_j/\beta_j$, i.e. for a fixed $j$, $d^j_k$ for $k = 0, ..., m_j$ is the sequence of integers which satisfies the initial conditions

$$d^j_0 = \alpha_j,$$

$$d^j_1 = \beta_j,$$

the recurrence relation

$$d^j_k = d^j_{k-1} \left\lfloor \frac{d^j_{k-2}}{d^j_{k-1}} \right\rfloor - d^j_{k-2},$$

and the termination condition that

$$d^j_{m_j} = 1$$

for some integer $m_j$. Here, $\lfloor x \rfloor$ denotes the greatest integer no greater than $x$. Given Seifert data $\mathbf{e} = (e, \epsilon_1, ..., \epsilon_n)$, let $\xi^j_k$ be the minimal decomposition of $\epsilon_j$ with respect to the $d^j_1, ..., d^j_{m_j}$, i.e. for a fixed $j$, the $\xi^j_k$ for $k = 1, ..., m_j$ satisfy the recurrence relation

$$\xi^j_k = \left\lfloor \frac{x - \sum_{i=1}^{k-1} d^j_i \xi^j_i}{d^j_k} \right\rfloor .$$

**Definition 1.0.3.** *Given the above data, we define a rational quantity, the* $Y$-dimension of $\mathbf{e}$, *by the formula*

$$\begin{aligned} \dim_Y(\mathbf{e}) &= \sum_{j=1}^{n} \Big( \sum_{\ell,k=1}^{m_j} d^j_k \xi^j_k d^j_\ell \xi^j_\ell \sum_{i=1}^{\min(k,\ell)} \frac{1}{d^j_{i-1} d^j_i} + \sum_{\ell=1}^{m_j} d^j_\ell \xi^j_\ell \sum_{i=1}^{\ell} \frac{1}{d^j_{i-1} d^j_i} - \sum_{\ell=1}^{m_j} \xi^j_\ell \Big) \\ & \quad + e + \frac{\deg(E)}{\deg(Y)} (\deg(E) - \deg(K_\Sigma)), \end{aligned} \tag{2}$$

*where* $K_\Sigma$ *is the orbifold canonical bundle, and the degree of the orbifold line bundle* $E$ *with Seifert data* $\mathbf{e}$.

Combining Theorems 2 and 3 with the Riemann-Roch formula, we get the following Corollary:

**Corollary 1.0.4.** *The moduli space of flows from* $\mathcal{C}^\pm(\mathbf{e}_1)$ *to* $\mathcal{C}^\pm(\mathbf{e}_2)$ *has dimension given by*

$$\dim_Y(\mathbf{e}_1) + \dim_{Y^{-1}}(\mathbf{e}_2),$$

*where* $Y^{-1}$ *is the inverse of* $Y$, *as a circle bundle over* $\Sigma$.

We specialize the above results first to the case where $Y$ fibers over a Riemann surface $\Sigma$ with no orbifold points, and then to the case where $Y = \Sigma(p, q, r)$ is a Brieskorn sphere (i.e. the Seifert fibered space which is an integral homology sphere fibering over the genus zero orbifold with three marked points with pairwise relatively prime multiplicities $p$, $q$, and $r$).

**Corollary 1.0.5.** *Let* $Y$ *be a non-trivial circle bundle over a Riemann surface* $\Sigma$ *of genus* $g$, *with first Chern number* $n$. *Then, the Seiberg-Witten moduli space is identified with:*

$$\left( \mathbb{Z}/n\mathbb{Z} \times (\mathbb{R}/\mathbb{Z})^{2g} \right) \amalg \left( \coprod_{e=0}^{g-2} \mathcal{C}^+(e) \right) \amalg \left( \coprod_{e=0}^{g-2} \mathcal{C}^-(e) \right),$$

*where now* $e$ *is simply a non-negative integer, and*

$$\mathcal{C}^+(e) \cong \mathcal{C}^-(e) \cong \mathrm{Sym}^e(\Sigma).$$

*The dimension of the space of flows connecting* $\mathcal{C}^+(e_1)$ *with* $\mathcal{C}^+(e_2)$ *is given by*

$$\dim \mathcal{M}(e_1, e_2) = \frac{(e_1 - e_2)}{n}(e_1 + e_2 - (2g-2)) + (e_1 + e_2).$$

When $Y = \Sigma(p, q, r)$ is a Brieskorn sphere, the Seiberg-Witten moduli space is a discrete collection of points, since

$$-\frac{\chi(\Sigma(p,q,r))}{2} = \frac{1}{2}(1 - \frac{1}{p} - \frac{1}{q} - \frac{1}{r}) < 1.$$

Following the usual construction of Floer homology, we can form the relatively graded chain complex freely generated by the irreducible critical points, with boundary maps induced by one-dimensional flow lines (see Section 13). In the Brieskorn sphere case, our results then specialize to the following statement about the homology groups of this complex, the *irreducible Seiberg-Witten Floer homology groups* $\mathrm{HF}^{irr}_*(\Sigma(p,q,r))$:

**Corollary 1.0.6.** *For* $p, q, r$ *pairwise relatively prime, the group* $\mathrm{HF}^{irr}_*(\Sigma(p,q,r))$ *is a free Abelian group generated by two generators for each triple* $(\epsilon_1, \epsilon_2, \epsilon_3)$ *of non-negative integers satisfying*

$$\frac{\epsilon_1}{p} + \frac{\epsilon_2}{q} + \frac{\epsilon_3}{r} < -\frac{\chi(\Sigma(p,q,r))}{2}.$$

*The Floer degree of either solution corresponding to* $(\epsilon_1, \epsilon_2, \epsilon_3)$ *is calculated by* $\dim_{\Sigma(p,q,r)}(0, \epsilon_1, \epsilon_2, \epsilon_3)$

This paper is organized as follows. Sections 2, 3, and 4 are introductory: Section 2 covers some of the basics of orbifold theory, Section 3 discusses some of the elementary properties of the ruled surfaces associated to an orbifold line bundle and its resolution, and Section 4 covers some of the basics of the Seiberg-Witten monopole equations. Section 5 contains the proof of the circle-invariance of monopoles over a Seifert fibered space $Y$, which gives Theorem 1. Section 6 gives the identification of finite-energy gradient flow-lines with holomorphic data over the cylinder. Section 7 shows how to construct a divisor on the ruled surface, given this data. Section 8 is concerned with the converse problem: it shows how a divisor on $R$ gives rise to a gradient flow-line, by reducing the problem to a Kazdan-Warner equation over the cylinder, which is then solved. Thus, Sections 7 and 8 together construct maps between the moduli spaces, which are shown to give induce identifications on the infinitesimal level in Section 9, proving Theorem 2. Section 10 covers the analogue of Theorem 2 for flows involving the reducibles. Section 11 proves the correspondence between orbifold divisors over the singular ruled surface and divisors over its de-singularization, Theorem 3. Section 12 derives the dimension formula in Corollary 1.0.4 from the Riemann-Roch theorem, together with the correspondence with the space of divisors in the de-singularized surface from Section 11. Finally, Section 13 is devoted to some examples.

The authors wish to thank P. B. Kronheimer for inspiring conversations and G. Matić for explaining her work with P. Lisca, which provided a fertile testing-ground for our results.

During the preparation of this manuscript, Baozhen Yu passed away. His untimely death has left us with a profound sense of loss.

## 2. Orbifolds and Seifert Fibered Spaces

A Seifert fibered space is a three-dimensional manifold $Y$ together with an $S^1$ action with finite stabilizers (see [31], [30], [10]); these spaces can be profitably viewed in terms of the theory of orbifold bundles. This section outlines some of the elements of this theory. A more extensive discussion can be found in [10].

Let $D$ denote the standard complex disk, on which $\mathbb{Z}/\alpha\mathbb{Z}$ acts by rotation. An orbifold is a Hausdorff space $|\Sigma|$ with a distinguished finite set of "marked points" $x_1, ..., x_n$ given with integral multiplicities $\alpha_1, ..., \alpha_n$ all greater than 1, equipped with an atlas of coordinate charts

$$\begin{array}{rll} \phi_i : (D,0) & \to (U_i, x_i) & i = 1, ..., n \\ \phi_x : D_x & \to U_x & x \in \Sigma - \{x_1, ..., x_n\}, \end{array}$$

where the $\phi_i$ induce homeomorphisms from $(D,0)/(\mathbb{Z}/\alpha_i\mathbb{Z})$ to $(U_i, x_i)$, the $\phi_x$ are homeomorphisms, the neighborhoods $U_i$ of the distinguished points are pairwise disjoint, and all transition functions (defined over the overlaps) are holomorphic. Of course, when there are no marked points, $\Sigma$ is simply a smooth holomorphic curve, and the theory which will be outlined below is simply gauge theory over that curve. In the presence of marked points, the underlying topological space $|\Sigma|$ still inherits the structure of a smooth complex curve, with local coordinate on $D/(\mathbb{Z}/\alpha_i\mathbb{Z})$-neighborhood of the marked point $x_i$ given by $w^{\alpha_i}$, where $w$ is a local (holomorphic) coordinate on $D$.

An $n$-dimensional *orbifold bundle* is a collection of $\mathbb{Z}/\alpha_i\mathbb{Z}$-equivariant $n$-dimensional vector bundles $E_i$ over $U_i$, and vector bundles $E_x$ over the $U_x$, together with a $1-$cocycle of transition functions over the overlaps. Over each $U_i$, a $\mathbb{Z}/\alpha_i\mathbb{Z}$-equivariant vector bundle is specified up to isometry by giving a representation

$$\rho_i : \mathbb{Z}/\alpha_i\mathbb{Z} \to \mathrm{Gl_n}(\mathbb{C}).$$

The usual notions from gauge theory generalize in a straightforward way to orbifold bundles. For example, an orbifold connection on an orbifold bundle is a collection $\nabla_i$ of $\mathbb{Z}/\alpha_i\mathbb{Z}$-equivariant connections over the disks $U_i$, and connections $\nabla_x$ over the $E_x$, which match up; an orbifold section is a collection of $\mathbb{Z}/\alpha_i\mathbb{Z}$-equivariant sections $\psi_i$ of the $E_i$ and sections $\psi_x$ of $E_x$, all of which match up; and a *holomorphic* orbifold bundle is a collection of (equivariant) holomorphic bundles over the charts, with a $1-$cocycle of holomorphic transition functions.

**Example 2.0.7.** *Holomorphic sections of the trivial orbifold bundle over $\Sigma$ correspond to holomorphic functions over the smooth curve $|\Sigma|$.*

**Example 2.0.8.** *The rotation action of $\mathbb{Z}/\alpha\mathbb{Z}$ on the disk $D$ naturally lifts to an action on the cotangent bundle of $D$, giving the cotangent bundle the structure of a $\mathbb{Z}/\alpha\mathbb{Z}$-equivariant line bundle, or an orbifold line bundle over $D/(\mathbb{Z}/\alpha\mathbb{Z})$. Moreover, any orbifold $\Sigma$ is naturally endowed with a distinguished holomorphic orbifold line bundle the* canonical bundle, *denoted $K_\Sigma$, defined by gluing the cotangent bundles of the $U_i$ and the $U_x$ via the complex derivative of the transition functions.*

**Example 2.0.9.** *Given an orbifold $\Sigma$ and a distinguished point $x_i \in \Sigma$ whose neighborhood $U_i$ is isomorphic to $D/(\mathbb{Z}/\alpha\mathbb{Z})$, we can define a (holomorphic) orbifold line bundle $H_{x_i}$ as follows. The bundle is trivial away from $x_i$, and, over $U_{x_i}$, is given by the $\mathbb{Z}/\alpha_i\mathbb{Z}$-equivariant line bundle $D \times \mathbb{C}$, where $a \in \mathbb{Z}/\alpha_i\mathbb{Z}$ acts by*

$$a \times (w, z) \mapsto (e^{\frac{2\pi \mathbf{i} a}{\alpha_i}} w, e^{\frac{2\pi \mathbf{i} a}{\alpha_i}} z),$$

*with transition function $w$.*

**Definition 2.0.10.** *The* topological Picard group, *denoted* $\mathrm{Pic}^t(\Sigma)$, *is the group of topological isomorphism classes of orbifold line bundles over $\Sigma$, where the group law is tensor product.*

**Lemma 2.0.11.** *The line bundle $H_0$ generates* $\mathrm{Pic}^t(D/(\mathbb{Z}/\alpha\mathbb{Z}))$, *inducing an isomorphism*

$$\mathrm{Pic}^t(D/(\mathbb{Z}/\alpha\mathbb{Z})) \cong \mathbb{Z}/\alpha\mathbb{Z}.$$

*Moreover, multiplying by $w^{-\beta}$, where $w$ is a holomorphic coordinate on $D$, and $0 \leq \beta < \alpha$, gives an identification between holomorphic sections on $H_0^\beta$ and holomorphic functions on $D/(\mathbb{Z}/\alpha\mathbb{Z})$.*

**Proof.** Holomorphic sections of $H_0^\beta$ correspond to holomorphic functions

$$f \colon D \to \mathbb{C}$$

with

$$f(\zeta w) = \zeta^\beta f(w),$$

for each $\zeta$ an $\alpha^{th}$ root of unity; or equivalently, functions of the form

$$f(w) = w^\beta g(w^p),$$

where $g$ is a holomorphic function. □

Given an orbifold line bundle $E$ over a closed orbifold $\Sigma$ with exceptional points $x_1, ..., x_n$ and multiplicities $\alpha_1, ..., \alpha_n$, Lemma 2.0.11 gives a collection of local invariants $\beta_1, ..., \beta_n$ with

$$0 \leq \beta_i < \alpha_i$$

which describe $E$ near the singular points.

**Definition 2.0.12.** *The bundle $E \otimes H_{x_1}^{-\beta_1} \otimes ... \otimes H_{x_n}^{-\beta_n}$ is an orbifold line bundle which is naturally isomorphic to a smooth line bundle over the smooth curve $|\Sigma|$. This bundle is called the* de-singularization *of $E$ and is denoted $|E|$.*

**Definition 2.0.13.** *Given an orbifold line bundle $E$ over $\Sigma$, the collection of integers*

$$(b, \beta_1, ..., \beta_n),$$

*where $b$ is the first Chern number of the de-singularization $|E|$ over $|\Sigma|$, is called the* Seifert invariant *of $E$ over $\Sigma$.*

*The integer $b$ will be called the* background degree *of $E$.*

Lemma 2.0.11 naturally globalizes to give the following:

**Proposition 2.0.14.** *The holomorphic sections of a holomorphic orbifold line bundle $E$ over $\Sigma$ correspond naturally to the holomorphic sections of its de-singularization $|E|$ over $|\Sigma|$.*

Two familiar theorems from complex geometry readily generalize to the orbifold context:

**Theorem 2.0.15.** *(Riemann-Roch) Let $E$ be a holomorphic line bundle over a complex orbifold $\Sigma$, and $\mathcal{E}$ denote its sheaf of holomorphic sections. Then*

$$\chi(\mathcal{E}) - \chi(\mathcal{O}_\Sigma) = b, \tag{3}$$

*where $\chi(\mathcal{F})$ denotes the (complex) Euler characteristic of the sheaf $\mathcal{F}$, $\mathcal{O}_\Sigma$ denotes the sheaf of regular functions on $\Sigma$, and $b$ is the background degree of $E$. Equivalently,*

$$\chi(\mathcal{E}) = 1 - g + b.$$

**Theorem 2.0.16.** *(Serre duality) Let $E \to \Sigma$ be an orbifold line bundle over an orbifold equipped with an orbifold metric. Then the Hodge star operator induces an isomorphism*

$$H^0(\Sigma, E) \cong H^1(\Sigma, K_\Sigma \otimes E^{-1}).$$

An orbifold line bundle $E$ over $\Sigma$ has a naturally associated first Chern number, or *degree* $\deg(E)$, defined by

$$\deg(E) = b + \sum_i \frac{\beta_i}{\alpha_i}.$$

**Example 2.0.17.** *Since the pull-back of $dz$ over $D$ by the map which rotates through $\zeta = e^{2\pi i/\alpha}$ is $\zeta dz$, we see that the local invariants for the canonical bundle are $\alpha - 1$. In fact, since we can naturally identify the canonical bundle of a orbifold $\Sigma$ away from the singular points with the canonical bundle of a smooth curve of genus $g$, we have that the Seifert invariant of $K_\Sigma$ is*

$$(2g-2, \alpha_1 - 1, ..., \alpha_n - 1).$$

*Thus, the degree of $K_\Sigma$ is minus the Euler characteristic of $\Sigma$ (as defined in Equation 1).*

We recall the following result (see [10],[30]), which says that the Seifert invariants of an orbifold line bundle classify the bundle.

**Proposition 2.0.18.** *The map*

$$\mathrm{Pic}^t(\Sigma) \to \mathbb{Q} \oplus \bigoplus_{i=1}^{n} Z/\alpha_i$$

*given by*

$$E \mapsto (\deg(E), \beta_1, ..., \beta_r)$$

*is an injection, with image the set of tuples $(c, \beta_i)$ with*

$$c \equiv \sum \beta_i/\alpha_i \pmod{\mathbb{Z}}.$$

*In particular, if the $\alpha_i$ are mutually coprime, then $\mathrm{Pic}^t(\Sigma) \cong \mathbb{Z}$ is generated by a single line bundle $E_0$ with*

$$\deg(E_0) = \frac{1}{\alpha_1 ... \alpha_n}.$$

Suppose $N$ is an orbifold bundle whose local invariants $\beta_i$ are all relatively prime to $\alpha_i$. Then, the unit circle bundle of $N$, denoted $S(N)$, is naturally a smooth three-manifold. Such a circle bundle is called a *Seifert fibered space.*

In this case, the orbifold invariants of $N$ are reflected in the topological invariants of its unit sphere bundle, according to the following theorem (see [10]).

**Theorem 2.0.19.** *If $S(N)$ is a Seifert fibered space, then for $1 \leq i \leq n$, $1 \leq j \leq g$,*

$$\pi_1(S(N)) = \langle a_i, b_i, g_i, h \big| [a_j, h] = [b_j, h] = [g_i, h] = 1 = g^{\alpha_i} h^{\beta_i} = h^{-b} \prod [a_j, b_j] \prod g_i, \rangle$$

*and*

$$\begin{aligned} H^1(S(N)) &\cong \left\{ \begin{array}{ll} H^1(\Sigma) & \deg(N) \neq 0 \\ H^1(\Sigma) \oplus \mathbb{Z} & \deg(N) = 0, \end{array} \right. , \\ H^2(S(N)) &\cong (\mathrm{Pic}^t(\Sigma)/\mathbb{Z}[N]) \oplus \mathbb{Z}^{2g}. \end{aligned}$$

**Remark 2.0.20.** *The subgroup* $\mathrm{Pic}^t(\Sigma)/\mathbb{Z}[N] \subset H^2(Y;\mathbb{Z})$ *corresponds to the image of the pull-back map*

$$\mathrm{Pic}^t(\Sigma) \xrightarrow{\pi^*} [\textit{Line bundles over } Y] \xrightarrow{c_1} H^2(Y;\mathbb{Z}).$$

*When* $\deg(N) \neq 0$*, this is a finite Abelian group, by Proposition 2.0.18.*

**Corollary 2.0.21.** *If* $N$ *is as above, with* $\deg(N) \neq 0$*, then the representation space of flat* $S^1$ *bundles over* $S(N)$*,* $\mathrm{Hom}(\pi_1(Y), S^1)$*, fits into a split exact sequence*

$$0 \longrightarrow \frac{H^1(|\Sigma|;\mathbb{R})}{H^1(|\Sigma|;\mathbb{Z})} \longrightarrow \mathrm{Hom}(\pi_1(Y), S^1) \xrightarrow{c_1} \mathrm{Pic}^t(\Sigma)/\mathbb{Z}[N] \longrightarrow 0.$$

## 3. The Ruled Surface

In looking at the flow equations for the Seiberg-Witten equation for a three-manifold $Y$, one is naturally led to consider the cylinder $R^o = \mathbb{R} \times Y$. When $Y = S(N)$ has the structure of a Seifert fibered space associated to the orbifold line bundle $N$, the cylinder has a natural compactification, the space obtained by attaching two copies of $\Sigma$ (attached to $Y$ via $\pi$), one at each end of the real line, $+\infty$ and $-\infty$. These curves will be denoted $\Sigma_-$ and $\Sigma_+$ respectively, and the compactified space will be denoted $R$.

Equivalently, the space $R$ can be thought of as the orbifold sphere bundle obtained by projectivizing the (orbifold) complex plane bundle $\mathbb{C} \oplus N$ over $\Sigma$, with $\Sigma_-$, $\Sigma_+$ the projectivizations of $\mathbb{C} \oplus 0$ and $0 \oplus N$ respectively. The projection

$$\pi \colon R \to \Sigma$$

naturally extends the projection map of $\mathbb{R} \times Y$ to $\Sigma$, via the obvious identification of $R - \Sigma_* \cong \mathbb{R} \times Y$, where $\Sigma_* = \Sigma_+ \coprod \Sigma_-$. It is clear from this model that a holomorphic structure on $N$ induces on $R$ the structure of a (possibly singular) complex surface, for which $\pi$ is a holomorphic map.

In fact, the singularities in $R$ are all quotient singularities, which are in a two-to-one correspondence with the marked points of $\Sigma$. More precisely, let $a \in \mathbb{Z}/\alpha\mathbb{Z}$ act on $\mathbb{C}^2$ by

$$a \times (w, z) = (\zeta^a w, \zeta^{\beta a} z),$$

where $\zeta$ is a primitive $\alpha^{th}$ root of unity. Denote the quotient by

$$C_{\alpha,\beta} = \mathbb{C}^2/(\mathbb{Z}/\alpha\mathbb{Z}).$$

Now, if $x \in \Sigma$ is a singular point, then there are two corresponding singular points in the fiber of $\pi$,

$$x_\pm = \Sigma_\pm \cap \pi^{-1}(x),$$

corresponding to the two fixed points of the $S^1$ action on $S^2$. If the local invariant of $N$ at $x$ is $\beta$, then a neighborhood of $x_\pm$ in $R$ is biholomorphic to a neighborhood of the origin in $C_{\alpha,\pm\beta}$.

The minimal resolution, $\widehat{R}$, of $R$ can now be described quite concretely. First recall (see for example [9]; we will revisit this construction in Section 11) that the minimal resolution $\widehat{C}_{\alpha,\beta}$ is constructed by writing $\alpha/\beta$ in its "Hirzebruch-Jung" continued fraction expansion

$$\frac{\alpha}{\beta} = a_1 - \cfrac{1}{a_2 - \cfrac{1}{\ddots - \frac{1}{a_m}}}, \tag{4}$$

(with $a_i \in \mathbb{Z}$, $a_i \geq 2$) which we will abbreviate by writing

$$\frac{\alpha}{\beta} = \langle a_1, ..., a_m \rangle.$$

Then, $\widehat{C}_{\alpha,\beta}$ is obtained by plumbing together $m$ bundles over $S^2$ according to the following diagram:

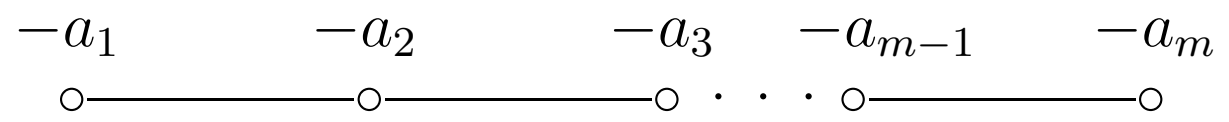

to obtain a manifold with $m$ two-spheres $\{S_i\}_{i=1}^m$, where $S_i$ has self-intersection number $-a_i$ and intersects only $S_{i-1}$ and $S_{i+1}$ (when $i-1$, $i+1$ are still in the range from 1 to $m$) transversally in a single positive point. The map

$$r\colon \widehat{C}_{\alpha,\beta} \to C_{\alpha,\beta}$$

collapses the chain $D = \bigcup_{i=1}^m S_i$ to $[0,0]$, and is biholomorphic away from the chain. Now, the minimal resolution $\widehat{R}$ is obtained by inserting such a configuration of spheres over each of the $2n$ singular points in $R$.

Since the restriction of an orbifold line bundle $N$ over $\Sigma$, $N|_{\Sigma - \{x_1,...,x_n\}}$, is isomorphic to the restriction of $|N|$ to $|\Sigma| - \{x_1, ..., x_n\}$, we have that a Zariski open subset of $R$ is isomorphic to a Zariski open subset of the smooth ruled surface obtained by projectivizing $|N| \oplus \mathbb{C}$ over $|\Sigma|$. Thus, we have:

**Proposition 3.0.22.** *The resolution $\widehat{R}$ is the blowup of a ruled surface fibering over a surface of genus $g$. In particular, if $g = 0$, $\widehat{R}$ is rational.*

## 4. The Seiberg-Witten Equations

We give a brief discussion of the Seiberg-Witten equations, mainly to set up the notation. For a thorough discussion see [29].

Let $Y$ be a Riemannian three-manifold. A $\mathrm{Spin}_c(3)$ structure $(W, \rho)$ is determined by a Hermitian two-plane bundle, the *spinor bundle* $W$ over $Y$, and a positive Clifford module structure

$$\rho\colon T^*Y \to su(W)$$

i.e. a skew-symmetric action of $T^*Y$ on $W$ with satisfying the Clifford relation

$$\rho(\theta) \circ \rho(\theta) = -|\theta|^2 \mathbb{1}_W,$$

for any $\theta \in \Omega^1(Y, \mathbb{R})$, with the property that

$$\rho(\mu_Y) = \mathbb{1}_W,$$

where $\mu_Y$ is the volume form of $Y$, and $\mathbb{1}_W$ is the identity endomorphism of the bundle $W$. In this second identity, we have used the action of $\Lambda^* Y$ on $W$, also denoted $\rho$, extending the given action of $T^*Y$ (this extension exists thanks to the Clifford relation). We will usually abbreviate $\rho(\theta)\Phi$ by $\theta \cdot \Phi$.

Fix an $SO(3)$-connection ${}^\circ\nabla$ on the cotangent bundle $T^*Y$.

**Definition 4.0.23.** *A Hermitian connection $\nabla$ on $W$ is called* spinorial with respect to ${}^\circ\nabla$*, if ${}^\circ\nabla$ is the connection $\nabla$ induces on the cotangent bundle; equivalently, if for all vector fields $X$ and cotangent vector fields $\theta$, we have*

$$[\nabla_X, \rho(\theta)] = -\rho({}^\circ\nabla_X \theta) \tag{5}$$

Given an orthonormal coframe $\{\theta^1, \theta^2, \theta^3\}$ over some patch, and a trivialization for $W$, the connection matrix of any ${}^\circ\nabla$-spinorial connection $\nabla$ can be written with respect to this trivialization as

$$\frac{1}{4} \sum \omega^i_j \otimes \rho(\theta^i \wedge \theta^j) + b\mathbb{1}_W \tag{6}$$

where $\omega^i_j$ are the connection matrices for ${}^\circ\nabla$ with respect to the trivialization $\{\theta^1, \theta^2, \theta^3\}$ of $T^*Y$, and $b \in \Omega^1(Y, i\mathbb{R})$ is any purely imaginary 1−form.

The *pre-configuration space associated to* $W$, written $\mathcal{C}(W)$, is the space

$$\mathcal{C}(W) = \mathcal{A}(W) \times \Gamma(Y, W)$$

formed from the space $\mathcal{A}(W)$ of ${}^\circ\nabla$-spinorial Hermitian connections on $W$, together with the space of spinor fields over $Y$. Since $\mathcal{A}(W)$ has an affine structure for the vector space $\Omega^1(Y, i\mathbb{R})$, $\mathcal{C}(W)$ naturally inherits a manifold structure with

$$T_{(B,\Psi)}\mathcal{C}(W) \cong \Omega^1(Y, i\mathbb{R}) \oplus \Gamma(Y, W).$$

The gauge group $\mathcal{G} \cong \mathrm{Map}(Y, S^1)$ acts on $W$ by scalar multiplication, and hence on the space $\mathcal{C}(W)$ by

$$u \times (\nabla_B, \Psi) = (u\nabla_B u^{-1}, u\Psi),$$

an action whose quotient,

$$\mathcal{B}(W) = \mathcal{C}(W)/\mathcal{G},$$

is the *configuration space associated to* $W$. If $\Phi \not\equiv 0$, one calls the pair $(B, \Psi)$ *irreducible*, and denotes the space of irreducibles by $\mathcal{C}^*(W)$. The gauge group acts freely on the irreducibles. The linearization for the $L^2$ slice condition for the group action at $(B, \Psi)$ is the map

$$T_{(B,\Psi)}\mathcal{C}(W) \to T\mathcal{G}_e \cong \Omega^0(Y, i\mathbb{R})$$

given by

$$(b, \psi) \mapsto d^* b + iIm(\Psi, \psi), \tag{7}$$

where the second term is the imaginary part of the Hermitian inner product of $\Psi$ with $\psi$. With the help of this slice, one can realize the quotient

$$\mathcal{B}^*(W) = \mathcal{C}^*(W)/\mathcal{G}$$

as a Hausdorff Hilbert manifold, after giving $\mathcal{C}(W)$ and $\mathcal{G}$ suitable Sobolev topologies.

Following [22], there is a functional whose gradient flow equations are the four-dimensional Seiberg-Witten equations, and whose critical points are the three-dimensional Seiberg-Witten equations. More precisely, if we choose some reference connection $B_0 \in \mathcal{A}(W)$, there is a functional

$$\mathrm{cs}\colon \mathcal{C}(W) \to \mathbb{R}$$

defined by

$$\mathrm{cs}(B, \Psi) = \int_Y (B - B_0) \wedge \mathrm{Tr}(F_B + F_{B_0}) + \int_Y \langle \Psi, \mathbf{D}_B \Psi \rangle,$$

where

$$\mathbf{D}_B \colon \Gamma(Y, W) \to \Gamma(Y, W)$$

is the Dirac operator of $Y$ induced by the spin connection $B$ on $W$, angle brackets $\langle , \rangle$ here denote the real inner product on $W$ (real part of the Hermitian inner product), and Tr denotes the map induced by taking traces of matrices in $\mathrm{ad}(W)$. If $u \in \mathcal{G}$, then

$$\mathrm{cs}(B, \Psi) - \mathrm{cs}(u(B, \Psi)) = 8\pi^2 \langle c_1(W) \cup [u], [Y] \rangle,$$

where $[u]$ denotes the one-dimensional cohomology class obtained by pulling back the fundamental class on $S^1$. Thus, in general cs descends as a well-defined circle-valued function on $\mathcal{B}(W)$. When, for example $c_1(W)$ is a torsion class, cs is naturally $\mathbb{R}$-valued.

When $\mathbf{D}_B$ is self-adjoint, a property which depends only on $^\circ\nabla$, a point $(B, \Psi) \in \mathcal{B}(W)$ is a critical point for this functional, i.e. the gradient vector field $\nabla \mathrm{cs}$ vanishes at $(B, \Psi)$, if and only if $(B, \Psi)$ satisfies

$$*\mathrm{Tr}(F_A) - i\tau(\Psi) = 0 \tag{8}$$

$$\mathbf{D}_A \Psi = 0, \tag{9}$$

the three-dimensional Seiberg-Witten equations. Here

$$\tau \colon \Gamma(Y, W) \to \Omega^1(Y, \mathbb{R})$$

is adjoint to Clifford multiplication, in the sense that for all $b \in \Omega^1(Y,\mathbb{R})$, $\Psi \in \Gamma(Y,W)$, we have

$$\frac{1}{2}\langle ib\cdot\Psi,\Psi\rangle_W = -\langle b,\tau(\Psi)\rangle_{\Lambda^1}. \tag{10}$$

The critical points for $\nabla$cs in $\mathcal{B}(W)$, the *Seiberg-Witten moduli space in* $W$ is denoted $\mathcal{M}_{sw}(W)$, and its irreducible part is denoted $\mathcal{M}^*_{sw}(W)$. The *Seiberg-Witten moduli space*, $\mathcal{M}_{sw}$ is the union over all distinct $\mathrm{Spin_c}(3)$ structures $W$ of the various moduli spaces $\mathcal{M}_{sw}(W)$.

Given a Riemannian four-manifold $X$ carrying a $\mathrm{Spin_c}(4)$-structure $(W^+,W^-,\rho)$ we can consider the four-dimensional Seiberg-Witten equations

$$\mathrm{Tr}F^+_A - i\sigma(\Phi) = 0 \tag{11}$$

$$\mathbf{D}_A\Phi = 0, \tag{12}$$

where $(A,\Phi)\in\mathcal{A}(W^+)\oplus\Gamma(X,W^+)$. When $X=\mathbb{R}\times Y$ (as Riemannian manifolds) is endowed with the $\mathrm{Spin_c}(4)$-structure induced by $\pi^*(W)$, these equations are gauge equivalent to the (upward) gradient flow lines for cs.

## 5. Circle invariance of Solutions

In this section, we prove Theorem 1. Throughout the section,

$$\pi\colon Y\to\Sigma$$

will denote the Seifert fibered space obtained as the unit circle bundle of an orbifold line bundle $N$ over $\Sigma$. Unless explicitly stated otherwise, we assume that $\deg(N)\neq 0$.

### 5.1. Line Bundles over Seifert Fibered Spaces.

Orbifold line bundles $E$ over $\Sigma$ induce in a natural way (by "pull-back," $\pi^*(E)$) line bundles over $Y$. This does not give a faithful correspondence between isomorphism classes, except when $Y\cong|\Sigma|\times S^1$ (see Remark 2.0.20). However, if one equips the line bundle with a connection, one gets a faithful correspondence, as described below (Proposition 5.1.3).

In order to state properly the image of the correspondence, we must introduce some notions.

**Definition 5.1.1.** *Let $E$ be a line bundle over $Y$. A connection $A$ in $E$ is said to have* trivial fiberwise holonomy *if for any $x\in\Sigma-\{x_i\}_{i=1}^n$, the holonomy of $A$ around $\pi^{-1}(x)$, $\mathrm{Hol_A}\pi^{-1}(\mathrm{x})$, is trivial.*

When the curvature two-form $F(A)$ pulls up from $\Sigma$, i.e. $\iota_{\frac{\partial}{\partial\varphi}}F_A=0$, so that $F_A=\pi^*(F_0)$ (here $\frac{\partial}{\partial\varphi}$ is an infinitesimal generator for the circle action on $Y$), the holonomy around a (generic) fiber is independent of the particular fiber chosen. This is seen by connecting any two points $x, y \in \Sigma - \{x_i\}_{i=1}^n$ by a path $\gamma$ missing the $x_i$, considering its preimage, a cylinder, and noticing that

$$\mathrm{Hol_A}\pi^{-1}\mathrm{x}(\mathrm{Hol_A}\pi^{-1}\mathrm{y})^{-1} = \mathrm{e}^{\int_{\pi^{-1}\gamma}\mathrm{F(A)}} = \mathrm{e}^{\int_\gamma \mathrm{F_0}} = 1.$$

Thus, we have the following:

**Lemma 5.1.2.** *A connection has trivial fiberwise holonomy iff its curvature form pulls up from $\Sigma$, and there is a point $x \in \Sigma - \{x_i\}_{i=1}^n$ whose fiber has trivial holonomy.*

Now we can state the correspondence:

**Proposition 5.1.3.** *There is a natural one-to-one correspondence between pairs (orbifold) bundles-with-connection over $\Sigma$ and (usual) bundles-with-connection over $Y$, whose curvature forms pull up from $\Sigma$ and whose fiberwise holonomy is trivial. Furthermore, this correspondence induces an identification between orbifold sections of the orbifold bundle over $\Sigma$ with fiberwise constant sections of its pull-back over $Y$.*

**Remark 5.1.4.** *It is important that the triviality of the fiberwise holonomy be tested over a point* not *in $\{x_i\}_{i=1}^n$. The holonomy over the exceptional circles $\pi^{-1}(x_i)$ is not expected to be trivial; indeed, for an orbifold line bundle $E$ the holonomy around the exceptional circle $\pi^{-1}(x_i)$ is multiplication by $e^{\frac{2\pi \mathbf{i}}{\alpha_i}(\beta_i(E)-\beta_i(Y))}$.*

**Proof.** We first discuss the correspondence between the bundles-with-connection. Given a bundle-with-connection over $\Sigma$, pull-back clearly induces a bundle-with-connection over $Y$ with trivial fiberwise holonomy. To see that the correspondence is one-to-one, we invert this construction.

Suppose we have a line bundle $E$ over $Y$ with connection $A$ with trivial fiberwise holonomy. We must specify, over each $U_i$, a $\mathbb{Z}/\alpha_i\mathbb{Z}$-equivariant line bundle. Notice that over $U_i$, the projection $\pi$ is modeled on

$$\pi_i \colon D \times_{\mathbb{Z}/\alpha_i\mathbb{Z}} S^1 \to D/(\mathbb{Z}/\alpha_i\mathbb{Z}),$$

so the local charts $\phi_i : D \to D/(\mathbb{Z}/\alpha_i\mathbb{Z})$ defining the orbifold structure on $\Sigma$ can be factored through $\pi_i$, using the local trivializations for the orbifold line bundle defining the Seifert fibered space. Any two such lifts

$$\tau, \tau' \colon D \to D \times_{\mathbb{Z}/\alpha_i\mathbb{Z}} S^1,$$

differ by the action of the circle; i.e. there is a function

$$f \colon D \to S^1$$

uniquely specified by the property that

$$\tau(w) = \tau'(w) f(w).$$

This observation allows us to use the connection $A$ to identify canonically the bundles $\tau^*(E)$ and $\tau'^*(E)$ over $D$. The identification is obtained by taking a homotopy from $f$ to the constant map $w \mapsto 1 \in S^1$ and then using parallel transport of $A$ along the tracks of the homotopy. This identification is canonical, in the sense that is independent of the particular homotopy connecting $f$ with the constant map. This is true because any two identifications differ by parallel transport around the fiber circles (or, more precisely, curves in the fiber of $\pi$ which are homologous to a generic fiber via a cylinder projecting to a path in $\Sigma$) whose holonomy was assumed to be always trivial.

With these observations in place, we can define our orbifold bundle by defining it over $U_i$ to be $\tau^*(E)$, for some lift $\tau$. The $\mathbb{Z}/\alpha_i\mathbb{Z}$ action is induced by the above construction, observing that if $\tau$ is a lift of the projection map $\pi_i$, so is

$$w \mapsto \tau(e^{\frac{2\pi \mathrm{i} a}{\alpha_i}} w),$$

a map we will call $\tau_a$. Thus, we have an identification of $\tau_a^*(E)$ with $\tau^*(E) = \tau_0^*(E)$, giving us the maps covering the $\mathbb{Z}/\alpha_i\mathbb{Z}$ action on $D$. The fact that this map actually induces a $\mathbb{Z}/\alpha_i\mathbb{Z}$ action on $\tau^*(E)$ (i.e. that its $\alpha_i^{th}$ power is trivial) follows from the uniqueness in the above construction, and the observation that $\tau_0 = \tau_{\alpha_i}$.

The statement about the correspondence between sections follows immediately from the definitions. □

5.2. **The geometry of the Seifert fibered space.** We give $Y$ the metric

$$g_Y = \eta^2 + \pi^*(g_\Sigma),$$

where $i\eta$ is a connection form for a constant curvature connection on $Y$, and $g_\Sigma$ is a constant curvature metric on $\Sigma$. The global $1-$form $\eta$ induces a reduction in the structure group of $TY$ to $SO(2)$; the kernel of $\eta$ is naturally identified with the pull-back of the orbifold tangent bundle of $\Sigma$, so that we have orthogonal splittings

$$TY \cong \mathbb{R}\frac{\partial}{\partial\varphi} \oplus \pi^*(T\Sigma) \quad \text{and} \quad T^*Y \cong \mathbb{R}\eta \oplus \pi^*(T^*\Sigma),$$

where $\frac{\partial}{\partial\varphi}$ is the vector field dual to $\eta$. Note that for any vector field $X$ dual to the pull-back of a 1-form from $\Sigma$, we have that

$$[X, \frac{\partial}{\partial\varphi}] = 0. \tag{13}$$

$T^*Y$ can be naturally given a connection compatible with this reduction. Letting $\nabla_\Sigma$ denote the Levi-Civita connection on $T^*\Sigma$, we can give $T^*Y$ the connection ${}^\circ\nabla = d \oplus \pi^*(\nabla_\Sigma)$; i.e. the $SO(3)$-connection which satisfies

$${}^\circ\nabla\eta = 0, \tag{14}$$

and, for any $\theta \in \Omega^1(\Sigma, \mathbb{R})$,

$${}^\circ\nabla\pi^*(\theta) = \pi^*(\nabla_\Sigma\theta). \tag{15}$$

Let ${}^\circ\widehat{\nabla}$ be the Levi-Civita connection on $Y$. The purpose of the next lemma is to compare ${}^\circ\nabla$-spinorial connections on $W$ with ${}^\circ\widehat{\nabla}$-spinorial connections.

Let $\xi$ denote the constant defined by

$$d\eta = 2\xi\mu_\Sigma, \tag{16}$$

where $\mu_\Sigma$ is the volume form on $\Sigma$. By Chern-Weil theory, we see that

$$\xi = -\frac{\pi\deg(N)}{\mathrm{Vol}(\Sigma)}.$$

Then, we have the following:

**Lemma 5.2.1.** *Let $(W,\rho)$ be a $\mathrm{Spin_c}(3)$ structure over $Y$, and $\nabla$ and $\widehat{\nabla}$ be a pair of connections in $W$ which are spinorial with respect to ${}^\circ\nabla$ and ${}^\circ\widehat{\nabla}$ respectively, with $\mathrm{Tr}(\nabla)=\mathrm{Tr}(\widehat{\nabla})$. Then, for any $b\in\Omega^1(Y,\mathbb{R})$, we have*

$$\widehat{\nabla}_{b^\flat} = \nabla_{b^\flat} + \xi(\frac{1}{2}\rho(b) - \langle b,\eta\rangle\rho(\eta)), \tag{17}$$

*where $b^\flat$ denotes the vector field which is $g_Y$-dual to $b$. If $\mathbf{D}$ and $\widehat{\mathbf{D}}$ are the Dirac operators $\nabla$ and $\widehat{\nabla}$ induce on $W$, then*

$$\widehat{\mathbf{D}} = \mathbf{D} - \frac{1}{2}\xi. \tag{18}$$

**Proof.** Pick an orthonormal coframe on a patch of $Y$

$$\{\theta^0,\theta^1,\theta^2\}$$

so that $\theta^0=\eta$, and $\{\theta^1,\theta^2\}$ is pulled back from an orthonormal coframe on $\Sigma$. The connection matrix $(\omega^i_j)$ for the covariant derivative $\widehat{\nabla}$ with respect to this coframe is given by

$$\left(\begin{array}{ccc} 0 & -\xi\theta^2 & \xi\theta^1 \\ \xi\theta^2 & 0 & \xi\theta^0+\omega^1_2 \\ -\xi\theta^1 & -\xi\theta^0-\omega^1_2 & 0 \end{array}\right),$$

which is an easy consequence of Cartan's structural equations,

$$d\theta^i = \sum_j \omega^i_j\wedge\theta^j,$$

and where $\omega^1_2$ is the (pull-back of the) connection $1-$form for the Levi-Civita connection on $\Sigma$. (Notice that the connection matrices here are for $T^*Y$, which are off by a sign from the connection matrices for $TY$.) Thus,

$${}^\circ\widehat{\nabla} - {}^\circ\nabla = \left(\begin{array}{ccc} 0 & -\xi\theta^2 & \xi\theta^1 \\ \xi\theta^2 & 0 & \xi\theta^0 \\ -\xi\theta^1 & -\xi\theta^0 & 0 \end{array}\right).$$

In local coordinates, Formula 6 says that

$$\begin{aligned}\widehat{\nabla} - \nabla &= \frac{1}{4}\sum \omega^i_j \otimes \rho(\theta^i \wedge \theta^j) \\ &= \frac{1}{2}\xi(-\theta^0 \otimes \rho(\theta^0) + \theta^1 \otimes \rho(\theta^1) + \theta^2 \otimes \rho(\theta^2))\rho(\mu_Y),\end{aligned}$$

from which both formulas follow easily. □

**Remark 5.2.2.** *As a consequence of the above result, we have that the Dirac operator* $\mathbf{D}$ *is self-adjoint, since both* $\widehat{\mathbf{D}}$ *and multiplication by a real scalar-valued function are self-adjoint.*

**Corollary 5.2.3.** *The anticommutator of the Dirac operator with Clifford multiplication by a* 1*−form is given by the formula*

$$\{\mathbf{D}, \rho(b)\} = -\rho((*d + d*)b) - 2\nabla_{b^\flat} + 2\xi\langle b, \eta\rangle. \tag{19}$$

**Proof.** This follows immediately from the analogous formula for the Levi-Civita Dirac operator (see for example [4]):

$$\{\widehat{\mathbf{D}}, \rho(b)\} = \rho((*d + d*)b) - 2\widehat{\nabla}_{b^\flat},$$

together with the Lemma. □

5.3. **Weitzenböck Decomposition of the Dirac Operator.** Given a $\mathrm{Spin}_c(3)$-structure $(W, \rho)$ on $Y$, let $\nabla$ be a $^\circ\nabla$-spinorial connection on $W$, where $^\circ\nabla$ is the reducible connection considered above. Corresponding to the splitting of $T^*Y$, there is a canonical decomposition of the Dirac operator into two terms

$$\mathbf{D} = \eta \cdot \nabla_{\frac{\partial}{\partial\varphi}} + \mathrm{D}_2,$$

which can be taken as the defining equation for $\mathrm{D}_2$. Since (taking $L^2$-adjoint)

$$(\eta \cdot \nabla_{\frac{\partial}{\partial\varphi}})^*\psi = \nabla_{\frac{\partial}{\partial\varphi}}\eta \cdot \psi = \eta \cdot \nabla_{\frac{\partial}{\partial\varphi}}\psi + (\nabla_{\frac{\partial}{\partial\varphi}}\eta) \cdot \psi = \eta \cdot \nabla_{\frac{\partial}{\partial\varphi}}\psi$$

(here we have used Equation 14, together with spinoriality, Equation 5), we have decomposed the Dirac operator as a sum of self-adjoint operators. Correspondingly, we can give a Weitzenböck-type decomposition of the square of the Dirac operator:

$$\mathbf{D}^2 = (\eta \cdot \nabla_{\frac{\partial}{\partial\varphi}})^2 + (\mathrm{D}_2)^2 + \{\eta \cdot \nabla_{\frac{\partial}{\partial\varphi}}, \mathrm{D}_2\}, \tag{20}$$

as a sum of three operators, the first two of which are manifestly non-negative. Since the sum of the first two terms on the right has the same symbol as $\mathbf{D}^2$, the third term on the right is *a priori* a first-order operator. In fact, it turns out that this third term is a zeroth order operator, which can be written in terms of Clifford multiplication, as follows.

Let

$$\Pi_{\eta^\perp} \colon T^*Y \to \eta^\perp \subset T^*Y$$

denote orthogonal projection of $T^*Y$ to the space perpendicular to $\eta$, a subbundle which is naturally isomorphic to $\pi^*(T^*\Sigma)$. Then we have the following:

**Lemma 5.3.1.** *The cross-term appearing Formula 20 is a zeroth-order operator. More precisely,*

$$\{\eta \cdot \nabla_{\frac{\partial}{\partial\varphi}}, \mathrm{D}_2\} = -\frac{1}{2}\rho(\Pi_{\eta^\perp} * \mathrm{Tr} F_B).$$

**Proof.** Let $\theta^1, \theta^2$ be a pair of orthonormal covectors over a neighborhood of $Y$ pulled back from (a neighborhood in) $\Sigma$, and $e_1, e_2$ their duals. Then we can write $\mathrm{D}_2 = \theta^1 \cdot \nabla_{e_1} + \theta^2 \cdot \nabla_{e_2}$, so

$$\{\eta \cdot \nabla_{\frac{\partial}{\partial\varphi}}, \mathrm{D}_2\} = \eta \cdot \theta^1 \cdot [\nabla_{\frac{\partial}{\partial\varphi}}, \nabla_{e_1}] + \eta \cdot \theta^2 \cdot [\nabla_{\frac{\partial}{\partial\varphi}}, \nabla_{e_2}].$$

Here we have used the properties of $\nabla$ from Equations 14, 15, and 5, together with the Clifford relations $\eta \cdot \theta^i \cdot = -\theta^i \cdot \eta \cdot$. Now, for general vector fields $u$ and $v$,

$$[\nabla_u, \nabla_v] - \nabla_{[u,v]} = F(u,v);$$

so using Equation 13, we see that

$$\{\eta \cdot \nabla_{\frac{\partial}{\partial\varphi}}, \mathrm{D}_2\} \;=\; \eta \cdot \theta^1 \cdot \widehat{F}_B(\frac{\partial}{\partial\varphi}, e_1) + \eta \cdot \theta^2 \cdot \widehat{F}_B(\frac{\partial}{\partial\varphi}, e_2).$$

Since ${}^\circ\nabla$ pulls up from $\Sigma$, we have that the action of $\widehat{F}_B(\frac{\partial}{\partial\varphi}, e_i)$ for $i = 1, 2$ must commute with Clifford multiplication, hence it must be a scalar endomorphism, so

$$\widehat{F}_B(\frac{\partial}{\partial\varphi}, e_i) = \frac{1}{2}\mathrm{Tr} F_B(\frac{\partial}{\partial\varphi}, e_i).$$

So, by definition, we see that

$$\{\eta \cdot \nabla_{\frac{\partial}{\partial\varphi}}, \mathrm{D}_2\} = \frac{1}{2}\rho(\eta \wedge (\iota_{\frac{\partial}{\partial\varphi}} \mathrm{Tr}(F_B))).$$

The lemma now follows, since for any two-from $\omega$, we have

$$\Pi_{\eta^\perp} * \omega = *\eta \wedge (\iota_{\frac{\partial}{\partial\varphi}} \omega),$$

and

$$\rho(\omega) = -\rho(*\omega),$$

by the assumption that $\rho(\mu_Y)$ is the identity map on $W$. □

5.4. **The Squaring Map.** Since the square of Clifford multiplication by $\pi^*(\mu_\Sigma)$ is the identity map, this operator breaks $W$ into its (orthogonal) $+1$ and $-1$ eigenbundles.

**Definition 5.4.1.** *$S^+$ (resp. $S^-$) is the bundle of spinors $\Psi$ satisfying*

$$i\pi^*(\mu_\Sigma)\cdot\Psi = \Psi \quad (\textit{resp.}\quad i\pi^*(\mu_\Sigma)\cdot\Psi = -\Psi). \tag{21}$$

Clifford multiplication by $\eta$ preserves the splitting $W \cong S^+ \oplus S^-$, so Clifford multiplication by anything in $\eta^\perp \cong \pi^*(T^*\Sigma)$ must reverse it. Indeed, it is straightforward to verify that Clifford multiplication induces an isometry between $\pi^*(T^*\Sigma)$ and the skew-symmetric endomorphisms of $S^\pm$ which reverse chirality or, equivalently, it induces an isometry

$$\pi^*(T^*\Sigma) \cong \mathrm{Hom}_{\mathbb{C}}(S^+, S^-).$$

We will now investigate the manner in which the bilinear map $\tau$ of Equation 10 behaves under this splitting.

**Lemma 5.4.2.** *Let $\alpha \in S^+$, $\beta \in S^-$ be spinors. Then,*

$$\langle i\tau(\alpha+\beta), \eta\rangle = i\frac{1}{2}(|\alpha|^2 - |\beta|^2) \tag{22}$$

$$\Pi_{\eta^\perp} i\tau(\alpha+\beta) = -(\alpha\otimes\beta^* + \beta\otimes\alpha^*). \tag{23}$$

**Proof.** To see the first equation,

$$\begin{aligned}\langle \eta, \tau(\alpha+\beta)\rangle &= -\langle i\rho(\eta)(\alpha+\beta), \alpha+\beta\rangle \\ &= \langle i\rho(*\eta)(\alpha+\beta), \alpha+\beta\rangle \\ &= |\alpha|^2 - |\beta|^2.\end{aligned}$$

For the second, we have

$$\begin{aligned}\langle \tau(\alpha+\beta), i(\alpha^*\otimes\beta + \beta^*\otimes\alpha)\rangle &= \frac{1}{2}\langle(\alpha^*\otimes\beta + \beta^*\otimes\alpha)(\alpha+\beta), \alpha+\beta\rangle \\ &= |\alpha|^2|\beta|^2.\end{aligned}$$

The result follows from this, together with the fact that the transformation $i(\alpha^*\otimes\beta + \beta^*\otimes\alpha) \in \mathrm{Hom}_{\mathbb{C}}(S^+, S^-)$ has norm $|\alpha||\beta|$. □

The Kähler structure on $\Sigma$ induces a canonical (orbifold) $\mathrm{Spin_c}(2)$ structure (see [33],[24]), $\mathbb{C}\oplus K_\Sigma^{-1}$, with Clifford module structure given by the symbol of $\sqrt{2}(\overline{\partial} + \overline{\partial}^*)$. This in turn endows $Y$ with a canonical $\mathrm{Spin_c}(3)$ structure with $W_c = \mathbb{C}\oplus K_\Sigma^{-1}$, given by defining $\rho(\eta)|_{S^\pm} = \pm i$. Moreover, $W_c$ has a canonical ${}^\circ\nabla$-spinorial connection $\nabla_c$, whose associated Dirac operator is given by

$$(\alpha, \beta) \mapsto (-i\frac{\partial}{\partial\varphi}\alpha - \sqrt{2}\overline{\partial}\beta, \sqrt{2}\overline{\partial}^*\alpha + i\frac{\partial}{\partial\varphi}\beta). \tag{24}$$

This $\mathrm{Spin_c}(3)$-structure gives rise to a one-to-one correspondence between the Hermitian line bundles $E$ over $Y$ and the $\mathrm{Spin_c}(3)$ structures over $Y$ under $E \mapsto E \otimes W_c$, a $\mathrm{Spin_c}(3)$ structure whose $S^+$-bundle (in the sense of Definition 5.4.1) is $E$. Indeed, tensor product with $(W_c, \nabla_c)$ induces an identification between $\mathcal{A}(E)$, the space of Hermitian connections in a line bundle $E$, and $\mathcal{A}(W_c \otimes E)$, the space of $^\circ\nabla$-spinorial connections in the $\mathrm{Spin_c}(3)$ structure $W_c \otimes E$.

5.5. **The Vanishing Spinor Argument.** We recall the Kähler vortex equations over $\Sigma$. Given a Hermitian orbifold line bundle $E_0$ over $\Sigma$, let $\alpha_0, \beta_0$, be orbifold sections in $\Gamma(\Sigma, E_0)$ and $\Gamma(\Sigma, K_\Sigma^{-1} \otimes E_0)$ respectively, and $B_0$ be a compatible orbifold connection in $E_0$. The triple $(B_0, \alpha_0, \beta_0)$ is called a *Kähler vortex in* $E_0$ if it satisfies the equations

$$2F_{B_0} - F_{K_\Sigma} = i(|\alpha_0|^2 - |\beta_0|^2)\mu_\Sigma \tag{25}$$

$$\overline{\partial}_{B_0}\alpha_0 = 0 \quad \text{and} \quad \overline{\partial}^*_{B_0}\beta_0 = 0 \tag{26}$$

$$\alpha_0 = 0 \quad or \quad \beta_0 = 0. \tag{27}$$

Once again, these equations are preserved by the gauge group of unitary transformations of $E_0$, so we can form a moduli space of vortices, written $\mathcal{M}_v(\Sigma)$.

The pull-back of a Kähler vortex over $\Sigma$ is naturally the solution of the Seiberg-Witten Equations 8 and 9 over $Y$. The holomorphicity conditions on $\alpha_0$ and $\beta_0$, along with the their circle invariance, ensure that the pull-back spinor is harmonic, according to Equation 24. As explained in Lemma 5.4.2, the squaring map is equivalent to the pair of Equations 22 and 23, the first of which is Equation 25, the second is automatically satisfied because both sides vanish: the left hand side because $F_B$ pulls up from $\Sigma$, the right because of Equation 27.

Since a gauge transformation on $E_0$ pulls back to a gauge transformation of $W$, we get a well-defined map on the level of moduli spaces

$$\pi^*\colon \mathcal{M}_v(\Sigma) \to \mathcal{M}_{sw}(Y).$$

Letting $\mathcal{M}^*_v(\Sigma) \subset \mathcal{M}_v(\Sigma)$ be the subspace represented by vortices for which $\alpha_0$ and $\beta_0$ do not simultaneously vanish, and $\mathcal{M}^*_{sw}(Y) \subset \mathcal{M}_{sw}(Y)$ denote the moduli space of irreducible critical points, we will prove the following:

**Theorem 4.** *The pull-back map $\pi^*$ induces a diffeomorphism between $\mathcal{M}^*_v(\Sigma)$ and $\mathcal{M}^*_{sw}(Y)$.*

**Remark 5.5.1.** *Note that the spaces are identified as differentiable manifolds, not as "moduli spaces"; i.e. whereas $\mathcal{M}^*_v(\Sigma)$ is always cut out by a map whose differential is surjective (according to Theorem 5), $\mathcal{M}^*_{sw}(Y)$ is cut out by a map which always has index zero, as it is cut out by the gradient of a function (which always has a self-adjoint linearization on its critical set).*

The proof occupies the rest of this subsection, and Subsection 5.6. As a first step, we show that $\pi^*$ induces a homeomorphism. First, we see that $\pi^*$ is injective. That is, if a gauge transformation $u$ preserves an irreducible solution which pulls back from $\Sigma$, then it must actually be the pull-back of a gauge transformation over $E_0$. Viewing $u$ as a section of $\mathrm{End}(\pi^*(E_0)) = \pi^*(\mathrm{End}(E_0))$, we have, by the definition of the covariant derivative on the endomorphism bundle, that

$$(\nabla_{\frac{\partial}{\partial\varphi}} u)\Psi = [\nabla_{\frac{\partial}{\partial\varphi}}, u]\Psi = 0.$$

But $\Psi \neq 0$ over a dense open subset of $\Sigma$, so $\nabla_{\frac{\partial}{\partial\varphi}} u = 0$ over all of $Y$. Thus $u$ is the pull-back of a section in $\Gamma(\Sigma, \mathrm{End}(E_0))$, by the correspondence in Proposition 5.1.3.

We will now show that $\pi^*|_{\mathcal{M}_v^*}$ surjects onto $\mathcal{M}_{sw}^*$. Consider any solution $(B, \Psi)$. As above we can decompose $\Psi = \alpha + \beta$ into its components in $S^{\pm}$.

Combining Equation 20, Equation 9, and Equation 5.3.1, we get that

$$0 = (\eta \cdot \nabla_{\frac{\partial}{\partial\varphi}})^2\Psi + (\mathrm{D}_2)^2\Psi - \frac{1}{2}\rho(\Pi_{\eta^\perp} * \mathrm{Tr}F_B)\Psi.$$

Taking the inner product of the above with $\Psi$, and integrating over $Y$, we get

$$0 = |\nabla_{\frac{\partial}{\partial\varphi}}\Psi|^2 + |\mathrm{D}_2{}^2\Psi|^2 - \frac{1}{2}\langle\rho(\Pi_{\eta^\perp} * \mathrm{Tr}F_B)\Psi, \Psi\rangle. \tag{28}$$

Clearly the first two terms are individually non-negative. To analyze the third term, combine

Equation 8, with its interpretation from Lemma 5.4.2, to get that

$$\begin{aligned} -\frac{1}{2}\langle\rho(\Pi_{\eta^\perp}F_B)\Psi, \Psi\rangle &= -\frac{1}{2}\langle i\tau(\alpha+\beta)\alpha, \beta\rangle \\ &= \frac{1}{2}\langle(\alpha\otimes\beta^* + \beta\otimes\alpha^*)(\alpha+\beta), \alpha+\beta\rangle \\ &= |\alpha|^2|\beta|^2, \end{aligned}$$

another manifestly non-negative quantity, which, by Equation 28 must vanish identically.

Thus, we can draw the following conclusions

$$\Pi_{\eta^\perp}(*\mathrm{Tr}F_B) = 0, \tag{29}$$

$$\nabla_{\frac{\partial}{\partial\varphi}}\alpha = \nabla_{\frac{\partial}{\partial\varphi}}\beta = 0 \tag{30}$$

$$|\alpha||\beta| = 0. \tag{31}$$

Since $\Psi \not\equiv 0$, one of $\alpha, \beta \not\equiv 0$, so Equation 30 shows that the fiberwise holonomy in $E_0$ must be trivial. By the correspondence in Proposition 5.1.3, which we can apply thanks to Equations 29 and 30, we see that the sections $\alpha$ and $\beta$ of $S^{\pm}$, together with the connection $B$, correspond to orbifold sections $\alpha_0$, $\beta_0$ of the orbifold bundles $E_0$ and $K_\Sigma^{-1} \otimes E_0$, together with the connection $B_0$ .

Now, the Equation 9 together with Equation 30 gives us Equation 26. The unique continuation theorem for the $\overline{\partial}$ operator, together with Equation 31 gives us Equation 27.

To get an identification between actual moduli spaces, we need to show that $\pi^*$ identifies the kernels, at each point in the respective moduli spaces, of the linearizations of the maps which cut out those moduli spaces.

5.6. **The Linearizations.** We wish to describe the linearization of the map $\nabla\mathrm{cs}$, which cuts out the moduli space $\mathcal{M}_{sw} \subset \mathcal{B}$. First, consider the bilinear map

$$\widehat{\tau}\colon W \otimes W \to (\Omega^1 \oplus \Omega^3)(Y, \mathbb{R})$$

characterized by the property that for any $b \in (\Omega^1 \oplus \Omega^3)(Y, \mathbb{R})$ and $\psi_1, \psi_2 \in \Gamma(Y, W)$,

$$\langle i\rho(b)\psi_1, \psi_2\rangle_W = -\langle b, \widehat{\tau}(\psi_1, \psi_2)\rangle_{\Lambda^1 \oplus \Lambda^3}.$$

The $\Omega^1$ component of $\widehat{\tau}$ is the bilinear form associated to the quadratic map $\tau$, and the $\Omega^3$ component of is the imaginary part of the Hermitian inner product on $W$. Thus, the linearization of $\nabla\mathrm{cs}$ at a solution $(B, \Psi)$ can be described as a map

$$D_{[B,\Psi]}\nabla\mathrm{cs}\colon \Omega^1(Y, i\mathbb{R}) \oplus \Gamma(Y, W) \to \Omega^3(Y, i\mathbb{R}) \oplus \Omega^1(Y, i\mathbb{R}) \oplus \Gamma(Y, W),$$

given by the formula:

$$(b, \psi) \mapsto ((*d + d*)b - i\widehat{\tau}(\Psi, \psi), (\mathbf{D}_B\psi) + b \cdot \Psi). \tag{32}$$

The vortex moduli space, on the other hand, can be thought of as the zero set of a map defined on

$$\mathcal{B}_0^{\pm} = \mathcal{A}_0 \times \Gamma(\Sigma, S^{\pm})/\mathcal{G},$$

$$\mathcal{V}^{\pm}\colon \mathcal{B}_0^{\pm} \to i\Lambda^2(\Sigma) \times_{\mathcal{G}} \Gamma(\Sigma, S^{\mp}),$$

$$\mathcal{V}^{\pm}(B_0, \Psi) = (F_{B_0} \pm i|\Psi|^2, (\mathrm{D}_{2B_0})\Psi).$$

Pulling back gives a map from the vortex configuration space to the Seiberg-Witten configuration space:

$$\pi^*\colon \mathcal{B}_0^{\pm} \to \mathcal{B}.$$

We have shown that $\mathcal{M}^*_{sw}$ is contained in the image of this map; indeed, that $\mathcal{M}^*_{sw} = \pi^*\mathcal{M}^*_v$. We will now show that $\pi^*$ naturally induces identifications

$$\mathrm{Ker} D_{(B,\Psi)}\nabla\mathrm{cs} = \mathrm{Ker} D_{(\pi^*)^{-1}(B,\Psi)}\mathcal{V}^{\pm}, \tag{33}$$

when $(B, \Psi) \in \mathcal{M}^*_{sw}$. More explicitly, this is saying that for any $(B, \Psi) \in \mathcal{M}^*_{sw}$, if $(b, \psi)$ satisfies

$$\mathbf{D}_B\psi + b \cdot \Psi = 0 \tag{34}$$

$$(*d + d*)b - i\widehat{\tau}(\Psi, \psi) = 0, \tag{35}$$

then

$$\Psi \in \Gamma(Y, S^{\pm}) \implies \psi \in \Gamma(Y, S^{\pm}). \tag{36}$$

$$\nabla_{\frac{\partial}{\partial\varphi}}\psi = 0 \tag{37}$$

$$b \in \pi^*(\Omega^1(\Sigma, i\mathbb{R})) \tag{38}$$

We will assume that $\Psi \in \Gamma(Y, S^+)$; the case where $\Psi \in \Gamma(Y, S^-)$ is analogous. Applying $\mathbf{D}_B$ to Equation 34, and applying Equation 19, together with the fact that $\mathbf{D}_B\Psi = 0$, we see that

$$\begin{aligned} &\mathbf{D}_B\mathbf{D}_B\psi + \mathbf{D}_B(b \cdot \Psi) \\ &= \mathbf{D}_B\mathbf{D}_B\psi + \mathbf{D}_B(b \cdot \Psi) \\ &= \mathbf{D}_B^2\psi - (*d + d*)b \cdot \Psi - 2\nabla_{b^\flat}\Psi + 2\xi\langle b, \eta\rangle\eta \cdot \Psi. \end{aligned} \tag{39}$$

Let $\psi = \alpha + \beta$ be the decomposition of $\psi$ in $\Gamma(Y, S^{\pm})$. We will show that $\beta = 0$, by taking the inner product of the above equation with $\beta$ and making four observations.

1. By Equation 29 combined with Equation 20, we see that $\mathbf{D}_B^2$ preserves the grading on $S^{\pm}$, so
$$\langle \mathbf{D}_B^2\psi, \beta\rangle = \langle \mathbf{D}_B^2\beta, \beta\rangle = \langle (\eta \cdot \nabla_{\frac{\partial}{\partial\varphi}})^2\beta, \beta\rangle + \langle (\mathrm{D}_2)^2\beta, \beta\rangle.$$
2. Combining Equation 35 with the definition of $\widehat{\tau}$, we see that
$$\begin{aligned} \langle -(*d + d*)b \cdot \Psi, \beta\rangle &= -\langle i\widehat{\tau}(\Psi, \psi) \cdot \Psi, \beta\rangle \\ &= \langle \widehat{\tau}(\Psi, \psi), \widehat{\tau}(\Psi, \beta)\rangle \\ &= \langle \widehat{\tau}(\Psi, \beta), \widehat{\tau}(\Psi, \beta)\rangle. \end{aligned}$$
The last equality follows from the fact that $\widehat{\tau}(\Psi, \alpha)$ preserves the splitting of $W$ into $S^{\pm}$, whereas $\widehat{\tau}(\Psi, \beta)$ reverses it, so they must be orthogonal.
3. Since the splitting of $W$ into $S^{\pm}$ is induced by the eigenspaces of Clifford multiplication by a covariantly constant 2-form, $i\mu_\Sigma$, spinoriality of $\nabla$ ensures that covariant differentiation respects the splitting. Also, Clifford multiplication by $\eta$ preserves the spitting, so
$$\langle -2\nabla_{a^*}\Psi + 2\xi\langle a, \eta\rangle\eta \cdot \Psi, \beta\rangle = 0.$$

Taking the inner product of Equation 39 with $\beta$ and applying these observations (the vanishing of four *a priori* non-negative terms), we can conclude that

$$\langle \widehat{\tau}(\Psi, \beta), \widehat{\tau}(\Psi, \beta)\rangle = 0;$$

but this implies that $\beta = 0$, and we have $\psi \in \Gamma(Y, S^+)$ as in Equation 36.

Next, we turn our attention to deriving Equation 37. Since, as we have just established, both $\Psi, \psi \in \Gamma(Y, S^+)$, it follows that $\widehat{\tau}(\Psi, \psi)$ must preserve the splitting of $S^{\pm}$, so that $\Pi_{\eta^\perp}\widehat{\tau}(\Psi, \psi) = 0$; i.e. by Equation 35,

$$\iota_{\frac{\partial}{\partial\varphi}} db = 0.$$

Hence, $db$ pulls up from $\Sigma$, and, once again by Equation 35, together with Equation 15, we get that

$$0 = {}^{\circ}\nabla_{\frac{\partial}{\partial\varphi}}\widehat{\tau}(\Psi,\psi);$$

which implies, after taking the $\eta$-component, that

$$0 = \frac{\partial}{\partial\varphi}\langle\Psi,\psi\rangle. \tag{40}$$

Decompose $b$ according to the splitting of the cotangent bundle,

$$b = f\eta + c,$$

i.e. $f \in \Omega^0(Y, i\mathbb{R})$ and $c \in \Omega^1(Y, i\mathbb{R})$. Then, the $S^+$-component of Equation 34 gives us that

$$i\nabla_{\frac{\partial}{\partial\varphi}}\psi + if\Psi = 0. \tag{41}$$

Hitting this equation with $i\nabla_{\frac{\partial}{\partial\varphi}}$ and integrating against $\psi$, we get

$$\begin{aligned} 0 &= \langle i\nabla_{\frac{\partial}{\partial\varphi}}\psi, i\nabla_{\frac{\partial}{\partial\varphi}}\psi\rangle - \int(\frac{\partial}{\partial\varphi}f)\langle\Psi,\psi\rangle = \langle i\nabla_{\frac{\partial}{\partial\varphi}}\psi, i\nabla_{\frac{\partial}{\partial\varphi}}\psi\rangle + \int f\frac{\partial}{\partial\varphi}\langle\Psi,\psi\rangle \\ &= \langle i\nabla_{\frac{\partial}{\partial\varphi}}\psi, i\nabla_{\frac{\partial}{\partial\varphi}}\psi\rangle, \end{aligned}$$

by Equation 40, together with the fact that $\nabla_{\frac{\partial}{\partial\varphi}}\Psi = 0$. Hence, we have obtained Equation 37.

To show that the form $b$ pulls up from $\Sigma$, notice that Equation 41 combined with Equation 37 shows that $f = 0$ on the (dense, open) set where $\Psi \neq 0$, so it must vanish identically. Thus,

$$\iota_{\frac{\partial}{\partial\varphi}}b = 0.$$

Combined with the fact that $db$ pulls up from $\Sigma$, we get that $b$ itself must pull up, giving Equation 38.

This completes the proof of Equation 33.

5.7. **Vortices over Orbifolds.** The complex interpretation of irreducible solutions in $\mathcal{M}^*_{sw}(Y)$ appearing in Theorem 1 follows from Theorem 4, together with the following theorem, which is essentially contained in [17],[5], about vortices over orbifolds.

**Definition 5.7.1.** *A vortex with $\beta_0 = 0$ (resp. $\alpha_0 = 0$) is called a* positive vortex *(resp. negative vortex). The moduli space of positive vortices (resp. negative vortices) in a line bundle $E$ over $\Sigma$ will be denoted $\mathcal{M}^+_v(E)$ (resp $\mathcal{M}^-_v(E)$).*

**Theorem 5.** *Let $E \to \Sigma$ be an orbifold line bundle with background Chern number $e$. The moduli space of vortices $\mathcal{M}^+_v(E) \cong \mathcal{M}^-_v(K_\Sigma \otimes E^{-1})$ is empty if*

$$\deg(E) > \frac{\deg(K_\Sigma)}{2},$$

*and it is naturally diffeomorphic to the e-fold symmetric product of* $|\Sigma|$, $\mathrm{Sym}^e(|\Sigma|)$, *if*

$$\deg(E) < \frac{\deg(K_\Sigma)}{2}. \tag{42}$$

**Proof.** For the proof, we refer to [5],[17]. We give the outline of the argument.

The identification of $\mathcal{M}_v^+(\Sigma, E)$ with $\mathcal{M}_v^-(\Sigma, E)$ arises from Serre duality, which gives a map

$$*\colon \Gamma(\Sigma, E) \to \Gamma(\Sigma, K_\Sigma^{-1} \otimes (K_\Sigma \otimes E^{-1})),$$

with

$$\overline{\partial}^* * = - * \overline{\partial}.$$

Since a line bundle with negative degree supports no holomorphic sections, the vortex equations force any bundle $E$ with a non-zero vortex to satisfy Inequality 42.

When the Inequality 42 is satisfied, then the moduli space of vortices is identified with a moduli space of divisors, which consist of pairs $(B_0, \alpha_0)$ satisfying only $\overline{\partial}_{B_0}\alpha_0 = 0$, modulo the action of $\mathrm{Map}(\Sigma, \mathbb{C}^*)$. This identification is done by showing that for each divisor $(\overline{\partial}_{B_0}, \alpha_0)$, there is a unique real-valued function $u$ on $\Sigma$ such that $e^u(\overline{\partial}_{B_0}, \alpha_0)$ satisfies the curvature condition of Equation 25, by reexpressing this latter equation as a Kazdan-Warner equation for $u$. Properties of this equation guarantee the existence and uniqueness of $u$.

Finally, such a divisor is uniquely determined (as the name suggests) by the zero-set $|\alpha|^{-1}(0)$, with zeros counted with appropriate multiplicities. (Here $|\alpha|$ denotes the holomorphic section of the desingularization of $(E, \overline{\partial}_{B_0})$.) Over a (one-dimensional) orbifold, the moduli space of divisors is always cut out transversally: it is cut out by a map whose linearization has kernel $H^0(\alpha^{-1}(0), \mathcal{E}|_{\alpha^{-1}(0)})$, and and cokernel $H^1(\alpha^{-1}(0), \mathcal{E}|_{\alpha^{-1}(0)})$. The cokernel, being the first cohomology group of a sheaf over a zero-dimensional space, must vanish, and the index of the operator is given by the difference

$$\chi(\Sigma, \mathcal{E}) - \chi(\Sigma, \mathcal{O}_\Sigma) = e,$$

by the Riemann-Roch formula. □

5.8. **Reducibles.** We now turn our attention to the reducible critical points in $Y$. Reducible critical points in the $\mathrm{Spin}_c(3)$ structure determined by $E$ correspond to gauge equivalence classes of connections $B$ whose curvature form satisfies

$$2F_B - F_{K_\Sigma} = 0. \tag{43}$$

**Definition 5.8.1.** *A line bundle* $E_0$ *is called an* orbi-spin bundle on $\Sigma$ *if*

$$2\deg(E_0) = \deg(K_\Sigma),$$

Not every orbifold admits an orbi-spin bundle. For example, among the orbifolds with cyclic topological Picard group (i.e. orbifolds for which all the local invariants $\alpha_i$ are pairwise relatively prime), the ones which do admit an orbi-spin bundle are precisely those whose local invariants are all odd.

**Proposition 5.8.2.** *A line bundle $E$ over $Y$ supports a reducible solution if and only if $E \cong \pi^*(E_0)$ for some orbifold line bundle $E_0$ over $\Sigma$. In this case, the space of reducibles in the* $\mathrm{Spin_c}(3)$ *structure determined by $E$, denoted $\mathfrak{J}(E)$, is homeomorphic to the Jacobian torus $H^1(\Sigma;\mathbb{R})/H^1(\Sigma;\mathbb{Z})$. A line bundle $E$ supports a reducible solution with trivial fiberwise holonomy if and only if $E \cong \pi^*(E_0)$ for an orbi-spin bundle $E_0$ over $\Sigma$.*

**Proof.** Suppose $E \cong \pi^*(E_0)$. One can always find a constant curvature connection $B_0$ on $E_0$. Then, the connection

$$\pi^*(B_0) + i\left(\frac{\frac{1}{2}\deg(K_\Sigma) - \deg(E_0)}{\deg(N)}\right)\eta$$

satisfies Equation 43. Conversely, if $E$ supports a reducible solution, then it follows from Equation 43 and Theorem 2.0.19 that $c_1(E)$ is a torsion class, hence it must be isomorphic to $\pi^*(E_0)$ for an orbifold line bundle $E_0$ over $\Sigma$.

Any two connections $B_1, B_2$ in line bundles $E_1, E_2$ respectively, with the same curvature form must differ by tensoring with a flat line bundle. Thus, the identification of the reducibles with the Jacobian follows from Corollary 2.0.21.

The final statement is a direct consequence of Proposition 5.1.3. □

We wish to investigate conditions under which the identification between $\mathfrak{J}(E)$ and the Jacobian torus of $\Sigma$ is naturally a diffeomorphism.

**Definition 5.8.3.** *A point $[B,0] \in \mathfrak{J}(E)$ is called a* non-degenerate critical point *if the kernel of the linearization of the Seiberg-Witten equations $D_{[B,0]}\nabla\mathrm{cs}$ is isomorphic to the tangent space to the Jacobian torus. The reducible locus $\mathfrak{J}(E)$ is called* non-degenerate *when each reducible $[B,0] \in \mathfrak{J}(E)$ is a non-degenerate critical point.*

Around a reducible solution $(B,0)$, the linearization of the Seiberg-Witten equations takes the form

$$D_{[B,0]}\nabla\mathrm{cs}(b,\psi) = ((*d + d*)b, \mathbf{D}_B\psi),$$

so its kernel is isomorphic to $H^1(Y;\mathbb{R}) \oplus \mathrm{Ker}\mathbf{D}_B$; so $\mathfrak{J}(E)$ is non-degenerate when for all $[B,0] \in \mathfrak{J}(E)$, $\mathbf{D}_B$ has no kernel.

**Proposition 5.8.4.** *Let $(E,B)$ be a line-bundle-with-connection over $Y$ whose curvature form pulls back from $\Sigma$. If $B$ has non-trivial fiberwise holonomy, then the*

*Dirac operator $\mathbf{D}_B$ has no kernel. Otherwise, by Proposition 5.1.3, we can push $(E, B)$ forward to $\Sigma$, and we have a corresponding identification*

$$\mathrm{Ker}(\mathbf{D}_B) \cong H^0(\Sigma, \pi_*(E,B)) \oplus H^1(\Sigma, \pi_*(E,B)).$$

**Proof.** This follows from the Weitzenböck formula from Equation 20, together with the hypothesis that $F_B$ pulls back from $\Sigma$ and Lemma 5.3.1 give that the square of the Dirac operator splits as a pair of non-negative operators

$$\mathbf{D}_B^2 = (\eta \cdot \nabla_{\frac{\partial}{\partial \varphi}})^2 + (\mathrm{D}_2)^2,$$

so any kernel element must be fiberwise constant. Hence the first statement. The second follows from the complex interpretation of the Dirac operator, Equation 24. □

It follows from the above Proposition that if a line bundle $E$ admits a reducible solution $B$ whose Dirac operator has a non-trivial kernel, then $E \cong \pi^*(E_0)$ for an orbi-spin bundle $E_0$. When $\Sigma$ has cyclic topological Picard group, the orbifold line bundles over $\Sigma$ are uniquely specified by their degree. Hence, there is at most one orbi-spin bundle, which we write as $K_\Sigma^{1/2}$.

**Corollary 5.8.5.** *Suppose that $Y$ is a Seifert fibered space over an orbifold $\Sigma$ whose topological Picard group is cyclic, and $E$ is the line bundle corresponding to a $\mathrm{Spin}_c(3)$-structure which supports reducible solutions. Then, the reducible locus $\mathfrak{J}(E)$ is not non-degenerate if and only if the following three conditions are satisfied:*

1. *$\Sigma$ is orbi-spin,*
2. *$E \cong \pi^*(K^{1/2})$, and*
3. *$g(\Sigma) > 0$,*

**Proof.** By Proposition 5.8.4, if $\mathbf{D}_B$ has non-trivial kernel, $B$ must have non-trivial fiberwise holonomy, so that $\pi_*(E,B)$ is a bundle-with-connection over $\Sigma$. Equation 43 then implies that $\pi_*(E,B)$ lives in an orbi-spin bundle on $\Sigma$, proving Conditions 1 and 2. The Seifert invariants of the orbi-spin $K^{1/2}$ bundle are

$$(g-1, \frac{\alpha_1 - 1}{2}, ..., \frac{\alpha_n - 1}{2})$$

(all these numbers are integers, since all the $\alpha_i$ are odd); and as $B$ ranges over all connections with a fixed curvature, $\pi_*(E,B)$ ranges over all possible complex structures on $K^{1/2}$. According to Proposition 5.8.4, together with Serre duality (Theorem 2.0.16), the question of whether there is a reducible $B$ on $E$ with $\ker \mathbf{D}_B \neq 0$ is equivalent to the question of whether there is a complex structure on $K^{1/2}$ with holomorphic sections. This is in turn is topological question: $K^{1/2}$ has holomorphic sections if and only if its desingularization has non-negative degree (using Proposition 2.0.14 and the corresponding well-known fact for smooth curves), i.e. $g-1 \geq 0$. □

**Remark 5.8.6.** *If $g(\Sigma)=0$ and $\Sigma$ does not have a cyclic topological Picard group, $Y$ might admit a reducible solution $(B,0)$ whose Dirac operator has non-trivial kernel. Let $\Sigma$ be the genus zero orbifold with three marked points, each having multiplicity 3. Since $deg(K_\Sigma)=0$, the trivial bundle is an orbi-spin bundle with a holomorphic section.*

5.9. **The critical sets of the Chern-Simons functional.** Combining Theorem 4, Theorem 5, and Proposition 5.8.2, we obtain the following refined version of Theorem 1.

**Theorem 5.9.1.** *Let $\pi\colon Y = S(N) \to \Sigma$ be a Seifert fibered space corresponding to an orbifold line bundle $N$ with non-zero degree. The Seiberg-Witten moduli space for the $\mathrm{Spin_c}(3)$ structure $E \otimes W_c$ is non-empty if and only if $E \cong \pi^*(E_0)$ for some orbifold line bundle $E_0$ over $\Sigma$. The moduli space $\mathcal{M}_{sw}(\pi^*(E_0)\otimes W_c)$ consists of one reducible component $\mathfrak{J}(\pi^*(E_0))$, which is homeomorphic to the Jacobian torus of $|\Sigma|$, and a pair of components $\mathcal{C}^\pm(E)$ for each isomorphism class of orbifold line bundles $E$ over $\Sigma$ with*

$$0 \leq \deg(E) < \frac{\deg(K_\Sigma)}{2}$$

*and*

$$[E] \equiv [E_0] \pmod{\mathbb{Z}[N]}$$

*in $\mathrm{Pic}^t(\Sigma)$. The component $\mathcal{C}^+(E) \cong \mathcal{C}^-(E)$ is naturally diffeomorphic to the moduli space of effective divisors in $E$.*

**Remark 5.9.2.** *In the above statement, $\mathcal{C}^+(E)$ is obtained as $\pi^*(\mathcal{M}_v^+(E))$, while $\mathcal{C}^-(K_\Sigma \otimes E^{-1})$ is $\pi^*(\mathcal{M}_v^-(E))$. Sometimes (as in Section 1), we find it convenient to label the irreducible critical manifolds by Seifert data rather than isomorphism classes of line bundles (these data are equivalent, according to Proposition 2.0.18); i.e. $\mathcal{C}^\pm(\mathbf{e})$ denotes $\mathcal{C}^\pm(E)$, when $\mathbf{e}$ is the Seifert data for the line bundle $E$. Moreover, when it is clear from the context, we will drop the line bundle (indicating the $\mathrm{Spin_c}(3)$ structure) from the notation for the reducible locus $\mathfrak{J}$.*

Recall the following definition from Morse theory (which is the natural analogue of Definition 5.8.3 in the irreducible case).

**Definition 5.9.3.** *A critical point $[B,\Psi] \in \mathcal{M}^*_{sw}(E)$ is called a* non-degenerate critical point *if there is a neighborhood of $[B,\Psi]$ which is a submanifold of of $\mathcal{B}(E)$, and*

$$\mathrm{Ker} D_{[B,\Psi]} \nabla \mathrm{cs} = T_{[B,\Psi]} \mathcal{M}^*_{sw}(E).$$

*A component* $\mathcal{C} \subset \mathcal{M}^*_{sw}(E)$ *is called a* non-degenerate critical manifold *if every point* $[B,\Psi] \in \mathcal{C}$ *is a non-degenerate point.*

**Corollary 5.9.4.** *The irreducible critical manifolds* $\mathcal{C}^{\pm}(E)$ *are all non-degenerate.*

**Proof.** This follows from the fact that the moduli space of vortices is always a smooth manifold, cut out transversally according to Theorem 5, together with the identification of the kernels in Equation 33. □

The Chern-Simons invariant, cs, of a solution depends only on the component of that solution in the moduli space. Moreover, since only torsion classes arise as the first Chern class of a $\mathrm{Spin}_c(3)$ structure for which cs has critical values, we see that the Chern-Simons function is naturally a real-valued function (we are using here the assumption that $Y$ is associated to a line bundle with non-zero degree, see Remark 2.0.20).

**Proposition 5.9.5.** *Suppose that the function* cs *is normalized so that* $\mathrm{cs}(\mathfrak{J}) = 0$ *for* $\mathfrak{J}$ *a reducible solution to the Seiberg-Witten equations. Then,*

$$\mathrm{cs}(\mathcal{C}^{\pm}(E)) = \frac{4\pi^2}{\deg(N)}\left(\deg(E) - \frac{\deg(K_\Sigma)}{2}\right)^2. \tag{44}$$

**Proof.** Suppose $B$ is a Hermitian connection in a line bundle $E$ over $Y$. We define the *$Y$-degree of $E$*, $\deg_Y(B)$ by the formula

$$\deg_Y(B) = \frac{i}{4\pi^2}\int F_B \wedge \eta.$$

Let $B_1, B_2$ be two Hermitian connections in a line bundle $E$ over $Y$ whose curvature forms pull up from $\Sigma$. (This condition is equivalent to the condition that the fiberwise holonomy over the smooth locus of $\Sigma$ is constant, see the remarks preceding the proof of Proposition 5.1.3). We have that

$$\int_Y (B_1 - B_2) \wedge \mathrm{Tr}(F_{B_1\otimes W_c} + F_{B_2\otimes W_c}) \tag{45}$$

$$= \int_Y (B_1 - B_2) \wedge \pi^*(2F_{B_1} + 2F_{B_2} - 2F_{K_\Sigma})$$

$$= \int_\Sigma \pi_*(B_1 - B_2)(2F_{B_1} + 2F_{B_2} - 2F_{K_\Sigma}), \tag{46}$$

where

$$\pi_* : \Omega^1(Y) \to \Omega^0(\Sigma)$$

denotes the operation of integration along the fibers. Now, $\pi_*(B_1 - B_2)$ is a constant function over $\Sigma$, since $\pi_*(B_1 - B_2)$ measures the difference between the fiberwise holonomy of $B_1$ and the fiberwise holonomy of $B_2$, and both holonomies are constant. In fact,

$$\pi_*(B_1 - B_2) \equiv \frac{2\pi i(\deg_Y(B_1) - \deg_Y(B_2))}{\deg(N)},$$

since

$$d\eta = -\frac{2\pi \deg(N)}{\mathrm{Vol}(\Sigma)}\mu_\Sigma$$

so that

$$\begin{aligned}
\int_\Sigma \pi_*(B_1 - B_2)\mu_\Sigma &= \int_Y (B_1 - B_2)\pi^*(\mu_\Sigma) \\
&= -\frac{\mathrm{Vol}(\Sigma)}{2\pi \deg(N)} \int_Y (B_1 - B_2) \wedge d\eta \\
&= -\frac{\mathrm{Vol}(\Sigma)}{2\pi \deg(N)} \int_Y (F_{B_1} - F_{B_2}) \wedge \eta \\
&= \frac{2\pi i \mathrm{Vol}(\Sigma)}{\deg(N)}(\deg_Y(B_1) - \deg_Y(B_2)).
\end{aligned}$$

Thus, we have that

$$\int_Y (B_1 - B_2) \wedge \mathrm{Tr}(F_{B_1 \otimes W_c} + F_{B_2 \otimes W_c}) = \frac{8\pi^2}{\deg(N)} \left(\deg_Y(B_1) - \deg_Y(B_2)\right)^2.$$

Equation 44 follows from this equation, since Chern-Weil theory gives us that for any line-bundle-with-connection $(E_0, C_0)$ over $\Sigma$,

$$\deg_Y(\pi^*(E_0, C_0)) = \deg(E_0)$$

and Equation 43 guarantees that for any reducible solution $(B, 0)$,

$$\deg_Y(B) = \frac{\deg(K_\Sigma)}{2}.$$

□

## 6. Gradient Flow Lines

We now turn our attention to solutions to the Seiberg-Witten equations (Equations 11, 12) over the cylinder $R^o \cong \mathbb{R} \times Y$, given the metric

$$g_{cyl} = dt^2 + g_Y$$

and $\mathrm{Spin_c}(4)$ structure $p^*(W)$ induced by some given $\mathrm{Spin_c}(3)$ structure $(W, \rho)$ over $Y$. The configuration space now is given by

$$\mathcal{B}(p^*(W)) = \mathcal{A}(p^*(W)) \times \Gamma(R^o, p^*(W))/\mathrm{Map}(R^o, S^1),$$

where $\mathcal{A}(p^*(W))$ is the affine space of spinorial connections in $p^*(W)$, meaning those connections $\nabla$ which are spinorial with respect to the connection which is $\nabla_{cyl} = d \oplus {}^o\nabla$ in terms of the splitting $T^*R^o = \mathbb{R}(dt) \oplus T^*Y$. By elliptic regularity, the topology induced on the Seiberg-Witten moduli space from the $\mathcal{C}^\infty$ topology on $\mathcal{B}(p^*(W))$ is equivalent to the topology induced by the $L^2_k$ topology on the configuration space for any $k \geq 2$ (giving the gauge group $\mathrm{Map}(R^o, S^1)$ the $L^2_{k+1,loc}$ topology).

We will be interested in the moduli space of solutions to the Seiberg-Witten equations on the cylinder for which the Chern-Simons function cs has bounded variation along the slices $\{t\} \times Y$. In fact, arguments from from [27] show that for any solution $(A, \Phi)$ with bounded variation of cs along the slices, there are a pair of critical points $[B_\pm, \Psi_\pm]$ for cs for which

$$\lim_{t \mapsto \pm\infty} [(A, \Phi)|_{\{t\}\times Y}] \mapsto \pi^*[B_\pm, \Psi_\pm].$$

This allows us to partition the moduli space into components labeled by pairs of critical sets $\mathcal{C}_1, \mathcal{C}_2$ for cs, $\mathcal{M}(\mathcal{C}_1, \mathcal{C}_2)$, of solutions to the Seiberg-Witten equations over $\mathbb{R} \times Y$ which satisfy

$$\lim_{t \mapsto \pm\infty} [(A, \Phi)|_{\{t\}\times Y}] = \pi^*[B_\pm, \Psi_\pm],$$

for some pair of critical points $[B_-, \Psi_-] \in \mathcal{C}_1$, $[B_+, \Psi_+] \in \mathcal{C}_2$.

The following was observed in [22]:

**Proposition 6.0.6.** *Given a pair of critical sets for* cs*, $\mathcal{C}^\pm$, there is an identification between the moduli space of parameterized gradient flow lines for* cs *from $\mathcal{C}^-$ to $\mathcal{C}^+$ in $\mathcal{B}(W)$ and $\mathcal{M}(\mathcal{C}^+, \mathcal{C}^-)$.*

When $\pi\colon Y \to \Sigma$ is a Seifert fibered space, the cylinder $R^o$ can be given a complex structure as in Section 3. We recall the following basic result (see [34] and [20]) about the Seiberg-Witten equations on any complex surface.

**Theorem 6.0.7.** *Let $X$ be an almost-complex surface with a Riemannian metric, and denote its associated $(1,1)$ form by $\omega$. The* $\mathrm{Spin}_c(4)$ *structures on $X$ naturally correspond to isomorphism classes of Hermitian line bundles $E$ under the correspondence identifying*

$$W^+ \cong (\Lambda^{0,0} \oplus \Lambda^{0,2}) \otimes E \quad \text{and} \quad W^- \cong (\Lambda^{0,1}) \otimes E, \tag{47}$$

*where $\Lambda^{0,0}$ (resp. $\Lambda^{0,2}$) is the $2i$ (resp. $-2i$) eigenbundle of Clifford multiplication by $\omega$. Let $\nabla$ be a connection on $T^*X$ for which*

$$\nabla\omega \equiv 0. \tag{48}$$

*Then the Seiberg-Witten equations for $\nabla$-spinorial connections and spinors in $W^+$ reduce to the following equations for a compatible connection $A$ in the Hermitian line bundle $E$ and sections*

$$\alpha, \beta \in \Omega^{0,0}(X, \mathbb{C}), \Omega^{0,2}(X, \mathbb{C}) :$$

$$2\Lambda F_A - \Lambda F_{K_X} = i(|\alpha|^2 - |\beta|^2) \tag{49}$$

$$2F_A^{0,2} - F_{K_X}^{0,2} = \alpha^* \otimes \beta \tag{50}$$

$$\overline{\partial}_A \alpha + \overline{\partial}_A^* \beta = 0. \tag{51}$$

*In the above,* $F_{K_X}$ *denotes the curvature of the connection which* $\nabla$ *induces on the canonical bundle, and*

$$\Lambda\colon \Omega^2(X, \mathbb{C}) \to \Omega^0(X, \mathbb{C})$$

*denotes the operation*

$$F \mapsto \langle F, \omega \rangle.$$

**Remark 6.0.8.** *For a proof of the above result, see* [20]*. Equations 49 and 50 hold for any complex four-manifold. Equation 51 follows from the identification of the Dirac operator with the rolled-up* $\overline{\partial}$*-operator, which is usually stated for Levi-Civita-spinorial connections on a Kähler manifold. However, the proof uses only the fact that Clifford multiplication by* $\omega$ *commutes with the Dirac operator, which follows from Equation 48.*

For $(R^o, g_{cyl})$ with its complex structure as in Section 3, the $(1,1)$ form is:

$$\omega = dt \wedge \eta + \pi^*(\mu_\Sigma),$$

which is evidently covariantly constant for the connection $\nabla_{cyl}$, so we can apply the above theorem. Moreover, note that the above decomposition of spinors in $W^+$ is compatible with the decomposition of $W$ into $S^\pm$ (Equation 21), in the sense that

$$\pi^*(S^+) \cong \Lambda^{0,0}(E) \quad \text{and} \quad \pi^*(S^-) \cong \Lambda^{0,2}(E). \tag{52}$$

In the compact, Kähler case one uses Theorem 6.0.7, together with an integration-by-parts argument to show a correspondence between solutions to the Seiberg-Witten equations and holomorphic data (see [34]). (This was also the structure of the proof of Theorem 4.) To perform this integration-by-parts argument on the cylinder, we need to know about the behavior at infinity our solution. To this end, we apply the standard exponential decay estimates ([27],[32]) adapted to the context of Seiberg-Witten monopoles:

**Theorem 6.0.9.** *Given any pair of non-degenerate critical manifolds* $\mathcal{C}_-$, $\mathcal{C}_+$ *for* cs *(definition 5.9.3), there is a* $\delta > 0$ *with the property that any flow line between these manifolds corresponds to a solution* $(A, \Phi)$ *to the Seiberg-Witten equations over* $R^o$ *in the* $\mathrm{Spin_c}(4)$ *structure* $\pi^*(W)$ *for which there exist* $[B_\pm, \Psi_\pm] \in \mathcal{C}_\pm$*, with*

$$\lim_{t \mapsto \pm\infty} e^{\pm\delta t}(\|{\nabla_{cyl}}^{(k)}(A - \pi^*(B_\pm))\|_{\{t\}\times Y} + \|\nabla^{(k)}_{\pi^*(B_\pm)}(\Phi - \pi^*(\Psi_\pm))\|_{\{t\}\times Y}) \mapsto 0, \tag{53}$$

*for all* $k \in \mathbb{Z}$*. The norms appearing above are* $\mathcal{C}^0$ *norms on* $Y$*.*

We are now ready to apply the integration-by-parts argument, to obtain a correspondence between gradient flow lines and (certain) vortices over $R^o$.

**Proposition 6.0.10.** *Gradient flow lines between regular critical manifolds of* cs *correspond to configurations*

$$[A,\alpha,\beta]\in(\mathcal{A}(R^o,E)\times\Gamma(R^o,E)\times\Gamma(R^o,\pi^*(K_\Sigma^{-1})\otimes E))/\mathcal{G}(R^o)$$

*which satisfy satisfy the decay condition Equation 53, together with the generalized vortex equations:*

$$\begin{aligned}2\Lambda F_A-\Lambda\pi^*F_{K_\Sigma} &= i(|\alpha|^2-|\beta|^2)\mu_\Sigma\\ F_A^{0,2} &= 0\\ \overline{\partial}_A\alpha=0 \quad \textit{and} & \quad \overline{\partial}_A^*\beta=0\\ \alpha=0 \quad \textit{or} & \quad \beta=0.\end{aligned}$$

**Proof.** First notice that $\nabla_{cyl}=d\oplus d\oplus\pi^*(\widehat{\nabla}_\Sigma)$ with respect to the splitting

$$T^*R^o\cong\mathbb{R}(dt)\oplus\mathbb{R}(\eta)\oplus\pi^*(T^*\Sigma),$$

so that $F_{K_{R^o}}=\pi^*F_\Sigma$; in particular, $F_K^{0,2}=0$. Applying $\overline{\partial}_A$ to Equation 51, and then using Equation 51, we see that

$$0=F_A^{0,2}\alpha+\overline{\partial}_A\overline{\partial}_A^*\beta=|\alpha|^2\beta+\overline{\partial}_A\overline{\partial}_A^*\beta.$$

Taking the inner product with $\beta$, we get

$$0=|\alpha|^2|\beta|^2+\langle\overline{\partial}_A\overline{\partial}_A^*\beta,\beta\rangle. \tag{54}$$

By Stokes' theorem, for any $\gamma\in\Omega^{0,1}(Y,E)$, $\beta\in\Omega^{0,2}(Y,E)$, $T\in\mathbb{R}$,

$$\begin{aligned}\int_{\{T\}\times Y}\gamma\wedge\overline{\beta}-\int_{\{-T\}\times Y}\gamma\wedge\overline{\beta} &= \int_{[-T,T]\times Y}d(\gamma\wedge\overline{\beta})\\ &= \int_{[-T,T]}\overline{\partial}_A\gamma\wedge\overline{\beta}+\int_{[-T,T]}\gamma\wedge\overline{\overline{\partial}_A^*\beta}.\end{aligned}$$

Applying this equation to the integral of Equation 54 over the region $[-T,T]\times Y$, we get

$$\begin{aligned}0 &= \int_{[-T,T]\times Y}|\alpha|^2|\beta|^2+\int_{[-T,T]\times Y}\langle\overline{\partial}_A^*\beta,\overline{\partial}_A^*\beta\rangle\\ &\quad+\int_{\{t\}\times Y}\overline{\partial}_A^*\beta\wedge\overline{\beta}\Big|_{t=-T}^{t=T}\end{aligned}$$

But, applying Theorem 6.0.9 and Theorem 4, we get a convergance (of sections of $W$ over $Y$, in the $\mathcal{C}^\infty$ topology):

$$\begin{aligned}\lim_{t\mapsto\pm\infty}\overline{\partial}_{A_T}^*\beta_T &= \overline{\partial}_{A_{\pm\infty}}^*\beta_{\pm\infty}\\ &= 0.\end{aligned}$$

Thus, the right hand side of Equation 55 is a sum of *a priori* non-negative terms and some terms which go to zero as $T \mapsto \infty$. Hence, the non-negative terms must vanish, proving the theorem. □

**Remark 6.0.11.** *It follows from the above result, along with the compatibility condition Equation 52, that there are no flow lines connecting critical manifolds of opposite sign, i.e.* $\mathcal{C}^{\pm}(\mathbf{e}_1)$ *to* $\mathcal{C}^{\mp}(E_2)$.

**Remark 6.0.12.** *The Serre duality operator on* $R^o$ *induces an identification between the space of flows from* $\mathcal{C}^+(E_1)$ *to* $\mathcal{C}^+(E_2)$ *with the space of flows from* $\mathcal{C}^-(E_1)$ *to* $\mathcal{C}^-(E_2)$*. Similarly for flows between* $\mathcal{C}^{\pm}(E_1)$ *and the reducible locus* $\mathfrak{J}$.

## 7. Extending Holomorphic Data

Our goal in this section is to pass from the decaying vortex data over the cylinder $R^o$ from Proposition 6.0.10, to data over the closed ruled surface $R$.

In light of the duality between positive and negative vortices (Remark 6.0.12), we can assume without loss of generality that all Kähler vortices over curves are positive, i.e. vortices for which the $\alpha$-component does not vanish identically (hence $\beta \equiv 0$). We return to the reducible case in Section 10.

Let $E^o \to R^o$ be a (Hermitian) orbifold line bundle over the cylinder $\mathbb{R} \times Y$, and $E_{\pm}$ be a pair of Hermitian orbifold bundles over $\Sigma$, containing vortices $(A_{\pm}, \alpha_{\pm})$.

**Definition 7.0.13.** *A pair*

$$(A, \alpha) \in \mathcal{A}(E^o) \times \Gamma(R^o, E^o)$$

*is said to* connect the vortex $(A_-, \alpha_-)$ to $(A_+, \alpha_+)$ *if there are isomorphisms*

$$j_{\pm} \colon \pi^*(E_{\pm})|_{R^o} \to E^o$$

*and a real number* $\delta > 0$*, so that for all* $k \in \mathbb{Z}$ *with* $k \geq 0$,

(55)
$$\lim_{t \mapsto \pm\infty} e^{\pm\delta t} \left( \|{\nabla_{cyl}}^{(k)}(j_{\pm}^*(A) - \pi^*(A_{\pm}))\|_{\{t\}\times Y} + \|\nabla^{(k)}_{\pi^*(A_{\pm})}(j_{\pm}^*(\alpha) - \pi^*(\alpha_{\pm}))\|_{\{t\}\times Y} \right) \mapsto 0.$$

**Definition 7.0.14.** *A pair* $(A, \alpha)$ *is called a* generalized vortex *if it satisfies:*

$$\begin{aligned} 2\Lambda F_A - \Lambda\pi^*(F_{K_\Sigma}) &= i|\alpha|^2 \\ F_A^{0,2} &= 0 \\ \overline{\partial}_A \alpha &= 0. \end{aligned}$$

The moduli space of generalized vortices connecting $E_-$ to $E_+$ is obtained by dividing out the space of generalized vortex pairs which connect vortices in $E_-$ to vortices in $E_+$ by the action of the gauge group of Hermitian transformations whose derivatives exponentially decay to zero. When

$$0 \leq \deg(E_\pm) < \frac{\deg(K_\Sigma)}{2},$$

results of Section 6, allow us to view this space as the moduli space of flows $\mathcal{M}(\mathcal{C}^+(E_1), \mathcal{C}^+(E_2))$

Let $E$ be a Hermitian vector bundle over the closed ruled surface $R$, with

$$E|_{\Sigma_-} = E_- \quad \text{and} \quad E|_{\Sigma_+} = E_+.$$

Note that the isomorphism type of $E$ is uniquely determined by $E_-$ and $E_+$.

**Definition 7.0.15.** *A pair $(A, \alpha) \in \mathcal{A}(E) \times \Gamma(R, E)$ is called a* holomorphic pair connecting $E_-$ to $E_+$ *if it satisfies the equations*

$$\begin{aligned} F_A^{0,2} &= 0, \\ \overline{\partial}_A \alpha &= 0. \end{aligned}$$

*The* moduli space of effective divisors connecting $E_-$ to $E_+$*, denoted $\mathcal{D}(E_-, E_+)$, is obtained by dividing out the space of such holomorphic pairs by the natural action of the complex gauge group* $\mathrm{Map}(R, \mathbb{C}^*)$*.*

**Remark 7.0.16.** *Strictly speaking, the above definition endows the space of divisors with a topology, but not a deformation theory. The moduli space of divisors can be given a deformation theory, locally modeling a neighborhood of the divisor $[A, \alpha]$ on the zeros of a map*

$$H^0(\alpha^{-1}(0), E|_{\alpha^{-1}(0)}) \to H^1(\alpha^{-1}(0), E|_{\alpha^{-1}(0)}).$$

*We will return to this point in Section 9.*

**Definition 7.0.17.** *Given a subset $C \subset R$, a divisor $[A, \alpha] \in \mathcal{D}(E_-, E_+)$ is said to* contain $C$ *if it is represented by a pair $(A, \alpha)$ for which $\alpha|_C \equiv 0$.*

**Theorem 7.0.18.** *There is a natural identification between the moduli space of generalized vortices connecting $E_-$ to $E_+$, with the open subset of divisors in $\mathcal{D}(E_-, E_+)$ which do not contain $\Sigma_-$ or $\Sigma_+$.*

We break this theorem up into two parts. First (Theorem 7.0.19), we show that a generalized vortex which connects the vortex $(A_-, \alpha_-)$ to $(A_+, \alpha_+)$ completes canonically to a holomorphic pair interpolating between these two vortices. The methods are essentially in [15] and [3]. In fact, since the structure group in question here is $U(1)$ rather than $SU(2)$ bundles, the arguments are rather simpler; but we include them for completeness. Then, in Section 8, we will describe how to invert this operation.

The rest of this section will focus on the following analogue of Guo's theorem [15]:

**Theorem 7.0.19.** *Let $(A, \alpha)$ be a generalized vortex in $R^o$ connecting $(A_-, \alpha_-)$ to $(A_+, \alpha_+)$. Then $(A, \alpha)$ admits a canonical extension to a holomorphic pair in $R$ which connects $(A_-, \alpha_-)$ to $(A_+, \alpha_+)$.*

First, we concentrate on the easier problem of canonically extending the section $\alpha$.

**Lemma 7.0.20.** *Given a generalized vortex $(A, \alpha)$ in the line bundle $E^o$ over $R^o$ which connects vortices $(A_\pm, \alpha_\pm)$ in $E_\pm$, there is a line bundle $\widehat{E}$ over $R$, a continuous section $\widehat{\alpha} \in \Gamma(R, \widehat{E})$, and an isomorphism*

$$k \colon \widehat{E}|_{R^o} \to E$$

*with*

$$k^*(\alpha) = \widehat{\alpha}|_{R^o}.$$

*The isomorphism class of the triple $(\widehat{E}, \widehat{\alpha}, k)$ is uniquely determined by the isomorphism class of the generalized vortex $(A, \alpha)$.*

**Proof.** Let $\widehat{E}$ be the bundle obtained by gluing $\pi^*(E_\pm)$ to $E^o$ using the isomorphisms $j_\pm$. The section is defined to be $\alpha$ away from $\Sigma_\pm$ and it is $\alpha_\pm$ over $\Sigma_\pm$. Continuity of $\widehat{\alpha}$ is an easy consequence of Equation 55.

To verify uniqueness of the construction, let $(A^1, \alpha^1)$, $(A^2, \alpha^2)$ be a pair of generalized vortices in line bundles $E^1, E^2$ which connect vortices $(A^i_\pm, \alpha^i_\pm)$ $i = 1, 2$, and let $(\widehat{E}^i, \widehat{\alpha}^i, k^i)$ $i = 1, 2$ be the respective extensions. Suppose moreover that $E^1 \cong E^2$ via an isomorphism $\ell \colon E^1 \to E^2$. Then there is a continuous isomorphism

$$\widehat{\ell} \colon \widehat{E}^1 \to \widehat{E}^2$$

which extends the isomorphism

$$(k^2)^{-1} \circ \ell \circ k^1 \colon \widehat{E}^1|_{R^o} \to \widehat{E}^2|_{R^o},$$

and is compatible with the sections, in the sense that

$$\widehat{\ell^*}(\widehat{\alpha}_2) = \widehat{\alpha}_1.$$

It suffices to prove that the bundle map

$$\pi^*(E^1_+)|_{R^o} \to \pi^*(E^2_+)|_{R^o}$$

defined by

$$\lambda = j_2^{-1} \circ \ell \circ j_1$$

extends continuously over $\Sigma_+ \subset R$. Viewing $\lambda$ as a section of the Hermitian line bundle with connection $\pi^*((E^1_+)^* \otimes E^2_+)$ (the connections induced by $A^i_+$ $i = 1, 2$), Equation 55 implies that

$$\lim_{t \mapsto +\infty} e^{t\delta} \|\nabla\lambda\|_{\{t\}\times Y} \mapsto 0.$$

This latter condition ensures that $\lambda$ extends as a continuous section over $R-\Sigma_-$. Repeating this construction over $\Sigma_-$, we get the required isomorphism $\widehat{\ell}$. The compatibility of the sections follows from the fact that $\widehat{\ell}^*(\widehat{\alpha}_1)$ and $\widehat{\alpha}_2$ are continuous sections which agree on the dense subset $R^o \subset R$. □

We turn our attention to extending the holomorphic structure from from $E^o$ to $\widehat{E}$ in such a way as to make $\widehat{\alpha}$ holomorphic. Once again, this holomorphic structure is canonically determined by the original generalized vortex $(A,\alpha)$. Suppose that $(\widehat{E}_1,\widehat{\alpha}_1,k_1)$, $(\widehat{E}_2,\widehat{\alpha}_2,k_2)$ are extensions of generalized vortices $(A_1,\alpha_1)$, $(A_2,\alpha_2)$ in bundles $E_1,E_2$, which are identified by a bundle isomorphism

$$\widehat{\ell}\colon \widehat{E}_1 \to \widehat{E}_2.$$

Then, $\widehat{\alpha}_1 = \widehat{\ell}^*(\widehat{\alpha}_2)$ is holomorphic with respect to both connections $\widehat{A}_1$ and $\widehat{\ell}^*(\widehat{A}_2)$. Thus,

$$(\widehat{A}_1 - \widehat{\ell}^*(\widehat{A}_2))^{0,1} \otimes \widehat{\alpha}_1 = 0,$$

forcing $\widehat{A}_1 - \widehat{\ell}^*(\widehat{A}_2) = 0$ on the dense open set where $\widehat{\alpha}_1 \neq 0$, so that $\widehat{A}_1 = \widehat{\ell}^*(\widehat{A}_2)$ everywhere.

For notational simplicity, we will discuss the extension from $R^o \cong Y \times_{S^1} (\mathbb{C}-0)$ across its zero-section, $\Sigma_-$, to $Y \times_{S^1} \mathbb{C} \cong R - \Sigma_+$. (Though the discussion applies to both ends of $R^o$, giving the requisite extension to all of $R$.)

Let $g_\circ$ be the metric on on $\mathbb{R}\times Y$ given by the formula

$$g_\circ = r^2(dt^2+\eta^2) + g_\Sigma,$$

where $r=e^t$. Let $\nabla_\circ$ be the connection on $T^*R^o$ characterized by the properties that

$$\nabla_\circ dt = -dt\otimes dt + \eta\otimes\eta,$$
$$\nabla_\circ \eta = -dt\otimes\eta - \eta\otimes dt,$$
$$\nabla_\circ \pi^*\theta = \pi^*(\nabla_\Sigma\theta).$$

Then, we have the following:

**Lemma 7.0.21.** *Under the natural holomorphic identification*

$$R^o \cong Y\times_{S^1}(\mathbb{C}-0),$$

*the metric $g_\circ$ and the compatible connection $\nabla_\circ$ extend as an orbifold metric and compatible connection over $T^*(Y\times_{S^1}\mathbb{C})$, with orbifold singularities only at the orbifold points of $Y\times_{S^1}\mathbb{C}$.*

**Proof.** The connection on $\pi\colon Y\to\Sigma$ induces an splitting

$$T_yY \cong T_{\pi(y)}\Sigma \oplus T_0\mathbb{C}.$$

The metric $g_\circ$ and connection $\nabla_\circ$ are the metric and connection induced by this splitting. □

**Remark 7.0.22.** *There is a considerable amount of leeway in the choices here – there are many metrics and compatible connections which extend over $\Sigma$; we have made choices which require a minimum of calculation.*

Given a $T \in \mathbb{R}$, let $\mathcal{N}_T$ denote subset of $R^o$ given by

$$\mathcal{N}_T = \{(t,y)\big|t < T\}).$$

**Lemma 7.0.23.** *Given $k \in \mathbb{N}$, $T \in \mathbb{R}$, there is a constant $C$ such that for any $a \in \Omega^1(\mathcal{N}_T, \mathbb{C})$, we have*

$$|\nabla_\circ^{(k)} a|_\circ \leq C r^{-(k+1)} (\sum_{j=0}^{k} |\nabla_{cyl}^{(j)} a|_{cyl}).$$

*Here, $||_\circ$ and $||_{cyl}$ indicate the pointwise metrics $g_\circ$ and $g_{cyl}$ induce on the spaces of tensors, and*

$$\nabla_\circ^{(j)}, \nabla_{cyl}^{(j)} \colon \Omega^1(\mathcal{N}_T, \mathbb{C}) \to \Gamma(\mathcal{N}_T; (\Lambda^1)^{\otimes(j+1)})$$

*are the operators induced by iterating the covariant derivatives.*

**Proof.** It is clear from the formulas for $g_{cyl}$ and $g_\circ$ that there exists a constant $C_1$ (depending on $j$) such that for any $a \in \Gamma(\mathcal{N}_T; (\Lambda^1)^{\otimes(j)})$, we have

$$|a|_\circ \leq C_1 r^{-j} |a|_{cyl}. \tag{56}$$

We have that $\nabla_\circ - \nabla_{cyl}$ is a zeroth order operator

$$\Gamma(\mathcal{N}_T; (\Lambda^1)^{\otimes(j)}) \to \Gamma(\mathcal{N}_T; (\Lambda^1)^{\otimes(j+1)}).$$

When $j = 1$, the operator is given by

$$a \mapsto \langle a, dt \rangle_{cyl} (-dt \otimes dt + \eta \otimes \eta) + \langle a, \eta \rangle_{cyl} (-dt \otimes \eta - \eta \otimes dt);$$

so for any $j$ there is a constant $C_2$ such that for any $a \in \Gamma(\mathcal{N}_T; (\Lambda^1)^{\otimes(j)})$,

$$|(\nabla_\circ - \nabla_{cyl}) a|_{cyl} \leq C_2 |a|_{cyl}.$$

Since $\eta$ and $dt$ are $\nabla_{cyl}$-constant forms, we see that there is a constant $C_3$ such that

$$|[\nabla_{cyl}, \nabla_\circ - \nabla_{cyl}] a|_{cyl} \leq C_3 |\nabla_{cyl} a|.$$

Thus, we have

$$\begin{aligned} |\nabla_\circ^{(k)} a|_{cyl} &= |(\nabla_{cyl} + (\nabla_\circ - \nabla_{cyl}))^{(k)} a|_{cyl} \\ &\leq C_4 \sum_{j=0}^{k} |\nabla_{cyl}^{(j)} a|_{cyl} \end{aligned}$$

The formula now follows from this together with Equation 56. □

Suppose for simplicity that $\Sigma$ has no marked points. The bundle structure of $\pi\colon Y \to \Sigma$ provides us with an open cover $\{U_i\}$ of $\Sigma$, a corresponding cover $\{\mathcal{U}_i\}$ of $Y \times_{S^1} \mathbb{C}$, and complex trivializations

$$\tau_i\colon U_i \times \mathbb{C} \to \mathcal{U}_i.$$

Consider now the $m$-fold branched cover

$$\phi_m\colon \mathcal{U}_i^m \to \mathcal{U}_i,$$

where $\mathcal{U}_i^m$ is the complex space $U_i \times \mathbb{C}$, and $\phi_m$ is given by the formula

$$\phi_m(w,z) = (w, z^m).$$

The space $U_i \times S^1 \times \mathbb{R}$ admits two canonical metrics, a cylindrical metric given by the formula

$$g_{cyl} = dt^2 + (\frac{1}{m}\phi_m^*\eta)^2 + g_\Sigma,$$

and a disk metric

$$g_\circ = r^2(dt^2 + (\frac{1}{m}\phi_m^*\eta)^2) + g_\Sigma,$$

which completes over $\mathcal{U}_i^m$ under the natural holomorphic identification

$$U_i \times S^1 \times \mathbb{R} \cong U_i \times (\mathbb{C} - 0) \subset \mathcal{U}_i^m.$$

Since $\phi_m^*(r) = r^m$, if $a \in \Omega^1(\mathcal{U}_i, \mathbb{C})$, satisfies $a \in L^2_{k,\delta}$ with respect to the cylindrical metric for some real $\delta > 0$ and integral $k \geq 0$, then $\phi_m^* a \in L^2_{k,m\delta}$. Thus, given any $k$, $a \in \Omega^1(R^o, \mathbb{C})$, with

$$a \in L^2_{k,\delta},$$

Lemma 7.0.23 guarantees that we can find $m$ large enough that for each open set $\mathcal{U}_i^m$,

$$\phi_m^*(a) \in L^2_{k,loc}(\Omega^1(\mathcal{U}_i^m, \mathbb{C})).$$

Moreover, when $k \geq 2$, $\phi_m^*(a)$ continuously extends over $U_i \times \{0\} \subset \mathcal{U}_i^m$ by vanishing on this subset, thanks the inclusion $L^2_{k,loc} \subset \mathcal{C}^0$.

Recall the following version of the $\overline{\partial}$-Poincaré lemma:

**Theorem 7.0.24.** *Let $X$ be a complex manifold. If $a \in L^2_k(\Omega^1(X, \mathbb{C}))$ for $k \geq \dim_{\mathbb{C}} X$ satisfies $\overline{\partial} a^{0,1} = 0$, then there is an open cover $\{V_j\}$ of $X$ and a collection of invertible functions $\{g_j \in L^2_{k+1}(V_j)\}$ with the property that*

$$\overline{\partial} g_j = -a^{0,1} g_j. \tag{57}$$

*If $\Sigma \subset X$ is a complex submanifold with the property that $a|_\Sigma \equiv 0$, then we can arrange that the $g_i$ satisfy*

$$g_j|_\Sigma \equiv 1. \tag{58}$$

**Proof.** The existence of invertible functions $g_j$ satisfying Equation 57 is standard (see for example [6]). Since $a|_\Sigma = 0$, the $g_j$ are holomorphic over $\Sigma$, so we can find holomorphic functions $h_j$ defined over the $V_j$ which hit $\Sigma$ such that $h_j = g_j^{-1}$ over $\Sigma \cap V_j$. The functions $h_j g_j$ now satisfy both Equations 57 and 58. □

Given the generalized vortex $(A, \alpha)$ from Theorem 7.0.19, we have that the difference one-form $a = j_-^*(A) - \pi^*(A_-)$ lies in the Sobolev space $L^2_{k,\delta}(R^o)$ for any $k \in \mathbb{Z}$. Combining Lemma 7.0.23 (and the discussion following it), together with Theorem 7.0.24, we see that there is a natural number $m$, an open cover $\mathcal{U}_i$ of a neighborhood of $\Sigma$ (obtained by shrinking the original $\mathcal{U}_i$, and possibly reindexing), a collection of holomorphic maps

$$\phi_m \colon \mathcal{U}_i^m \to \mathcal{U}_i,$$

which are $m$-fold branched covers branching over $\mathcal{U}_i \cap \Sigma$, and invertible functions

$$g_i \colon \mathcal{U}_i^m \to \mathbb{C}$$

satisfying

$$\overline{\partial} g_i = -a^{0,1} g_i \tag{59}$$

and

$$g_i|_{\mathcal{U}_i^m \cap \Sigma} \equiv 1. \tag{60}$$

The $\mathbb{Z}/m\mathbb{Z}$ action on the $\mathcal{U}_i^m$ lifts naturally to an action on $\phi_m^*(E)$; indeed, as the connections $A$ and $A_-$ pull back from $\mathcal{U}_i \subset R$, both connections, and their difference $a$, are invariant under the $\mathbb{Z}/m\mathbb{Z}$ action. So, replacing $g_i$ by its average

$$\frac{1}{m} \sum_{\gamma \in \mathbb{Z}/m\mathbb{Z}} \gamma^* g_i,$$

we get a new collection of smooth functions, still denoted $g_i$, satisfying Equations 59 and 60, but which are now $\mathbb{Z}/m\mathbb{Z}$-invariant. By Equation 60, we see that, after perhaps shrinking the $\mathcal{U}_i$, we can arrange that these functions $g_i$ are invertible over $\mathcal{U}_i$. Since the $\mathbb{Z}/m\mathbb{Z}$-invariant continuous functions on $\mathcal{U}_i^m$ correspond to the continuous functions on the quotient $\mathcal{U}_i$, we can think of the $\{(\mathcal{U}_i^m, g_i)\}$ as defining a Čech cochain $\{(\mathcal{U}_i, \widehat{g}_i)\}$ defined over a neighborhood of $\Sigma \subset Y \times_{S^1} \mathbb{C} \subset R$ with values in the *continuous* invertible functions. Since the $\mathbb{Z}/m\mathbb{Z}$ action on $\mathcal{U}_i^m$ is free away from $\Sigma \cap \mathcal{U}_i^m$, we have that the $\widehat{g}_i$ are smooth away from $\mathcal{U}_i \cap \Sigma$. Equation 59 translates into the equation over $\mathcal{U}_i - \mathcal{U}_i \cap \Sigma$

$$\overline{\partial} \widehat{g}_i = -a \widehat{g}_i \tag{61}$$

and Equation 60 translates into

$$\widehat{g}_i|_{\mathcal{U}_i \cap \Sigma} \equiv 1. \tag{62}$$

**Lemma 7.0.25.** *The coboundary of the cochain* $\{(\mathcal{U}_i, \widehat{g}_i)\}$ *takes values in the sheaf of invertible holomorphic functions, i.e.*

$$\delta\{(\mathcal{U}_i, \widehat{g}_i)\} \in Z^1(\{\mathcal{U}_i\}, \mathcal{O}^*_R).$$

**Proof.** The functions $\widehat{g}_{i,j}$

$$\widehat{g}_{i,j} = \widehat{g}_i^{-1}\widehat{g}_j$$

are holomorphic over all of $\mathcal{U}_i \cap \mathcal{U}_j$, because they are continuous, and they satisfy

$$\overline{\partial}(\widehat{g}_i\widehat{g}_j^{-1}) = 0$$

over $\mathcal{U}_i \cap \mathcal{U}_j \cap (R - \Sigma)$, by Equation 61. So they are holomorphic over all of $\mathcal{U}_i \cap \mathcal{U}_j$ by the regularity of the $\overline{\partial}$ operator. □

Thus, the cocycle $\delta\{(\mathcal{U}_i, \widehat{g}_i)\}$ naturally induces a holomorphic line bundle $F$ over $Y \times_{S^1} \mathbb{C}$, and the cochain $\{(\mathcal{U}_i, \widehat{g}_i)\}$ glues together to induce a continuous trivialization $\widehat{g}$ of $F$.

**Lemma 7.0.26.** *The topological isomorphism*

$$j_- \otimes \widehat{g} \colon \pi^*(E_-) \otimes F|_{R^o} \to E^o$$

*induces an isomorphism of holomorphic bundles*

$$(E, \overline{\partial}_A) \cong (\pi^*(E_-) \otimes F, \overline{\partial}_{\pi^*(A_-)\otimes F})|_{R^o}.$$

**Proof.** This is an easy consequence of Equation 61:

$$\begin{aligned}
\widehat{g}_j \circ \overline{\partial}_{A_-} \circ \widehat{g}_j^{-1} - \overline{\partial}_A &= a^{0,1} + \widehat{g}_j(\overline{\partial}\widehat{g}_j^{-1}) \\
&= a^{0,1} - (\overline{\partial}\widehat{g}_j)\widehat{g}_j^{-1} \\
&= 0.
\end{aligned}$$

□

The pull-back of $\widehat{\alpha}$ to $\pi^*(E_-)\otimes F$ via the isomorphism induced by $\widehat{g}$ is a continuous section of $\pi^*(E_-)\otimes F$ which is holomorphic away from $\Sigma_-$. The regularity of $\overline{\partial}_{A_-\otimes F}$, in fact, ensures that it is a holomorphic section over all of $Y \times_{S^1} \mathbb{C}$.

The above discussion easily generalizes in the presence of orbifold points on $\Sigma$. The key point is that the construction of the $\widehat{g}_i$ can be made $\mathbb{Z}/\alpha_i\mathbb{Z}$-equivariantly, by averaging over this group.

This concludes the proof of Theorem 7.0.19. Our results can be rephrased as follows:

**Theorem 7.0.27.** *Let $E_1$, $E_2$ be a pair of orbifold line bundles over $\Sigma$ with*

$$0 \leq \deg(E_i) < \frac{\deg(K_\Sigma)}{2}.$$

*The canonical extension procedure used in Theorem 7.0.19 induces a continuous map*

$$\mathcal{M}(\mathcal{C}^+(E_1), \mathcal{C}^+(E_2)) \to \mathcal{D}(E_1, E_2).$$

**Proof.** The fact that the map induces a well-defined map on the level of gauge equivalence classes follows immediately from the uniqueness statement in Theorem 7.0.19.

Continuity follows from a corresponding continuity statement in the solution to the $\overline{\partial}$-problem from Theorem 7.0.24. □

## 8. Vortices over Cylinders

In Section 7, we showed how a generalized vortex could be canonically completed to give a divisor in the ruled surface $R$. Presently, we turn our attention to the converse problem:

**Theorem 8.0.28.** *Let $E_-, E_+$ be a pair of orbifold bundles over $\Sigma$, with*

$$0 \leq \deg(E_\pm) < \frac{\deg(K_\Sigma)}{2}.$$

*Suppose that $(A, \alpha)$ is a holomorphic pair connecting $E_-$ to $E_+$ such that $\alpha|_{\Sigma_-} \not\equiv 0$ and $\alpha|_{\Sigma_+} \not\equiv 0$, then there is a complex gauge transformation, $g$, defined over $R^o$ which carries the restriction $(A, \alpha)|_{R^o}$ to a generalized vortex connecting $E_-$ to $E_+$. The complex gauge transformation $g$ is unique up to unitary gauge tranformations defined on the cylinder.*

The first step in the proof is to show that a holomorphic pair as in the thoerem gives rise to a holomorphic pair on the cylinder which connects the corresponding vortices. Indeed, we can arrange, after a complex gauge transformation, that the restriction of $(A, \alpha)|_{\Sigma_\pm}$ satisfies the Kähler vortex equations over $\Sigma_\pm$, according to Theorem 5. Furthermore, any two such complex gauge transformation differ over $\Sigma_\pm$ by a unitary gauge transformation. The first step is completed by the following lemma.

**Lemma 8.0.29.** *Let $(A, \alpha)$ be a pair in $E \to R$, whose restrictions $(A_\pm, \alpha_\pm)$ to $\Sigma_\pm$ satisfy the Kähler vortex equations. Then there is a real number $\gamma > 0$ such that the pair $(A, \alpha)|_{R^o}$ connects $(A_-, \alpha_-)$ to $(A_+, \alpha_+)$ with exponent one.*

**Proof.** This follows easily from Taylor's theorem. □

Given a holomorphic pair $(A,\alpha)$ over $R^o$ which arises as above, we wish to find a real-valued function $u$ (defined over $R^o$) so that $e^u(A,\alpha)$ satisfies Equation 56. Watching how its terms transform under the action of $e^u$, we rewrite this equation as a second-order, Kazdan-Warner type equation for $u$:

$$-2i\Delta_{cyl} u + 2\Lambda F_A - \Lambda\pi^*(F_{K_\Sigma}) = ie^{2u}|\alpha|^2,$$

where

$$\Delta_{cyl} = 2i\Lambda\overline{\partial}\partial,$$

and $\Lambda$ denotes the operation of inner product with the canonical $(1,1)$ form for the cylindrical metric on $R^o$. Theorem 8.0.28 follows immediately from the next proposition concerning the existence and uniqueness of such a $u$ after making the subsitutions $h = |\alpha|^2$ and $k = -2i\Lambda F_A + i\Lambda\pi^*(F_{K_\Sigma}) - h$. The hypotheses are satisfied thanks to the exponential convergence in Lemma 8.0.29.

**Proposition 8.0.30.** *Suppose that $Y$ is oriented so that $\deg(Y) < 0$. Let*

$$h\colon R^o \to \mathbb{R}$$

*be a smooth, bounded, non-negative function which converges exponentially as $t \mapsto \pm\infty$ to a pair of functions*

$$h_\pm\colon Y \to \mathbb{R}$$

*with $h_- \not\equiv 0$, and let $\mathcal{N}$ be the operator defined by*

$$\mathcal{N}(u) = \Delta_{cyl}(u) + h(e^u - 1).$$

*Then, given any smooth function $k \in \mathcal{C}^\infty_\gamma$, with $\gamma > 0$ there is an $\epsilon > 0$ and a unique function $u \in \mathcal{C}^\infty_\epsilon(R^o)$ satisfying*

$$\mathcal{N}(u) = k. \tag{63}$$

**Remark 8.0.31.** *Recall that a smooth function $u\colon R^o \to \mathbb{R}$ is in $\mathcal{C}^\infty_\epsilon$ if for each $k$, the function $(\nabla^{(k)} f)e^{\epsilon|t|}$ is bounded.*

**Remark 8.0.32.** *If $Y$ is given the other orientation, Equation 63 can be solved when $h_+ \not\equiv 0$. Thus, for the purposes of Theorem 8.0.28, the orientation is irrelevant.*

**Remark 8.0.33.** *The above proposition should be compared with* [18]*, where a similar operator is solved over a closed manifold.*

The operator $\Delta_{cyl}$ can be understood quite concretely, comparing it with the Laplacian for the Kähler metric on $\mathbb{R}\times Y$, see for example [15] or [23], where the following lemma is proved:

**Lemma 8.0.34.** *The operator $\Delta_{cyl}$ can be written*

$$\Delta_{cyl} = -e^{-2\xi t}\frac{\partial}{\partial t}e^{2\xi t}\frac{\partial}{\partial t} + \Delta_Y, \tag{64}$$

*where $\Delta_Y$ is the ordinary Laplacian on the three-manifold $Y$.*

Here, $\xi$ is the constant defined in Equation 16.

First, we show uniqueness for the scale $u$ which carries the pair $(A, \alpha)$ to a generalized vortex, using a lemma which will be useful in the existence argument as well.

**Lemma 8.0.35.** *Let $u$ and $v$ be smooth functions vanishing at infinity, with the property that*

$$\mathcal{N}(u) \leq \mathcal{N}(v)$$

*pointwise, then*

$$u \leq v$$

*pointwise, as well.*

**Proof.** This follows easily from the maximum principle, together with the monotonicity of the exponential map. We wish to show that $u \leq v$ globally. If not, then there must be a point $x$ and a real number $\delta > 0$ such that $u(x) - v(x) > \delta$. Consider the set

$$\Omega_\delta = \{x \in R^o \big| u - v \geq \delta\}.$$

Since $u - v$ vanishes at infinity, each $\Omega_\delta$ is a bounded set. Moreover, from its definition and the choice of $\delta$, we see that the function $u - v$ cannot achieve its maximum on the boundary of $\Omega_\delta$. The assumption that $\mathcal{N}(u) \leq \mathcal{N}(v)$ forces $\Delta_{cyl}(u - v) \leq 0$ on $\Omega_\delta$, hence by the maximum principle for the Laplace-type operator $\Delta_{cyl}$, the function $u - v$ must achieve its maximum on the boundary, a contradiction. □

Of course, this lemma guarantees that for any function $k$, the equation $\mathcal{N}(u) = k$ can have at most one solution which vanishes at infinity.

The rest of this section is devoted to the proof of the existence statement in Proposition 8.0.30, using the continuity method. That is, we will show that for each $s \in [0, 1]$, there is a solution $u_s$ to the equation

$$\mathcal{N}(u_s) = sk, \tag{65}$$

by showing the set of $s$ for which a solution $u_s$ exists is closed and open. Since the set is non-empty (at $s = 0$, the function $u_0 \equiv 0$ is clearly a solution), we conclude that our desired solution at $s = 1$ exists.

We find it convenient to work in the Hilbert space $L^2_{3,\epsilon}$, a weighted Sobolev space. This space is defined as follows. Let

$$\tau \colon \mathbb{R} \to \mathbb{R}$$

be a positive, smooth function with

$$\tau(t) = |t| \quad \text{when} \quad |t| > 1.$$

Now, $L^2_{\ell,\epsilon}$ is by definition the completion of $C^\infty_c(R^o)$ under the norm

$$\|u\|_{L^2_{\ell,\epsilon}} = \int_{R^o} (|u|^2 + |\nabla u|^2 + |\nabla^{(2)} u|^2 + ... + |\nabla^{(\ell)} u|^2) e^{\epsilon\tau} \mu_{R^o}.$$

**Lemma 8.0.36.** *For all $\epsilon > 0$, the operator $\mathcal{N}$ extends to a smooth operator*

$$\mathcal{N}\colon L^2_{3,\epsilon} \to L^2_{1,\epsilon},$$

*with derivative*

$$D\mathcal{N}_u(v) = (\Delta_{cyl} + e^u h)(v).$$

**Proof.** Let $g_1$ be the analytic function on $\mathbb{R}$ defined by the property that

$$u g_1(u) = e^u - 1.$$

We notice first that any function $u \in L^2_{3,\epsilon}$ must be continuous, and must vanish at infinity, by the continuity of the embedding

$$L^2_{3,\epsilon}(\Omega) \subset \mathcal{C}^0(\Omega),$$

valid for any bounded domain $\Omega \subset R^o$. In fact, the geometry of the cylinder gives us bounds.

$$\|u\|_{\mathcal{C}^0_{\epsilon/2}} \leq C\|u\|_{L^2_{3,\epsilon}},$$

where the $\mathcal{C}^\ell_\delta$-norm is the weighted norm defined by

$$\|f\|_{\mathcal{C}^\ell_\delta} = \sum_{i=1}^{\ell} \sup_{R^o} e^{\delta|t|} |\nabla^{(i)} f|.$$

In particular, $g_1(u)$ must be globally bounded:

$$g_1(u) \leq \sup_{[-M,M]} g_1, \tag{66}$$

where

$$M = C\|u\|_{L^2_{3,\epsilon}}.$$

Thus,

$$\|h(e^u - 1)\|_{L^2_{0,\epsilon}} = \|u h g_1(u)\|_{L^2_{0,\epsilon}} \leq \|u\|_{L^2_{0,\epsilon}} \|h g_1(u)\|_{\mathcal{C}^0},$$

and similarly

$$\|\nabla h(e^u - 1)\|_{L^2_{0,\epsilon}} \leq \|u\|_{L^2_{0,\epsilon}} \|(\nabla h) g_1(u)\|_{\mathcal{C}^0} + \|\nabla u\|_{L^2_{0,\epsilon}} \|h e^u\|_{\mathcal{C}^0},$$

so that

$$u \mapsto h(e^u - 1)$$

maps $L^2_{3,\epsilon}$ to $L^2_{1,\epsilon}$. Continuity follows at once from Inequality 66, together with the above $L^2_{0,\epsilon}$ estimates on $h(e^u - 1)$ and its first derivative.

To verify the formula for the derivative, take one more term in the Taylor expansion of $e^u$, and use the same kinds of estimates. □

Moreover, ellipticity guarantees that any solution lying in this space is smooth.

**Lemma 8.0.37.** *Suppose $k$ is a $C^\infty$ function, and $u \in L^2_{3,\epsilon}$ satisfies*

$$\mathcal{N}(u) = k,$$

*then $u$ is also $C^\infty$.*

**Proof.** This is a standard bootstrapping argument using the usual elliptic estimates on $\Delta_{cyl}$. We will see arguments of this nature in the proof of Lemma 8.0.39.

□

The maximum principle from Lemma 8.0.35, combined with the following lemma, ensure bounds of the form

$$u_- \le u_s \le u_+,$$

for a pair of exponentially decaying functions $u_-$ and $u_+$.

**Lemma 8.0.38.** *Let $k \in \mathcal{C}^0_\gamma$, then there is some $\epsilon_0$ such that for all $0 < \epsilon < \epsilon_0$, there are functions $u_\pm \in \mathcal{C}^\infty_\epsilon$ with the property that for any $s \in [0, 1]$*

$$\mathcal{N}(u_-) \le sk \le \mathcal{N}(u_+).$$

The proof of the above lemma is deferred until the end of the section.

Fixing an $\epsilon > 0$ and functions $u_\pm$ as in the previous lemma, closedness of the set of $s \in [0, 1]$ satisfying 65 now follows from the standard machinery of Sobolev spaces, as we see in the following lemma.

**Lemma 8.0.39.** *Let $u_i$ be a sequence of functions in $L^2_{3,\epsilon}$ with $u_- \le u_i \le u_+$, and $u_\pm \in \mathcal{C}^\infty_\epsilon$, with*

$$\lim_{i \mapsto \infty} \mathcal{N}(u_i) = k;$$

*then there is a function $u \in L^2_{3,\epsilon}$, with*

$$\mathcal{N}(u) = k.$$

**Proof.** Bootstrapping and compact inclusions of Sobolev spaces, along with the *a priori* estimates given by $u_- \le u_i \le u_+$, together with the equations allow us to extract an $L^2_{3,loc}$-convergent subsequence, which by continuity must solve $\mathcal{N}(u) = k$. Elliptic estimates applied to the equation for $u$ then allow us to conclude that $u$ has the appropriate decay properties. We give the details.

The Friedrich lemma states that for any pair $\Omega \subset \Omega' \subset R^o$ of precompact open sets and positive integer $\ell$ there is a constant $C > 0$ so that

$$\|u\|_{L^2_\ell(\Omega)} \le C(\|u\|_{L^2(\Omega')} + \|\Delta_{cyl} u\|_{L^2_{\ell-2}(\Omega')}),$$

Moreover, we also have global estimates

$$\|u\|_{L^2_{\ell,\epsilon}} \leq C(\|u\|_{L^2_\epsilon} + \|\Delta_{cyl} u\|_{L^2_{\ell-2,\epsilon}}),$$

thanks to the $\mathbb{R}$-invariance of $\Delta_{cyl}$, and the geometry of the cylinder.

To be able to use compactness, we will first derive local uniform $L^2_4$ bounds on the sequence $u_i$. Fix a pair of pre-compact regions $\Omega \subset \Omega' \subset R^o$. The equation for $u_i$ together with the elliptic estimate gives

$$\begin{aligned} \|u_i\|_{L^2_2(\Omega')} &\leq C(\|u_i\|_{L^2} + \|\Delta_{cyl} u_i\|_{L^2}) \\ &= C(\|u_i\|_{L^2} + \|h(e^{u_i} - 1)\|_{L^2} + \|k_i\|_{L^2}). \end{aligned}$$

Combining this with the uniform bound on the $\|k_i\|_{L^2}$ and the uniform pointwise bound on $u_i$ by $u_\pm$, gives us a uniform bound on $\|u_i\|_{L^2_2(\Omega')}$. The Rellich lemma then provides a subsequence of the $u_i$ which converges in $L^2_1(\Omega')$. Hence, using the equation again, together with another elliptic estimate,

$$\begin{aligned} \|u_i\|_{L^2_4(\Omega)} &\leq C_1(\|u_i\|_{L^2(\Omega')} + \|\Delta_{cyl} u_i\|_{L^2_2(\Omega')}) \\ &\leq C_2(\|u_i\|_{L^2_1(\Omega')} + \|(\nabla^{(2)} u_i) h e^{u_i}\|_{L^2(\Omega')} + \|k_i\|_{L^2_2(\Omega')}), \end{aligned}$$

which is uniformly bounded, as we have just established. Thus, by the Rellich lemma once again, we have a subsequence of the $u_i$, so that

$$u_i \mapsto u \in L^2_3(\Omega).$$

The above method constructs a convergent subsequence over any precompact region $\Omega \subset R^o$. Exhausting $R^o$ by pre-compact regions, and extracting a diagonal subsequence, we find a subsequence $u_i$ converging in $L^2_{3,loc}$ to a function $u$.

Thus, we have our function $u$; we must verify that it has the desired properties. Continuity of $\mathcal{N}$ on $L^2_{3,loc}$, gives

$$\mathcal{N}(u) = k.$$

Continuity of $L^2_3 \hookrightarrow \mathcal{C}^0$, shows that the pointwise bounds $u_- \leq u_i \leq u_+$ hold in the limit as well, and $u_- \leq u \leq u_+$. Thus, we get that $u \in \mathcal{C}^0_\epsilon \subset L^2_{0,\epsilon}$. Moreover, the map

$$v \mapsto e^v - 1$$

clearly extends as a map from

$$\mathcal{C}^0_\epsilon \to \mathcal{C}^0_\epsilon.$$

These facts are enough to start the same bootstrap as above, only this time over the entire cylinder, to get $u \in L^2_{3,\epsilon}$ (and indeed, they are enough to get $u \in L^2_{\ell,\epsilon}$ for any $\ell$). □

Openness follows from an application of the inverse function theorem:

**Lemma 8.0.40.** *The set of $s \in [0,1]$ for which there is a $u_s$ solving Equation 65 is open.*

**Proof.** Consider the linearization of $\mathcal{N}$, $D\mathcal{N}_u$, given by

$$v \mapsto \Delta_{cyl} v + (he^u)v,$$

and viewed as a map

$$L^2_{3,\delta} \to L^2_{1,\delta},$$

where now $\delta$ is any weight with $0 < \delta < \epsilon$. The kernel of the operator is zero by the maximum principle. Also, its cokernel is zero for all sufficiently small positive weights. To see this, note that the Fredholm theory of [26] says $D\mathcal{N}_u$ is Fredholm on the doubly-weighted spaces $L^2_{k,\delta_-,\delta_+}$ (i.e. where the weight function looks like $e^{\delta_\pm |t|}$ for $\pm t > 1$) for all weights $\delta_-$ and $\delta_+$ satisfying

$$-2(\xi + \sqrt{\xi^2 + \lambda_-}) < \delta_- < 2(-\xi + \sqrt{\xi^2 + \lambda_-})$$

and

$$0 < \delta_+ < 4\xi,$$

where $\lambda_-$ is the smallest non-zero eigenvalue of the positive operator

$$\Delta_Y + h_-.$$

So, provided that

$$\delta < min(-\xi + \sqrt{\xi^2 + \lambda_-}, 2\xi),$$

we can connect the operator acting on $L^2_\delta = L^2_{\delta,\delta}$ to the same acting on $L^2_{\xi,-\xi}$ through a family of Fredholm operators. $D\mathcal{N}_u$ is easily seen to be self-adjoint on $L^2_{\xi,-\xi}$, so by the homotopy invariance of the index, $D\mathcal{N}_u$ on $L^2_{k,\delta}$ has index zero.

Thus, the inverse function theorem guarantees that around each $\sigma \in [0,1]$ for which there is some $u_\sigma \in L^2_{3,\epsilon}$ solving Equation 65, there is a real number $\delta \leq \epsilon$ and a neighborhood of values $s$ for which we can find $u_s \in L^2_{3,\delta}$ solving Equation 65. We must show that in fact each such $u_s \in L^2_{3,\epsilon} \subset L^2_{3,\delta}$.

But this follows from the *a priori* estimate $u_- < u_s < u_+$, which forces $u \in \mathcal{C}^0_\epsilon \subset L^2_{0,\epsilon}$. Now, we can proceed as in the end of the proof of Lemma 8.0.39.

□

The proof of Proposition 8.0.30 and hence Theorem 8.0.28 will be complete once we prove Lemma 8.0.38.

**Proof.** Fix a smooth, non-decreasing cut-off function $\lambda$ on $\mathbb{R}$ such that $\lambda(t) = 0$ for $t \leq 0$ and $\lambda(t) = 1$ for $t \geq 1$. and a smooth, non-negative function $b$ supported in $[0, \frac{1}{2}]$ with total integral 1.

We begin with two constructions which are used in the definition of $u_-$.

**Construction 8.0.41.** *Given $M > 0$, there is a function $\psi_1 \in \mathcal{C}^\infty(\mathbb{R})$ with*

$$\psi_1 \leq 0,$$

$$-(\frac{\partial}{\partial t}^2 + 2\xi\frac{\partial}{\partial t})\psi_1 \leq 0,$$

*and*

$$\psi_1(t) = \begin{cases} -\frac{M}{\xi}e^{-\xi t} & t \geq 0 \\ -\frac{M}{\xi} + Mt & t \leq -1 \end{cases}.$$

The two functions we are attempting to glue together (the linear function and the exponential) agree up to first order at $t = 0$. So, we smooth out the two derivatives using the bump function $\lambda$, and then integrate this smooth function to get $\psi_1$. Integrating the smoothing process decreases the value, which we compensate for by integrating the bump function $b$.

In terms of formulas, define

$$\psi_1(s,t) = -M\int_t^\infty \left(\lambda(st+s)e^{-\xi t} + (1-\lambda(st+s))\right)dt.$$

For all $s > 1$, we have that

$$\psi_1(s,t) = -\frac{M}{\xi} + Mt - \delta(s),$$

where

$$\lim_{s\mapsto\infty}\delta(s) = 0.$$

So, for large enough $s$, we can arrange that

$$\delta(s) < \max(\frac{M}{\sup((\frac{\partial}{\partial t} + 2\xi)b)}, \frac{M}{\sup(b)}).$$

If $s > 2$ is large enough to satisfy the above inequalities for $\delta(s)$, then

$$\psi_1(t) = \psi_1(s,t) + \delta(s)\int_t^0 b(t+1)$$

satisfies all requisite properties.

**Construction 8.0.42.** *Given $C, M, \epsilon$ positive reals, there is a function $\psi_2 \in \mathcal{C}^\infty(\mathbb{R})$ with*

$$\psi_2(t) \leq -Ce^{\epsilon t},$$

$$-(\frac{\partial}{\partial t}^2 + 2\xi\frac{\partial}{\partial t})\psi_2 \leq -(\frac{\partial}{\partial t}^2 + 2\xi\frac{\partial}{\partial t})\left(-(C+1)e^{\epsilon t}\right)$$

*and, for $t \geq 0$,*

$$\psi_2(t) = -C + Mt.$$

Given $\delta$, there is a unique value of $s$ such that the integral

$$-\int_{-\infty}^{t} \big(C\epsilon e^{\epsilon t}(1-\lambda(st+s)) + M\lambda(st+s)\big)\, dt$$

satisfies the requirements for $\psi_2$.

Now, we make a few preliminary choices.

Letting $m = \sup_Y h_-$, we can find a function $f \in \mathcal{C}^\infty(Y)$ with the property that on the set where $h_- \leq \frac{3m}{4}$, $\Delta_\Sigma f > 5\xi$. This can be done since $\Delta_\Sigma$ has a one-dimensional cokernel consisting of constant functions; we just have to find a function $g$ which is larger than $5\xi$ on the given subset and has total integral 0 over $\Sigma$, and we can solve $f = {\Delta_{cyl}}^{-1} g$. We can find such a function $g$ because we have assumed that $h_- \not\equiv 0$.

By translating towards $-\infty$, we can assume without loss of generality that

$$(h - h_-)|_{(-\infty,0]\times Y} < \frac{m}{4}. \tag{67}$$

Since $k \in \mathcal{C}^0_\gamma$, there is some constant $\kappa > 0$ with the property that

$$|k| \leq \kappa e^{-\gamma|t|}. \tag{68}$$

Our goal is to construct a function $u_-$ which satisfies

$$\mathcal{N}(u_-) \leq -\kappa e^{-\gamma|t|} \leq k.$$

Let $C\colon (0,1) \to (0,\infty)$ be any function satisfying

$$\lim_{\epsilon\mapsto 0} C(\epsilon)\epsilon = 0 \tag{69}$$

and

$$\lim_{\epsilon\mapsto 0} \frac{C(\epsilon)}{\log(\epsilon)} = -\infty \tag{70}$$

(for example, $C(\epsilon) = -\frac{\log(\epsilon)}{\sqrt{\epsilon}}$).

Letting

$$T(\epsilon) = \frac{2\log(\frac{2\kappa}{\epsilon\xi})}{\gamma} + \frac{1}{\xi}, \tag{71}$$

fix $\epsilon > 0$ small enough that it satisfies each of the following inequalities:

$$T(\epsilon) > 2+\frac{1}{\xi} \tag{72}$$

$$\log(\frac{3}{2}) > \epsilon \tag{73}$$

$$C(\epsilon) > \log(\sqrt{8}) \tag{74}$$

$$\xi > \epsilon, \tag{75}$$

$$\frac{\gamma}{2} > \epsilon, \tag{76}$$

$$\frac{1}{3\sup_Y |f|} > \epsilon, \tag{77}$$

$$\frac{m}{4\xi} > \epsilon, \tag{78}$$

$$\frac{C(\epsilon)}{2T(\epsilon)} > \kappa + \sup_{[1,2]\times Y} |\epsilon\Delta_{cyl}\lambda f e^{\epsilon t}| + \sup_{[0,1]} |\Delta_{cyl}\lambda e^{\epsilon t}|, \tag{79}$$

$$\frac{m}{8} > (\frac{4}{3}C(\epsilon)+1)(\epsilon^2+\xi\epsilon) + C(\epsilon)\epsilon|\sup \Delta_{cyl} f| \tag{80}$$

(Conditions 80 and 79 are satisfied for any $\epsilon$ sufficiently small, thanks to Equations 69 and 70 respectively.) Let $C = C(\epsilon), T = T(\epsilon)$ for the chosen value of $\epsilon$, let $\psi_1$ be the function provided by Construction 8.0.41, with $M = C/2T$, and let $\psi_2$ be the function provided by Construction 8.0.42. We will now verify that the function

$$u_-(t,x) = -\lambda(t)e^{-\epsilon t} - C\epsilon f(x)e^{\epsilon(t+T)}\lambda(T-t) + \left\{ \begin{array}{ll} \psi_1(t-T-\frac{1}{\xi}) & t \geq 0 \\ \psi_2(t+T) & t \leq 0 \end{array} \right.$$

satisfies

$$\mathcal{N}(u_-) \leq -\kappa e^{-\gamma|t|}.$$

We divide the problem into four cases, corresponding to four regions which cover $R^o$.

1. On the set where $h \geq \frac{m}{2}$ and $t \leq -T$, we have that

$$\begin{aligned}
\mathcal{N}(u_-) &\leq \Delta_{cyl}(-(C+1+\epsilon Cf)e^{\epsilon(t+T)}) + C\epsilon\Delta_{cyl} f e^{\epsilon(t+T)} \\
&\quad + \frac{m}{2}(e^{-C(1+\epsilon f)e^{\epsilon t}} - 1) \\
&\leq ((\frac{4}{3}C+1)(\epsilon^2+\xi\epsilon) + C\epsilon|\sup \Delta_{cyl} f|)e^{\epsilon(t+T)} \\
&\quad + \frac{m}{2}(e^{-\frac{2C}{3}e^{\epsilon(t+T)}} - 1) \\
&\leq \frac{m}{8}e^{\epsilon(t+T)} + \frac{m}{2}(e^{-\frac{2C}{3}e^{\epsilon(t+T)}} - 1).
\end{aligned}$$

(The first step follows from the assumption that $h \geq \frac{m}{2}$, together with the properties of Construction 8.0.42; the second step involves Inequality 77 and

Inequality 80.) The above quantity is in turn bounded by $-\frac{m}{8}e^{\epsilon(t+T)}$, since the function

$$x \mapsto \frac{x}{2} + (e^{-\frac{2C}{3}x} - 1)$$

is negative on $x \in [0,1]$, thanks to Inequality 74. In turn, combining the above inequalities with Inequalities 78, 76, and 68, we have that

$$\mathcal{N}(u_-) \leq -\frac{m}{8}e^{\epsilon(t+T)} < -\frac{m}{8}e^{\epsilon t} < -\frac{\epsilon\xi}{2}e^{\epsilon t}.$$

The desired inequality over this region then follows from the above, together with the observation that for $|t| < T - \frac{1}{\xi}$,

$$-\frac{\epsilon\xi}{2}e^{-\epsilon|t|} < -\kappa e^{-\gamma|t|}. \tag{81}$$

This is true because for $t$ in this range, Inequality 76 and 71 give

$$e^{(\gamma-\epsilon)|t|} \geq e^{\frac{\gamma}{2}t} \geq e^{\frac{\gamma}{2}(T-\frac{1}{\xi})} = \frac{\epsilon\xi}{2\kappa}.$$

2. When $t \leq -T$ and $h \leq \frac{m}{2}$, we have that $\Delta_{cyl} f \geq 5\xi$, due to our choice of $f$ and Inequality 67. Combining Inequalities 77 and 75, we have

$$(\epsilon^2 + 2\xi\epsilon)(1 + \epsilon f) < 4\xi\epsilon.$$

Combining this with Inequality 76, we get that over this region,

$$\begin{aligned} \mathcal{N}(u_-) &\leq C((\epsilon^2 + 2\epsilon\xi)(1+\epsilon f) - \epsilon\Delta_{cyl} f)e^{\epsilon(t+T)} \\ &\leq -C\xi\epsilon e^{\epsilon(t+T)} \\ &< -\xi\epsilon e^{\epsilon t}. \end{aligned}$$

(The last step follows from the fact that $T \geq 0$ and $C \geq 1$, which are of course weaker than Inequalities 72 and 74.) Once again, the desired inequality then follows from Inequality 81.

3. Over the region $[-T, T - \frac{1}{\xi}] \times Y$, $u_-$ is a pointwise non-positive function. This is clear over $[-T+1, T-\frac{1}{\xi}]$, since each term in the definition for $u_-$ is non-positive. For $t \in [-T, -T+1]$,we use Inequality 77 to get

$$|C\epsilon\lambda f e^{\epsilon(t+T)}| \leq \frac{Ce^{\epsilon}}{3} \leq \frac{C}{2}.$$

On the other hand, $\psi_2$ in the region $[-T, 0]$ is a linear function with slope $M$ and value $-C$ at $-T$, so in this region,

$$\psi_2(t) \leq -\frac{C}{2};$$

i.e. over all of $[-T, T - \frac{1}{\xi}] \times Y$,

$$u_- \leq 0,$$

and therefore

$$\mathcal{N}(u_-) \leq \Delta_{cyl}(u_-).$$

Now, the function which is $\psi_1(t-T)$ for $t \geq 0$ and $\psi_2(t+T)$ for $t \leq 0$ is a linear function with slope $M$, so

$$\begin{aligned} \mathcal{N}(u_-) \leq \Delta_{cyl}(u_-) &< -M + \sup(\Delta_{cyl}(\epsilon\lambda f e^{t\epsilon})) + \sup(\Delta_{cyl}\lambda) \\ &< -\kappa, \end{aligned}$$

which is exactly Inequality 79.

4. For all $t > T - \frac{1}{\xi}$, we have that $u_- = -e^{-\epsilon t} - \frac{M}{\xi} e^{-\xi(t-T)}$, so that

$$\begin{aligned} \mathcal{N}(u_-) &\leq (\epsilon^2 - 2\xi\epsilon)e^{-\epsilon t} \\ &\leq -\xi\epsilon e^{-\epsilon t} \\ &< -\kappa e^{-\gamma t} \end{aligned}$$

(The first step follows from the fact that $u_-$ is non-positive, the second follows from Inequality 75. and the third by Inequality 81.)

Verification that $u_+ = -u_-$ is an upper solution, i.e.

$$\mathcal{N}(u_+) \geq \kappa e^{\gamma|t|},$$

mirrors the above. The only conceptual difference occurs in Region 1, where we must verify an inequality of the form

$$-\frac{x}{2} + (e^{\frac{2C}{3}x} - 1) \geq 0$$

for $x \in [0,1]$. This is satisfied provided that $C \geq \frac{3}{4}$ (the condition which takes the place of Inequality 74). □

## 9. Deformation Theory on the Cylinder

Section 7 and Section 8 together set up an identification, on the level of topological spaces, between the moduli space of flows on $R^o$ and holomorphic curves in $R$. Now, we would like to compare the corresponding deformation theories. But before we can consider this problem, we embark on a digression about exponentially decaying cohomology for complex manifolds with cylindrical ends.

9.1. **Dolbeault Cohomology of Cylinders.** Let $\widehat{E}$ be a Hermitian line bundle over $R$, and let $\widehat{A}$ be a compatible connection on $\widehat{E}$ whose associated $\overline{\partial}$-operator is integrable. Denote the sheaf of $\overline{\partial}_A$-holomorphic sections of $\widehat{E}$ by $\mathcal{E}$. Given a real number $\delta$, consider the presheaf over $R^o$ of $(p,q)$-forms with values in $E = \widehat{E}|_{R^o}$, whose coefficients lie in

$$L^2_{\infty,\delta} := \bigcap_{k=1}^{\infty} L^2_{k,\delta}.$$

This can be pushed forward to $R$ and sheafified to form a sheaf over $R$ denoted by $A^{p,q}_\delta(E)$. Explicitly, a section of $A^{p,q}_\delta(U,E)$ for $U\subset R$ is a form $\omega\in\Omega^{p,q}(U-U\cap\Sigma_*,E)$ such that there is a covering $\{U_i\}$ of $U$ for which $\omega|_{U_i-U_i\cap\Sigma_*}$ has coefficients in $L^2_{\infty,\delta}$.

**Remark 9.1.1.** *Of course, the sheaf $A^{p,q}_\delta$ is a subsheaf of the sheaf of $\mathcal{C}^\infty_{\delta/2}$ forms, but the containment is proper: the function $e^{\frac{\delta\tau}{2}}$ is in $\mathcal{C}^\infty_{\delta/2}$, but it is not a section of $A^{0,0}_\delta$.*

Let

$$(E_\pm,A_\pm)=(\widehat{E},\widehat{A})|_{\Sigma_\pm}.$$

We can define another subsheaf of the sheaf of continuous $p,q$ forms $A^{p,q}_{\delta,\Sigma_*}(E)$ generated by sections of $A^{p,q}_\delta(E)$, (smooth sections of) $\pi^*(\Omega^{p,q}(\Sigma_+))$, and $\pi^*(\Omega^{p,q}(\Sigma_-))$. A section of $A^{p,q}_{\delta,\Sigma_*}(E)$ is a form $\omega\in\Omega^{p,q}(U-U\cap\Sigma_*,E)$ such that there is a covering $\{U_i\}$ of $U$ such that for each $i$, $U_i$ hits only one of the curves $\Sigma_-$, $\Sigma_+$, and, if $U_i\cap\Sigma_\mp=\emptyset$, then

$$\omega=\omega_\delta+\pi^*(\omega_\pm),$$

where $\omega_\delta\in A^{p,q}_\delta(E)$, and $\omega_\pm\in\Omega^{p,q}(U\cap\Sigma_\pm,E_\pm)$ is a (smooth orbifold) form.

Recall ([14]) that a sheaf $\mathcal{F}$ over a topological space $X$ is called *fine* if it admits partitions of unity. In this case, the sheaf cohomology groups $H^i(X;\mathcal{F})=0$ for all $i>0$. Moreover, if

$$0\longrightarrow\mathcal{E}\longrightarrow\mathcal{F}^0\xrightarrow{d^0}\mathcal{F}^1\xrightarrow{d^1}\dots$$

is a resolution of a sheaf $\mathcal{E}$ by fine sheaves, then there is a natural isomorphism between the sheaf cohomology of $\mathcal{E}$ with the cohomology of the chain complex of sections of $\mathcal{F}^i$,

$$H^i(X;\mathcal{E})\cong H^i(\Gamma(X,\mathcal{F}^*),d^*).$$

**Proposition 9.1.2.** *Let $\mathcal{E}\otimes\mathcal{I}_{\Sigma_*}$ denote the kernel of the restriction map*

$$\mathcal{E}\to\mathcal{E}|_{\Sigma_*}.$$

*For any $\delta$ with $0<\delta<2$, the complex $(A^{0,p}_\delta(E),\overline{\partial}_A)$ forms a fine resolution of the sheaf $\mathcal{E}\otimes\mathcal{I}_{\Sigma_*}$. Thus,*

$$H^i(\Omega^{0,*}_\delta(R^o,E),\overline{\partial}_A)\cong H^i(R,\mathcal{E}\otimes\mathcal{I}_{\Sigma_*}).$$

**Remark 9.1.3.** *Caution: $\mathcal{E}\otimes\mathcal{I}_{\Sigma_*}$ is* not *the tensor product of coherent sheaves over $R$; rather, it is a tensor product of orbifold sheaves. More precisely, if*

$$\phi\colon\widetilde{U}\to R$$

*is an orbifold coordinate chart with an action of $\mathbb{Z}/p\mathbb{Z}$ such that the orbifold sections of $\mathcal{E}$ correspond to the $\mathbb{Z}/p\mathbb{Z}$-invariant sections of $\widetilde{\mathcal{E}}$, then sections of $\mathcal{E}\otimes\mathcal{I}_{\Sigma_*}$ correspond to the $\mathbb{Z}/p\mathbb{Z}$-invariant sections of*

$$\widetilde{\mathcal{E}}\otimes_{\mathcal{O}_{\widetilde{U}}}\mathcal{I}_{\phi^{-1}(\Sigma_*)},$$

*where $\mathcal{I}_{\phi^{-1}(\Sigma_*)}$ is the ideal sheaf in $\mathcal{O}_{\widetilde{U}}$ of the pre-image of $\Sigma_*$.*

**Corollary 9.1.4.** *For sufficiently* $0 < \delta < 2$*, the complex* $(A^{0,*}_{\delta,\Sigma_*}, \overline{\partial}_A)$ *forms a fine resolution of the sheaf* $\mathcal{E}$*, so that*

$$H^i(A^{0,*}_{\delta,\Sigma_*}, \overline{\partial}_A) \cong H^i(R, \mathcal{E}).$$

**Proof.** The pair $(A^{0,*}_{\delta,\Sigma_*}(E), \overline{\partial}_A)$ is a clearly a complex of fine sheaves. Moreover, it naturally fits into a short exact sequence of sheaves

$$0 \longrightarrow A^{0,*}_{\delta}(E) \longrightarrow A^{0,*}_{\delta,\Sigma_*}(E) \longrightarrow A^{0,*}(\Sigma_*, \mathcal{E}|\Sigma_*) \longrightarrow 0.$$

Since both $A^{0,*}_{\delta}(E)$ and $A^{0,*}(\Sigma_*, E|_{\Sigma_*})$ are acyclic complexes of sheaves (Proposition 9.1.2), so is $A^{0,*}_{\delta,\Sigma_*}(E)$.

It remains therefore to identify the kernel $K$ of

$$\overline{\partial}\colon A^{0,0}_{\delta,\Sigma_*}(E) \to A^{0,1}_{\delta,\Sigma_*}(E).$$

Restricting holomorphic sections of $\mathcal{E}$ from $R$ to $R^o \subset R$ induces a sheaf map

$$\mathcal{E} \to K,$$

since $\mathcal{E}$ is generated by the subsheaves of holomorphic functions which vanish along $\Sigma_*$, $\mathcal{E} \otimes \mathcal{I}_{\Sigma_*}$, and the sheaf of holomorphic functions which are pull-backs of holomorphic functions defined over $\Sigma_*$. It follows then that $K \cong \mathcal{E}$, from the short exact sequence above and the proposition. □

The proof of Proposition 9.1.2 is tantamount to showing a version of the $\overline{\partial}$-Poincaré lemma, which we build to in the following sequence of lemmas. But first, some notation. Let $D_r$ denote the disk of radius $r$ in the complex plane, and $C_T = [T, \infty) \times S^1$ be the cylinder with its natural complex structure.

Given $T \in \mathbb{R}$, let $\mathcal{H}_\delta(D_r \times C_T)$ denote the space of holomorphic functions which lie in $L^2_{\infty,\delta}$ (i.e. functions in $L^2_{\infty,\delta}(D_r \times C_T)$ which can be holomorphically extended to $D_{r'} \times C_{T'}$ for some $r' > r$, $T' < T$). By elliptic regularity, a holomorphic function lies in this space if and only if it lies in $L^2_{0,\delta}$. Let

$$c\colon D_r \times C_T \to D_r \times D_{e^{-T}}$$

given by

$$(w, z) \mapsto (w, e^{-z}).$$

**Lemma 9.1.5.** *The space* $\mathcal{H}_\delta(D_r \times C_T)$ *is identified via precomposition with* $c$ *with the space of holomorphic functions on* $D_r \times D_{e^{-T}}$ *which vanish to order greater than* $\delta/2$ *along the subspace* $D_r \times \{0\}$*.*

**Proof.** Suppose $f \in \mathcal{H}_\delta(D_r \times C_T)$. Then, by the Rellich lemma,

$$\lim_{t \mapsto +\infty} |f| e^{t\delta/2} = 0,$$

so the function $\widetilde{f}$ we get by extending $c^*f$ by zero over $D_r \times \{0\}$ is continuous over $D_r \times D_{e^{-T}}$ and holomorphic away from $D_r \times \{0\}$. By the removable singularities theorem of complex analysis, it follows then that $\widetilde{f}$ is holomorphic over all of $D_r \times D_{e^{-T}}$. In fact, the same argument says that if $n \leq \frac{\delta}{2}$, then $\widetilde{f}(z)/z^n$ extends over all of $D_r \times D_{e^{-T}}$ as a function which vanishes at along $D_r \times \{0\}$.

Conversely, by Taylor's theorem, if $f$ is a holomorphic function defined over $D_r \times D_{e^{-T}}$ which vanishes to order greater than $\delta/2$ along $D_r \times \{0\}$, we can write $f = z^n g$ for $n > \delta/2$ and $g$ a holomorphic function over $D_r \times C_T$. It follows then that

$$\|f\|_{L^2_\delta} \leq \pi r^2 \|g\|_{L^\infty} (\int_0^\infty e^{-2n+\delta} dt) < \infty.$$

□

**Definition 9.1.6.** *Given a smooth function*

$$g\colon D_r \times C_T \to \mathbb{C},$$

*let*

$$Q_z g\colon D_r \times C_T \to \mathbb{C}$$

*denote the function given by*

$$Q_z g(w,z) = \frac{1}{2\pi i} \int_{C_T} \frac{f(w,\eta) d\eta \wedge d\overline{\eta}}{1 - e^{z-\eta}},$$

*whenever the above integral makes sense.*

**Lemma 9.1.7.** *If $\frac{\delta}{2} \notin \mathbb{Z}$, then $Q_z$ induces a bounded map*

$$Q_z\colon L^2_{\infty,\delta}(D_r \times C_T) \to L^2_{\infty,\delta}(D_r \times C_T)$$

*with the property that*

$$\frac{\partial}{\partial \overline{z}} Q_z f(w,z) = f(w,z).$$

**Proof.** For $f\colon C_T \to \mathbb{C}$, define

$$Qf(z) = \frac{1}{2\pi i} \int_{C_T} \frac{f(\eta) d\eta \wedge d\overline{\eta}}{1 - e^{z-\eta}}.$$

Write the Fourier expansion for a function $f$ as

$$f(t + i\theta) = \sum_{n \in \mathbb{Z}} \widehat{f}_n(t) e^{2\pi i n\theta}.$$

Using Cauchy's integral formula, we have for any $a \in \mathbb{C}$ with $|a| \neq 1$ that

$$\frac{1}{2\pi} \int \frac{e^{-in\theta} d\theta}{1 - ae^{i\theta}} = \begin{cases} a^n & \text{if } |a| < 1 \text{ and } n \geq 0 \\ -a^n & \text{if } |a| > 1 \text{ and } n \leq 0 \\ 0 & \text{otherwise} \end{cases},$$

so that, writing $\eta = s + i\varphi$

$$\begin{aligned}
\widehat{Qf}_n(t) &= \frac{1}{2\pi}\int e^{-in\theta}d\theta\left(\frac{1}{2\pi i}\int_{C_T}\frac{f(\eta)d\eta\wedge d\overline{\eta}}{1-e^{z-\eta}}\right)\\
&= \frac{1}{2\pi i}\int_{C_T} f(\eta)d\eta\wedge d\overline{\eta}\left(\frac{1}{2\pi}\int\frac{e^{-in\theta}d\theta}{e^{t+i\theta-\eta}-1}\right)\\
&= \frac{1}{2\pi i}\int_{C_T} f(s+i\varphi)(-2ids\wedge d\varphi)\begin{cases} e^{n(t-s-i\varphi)}\chi_{(t,\infty)}(s) & \text{if } n\geq 0\\ -e^{n(t-s-i\varphi)}\chi_{(T,t)}(s) & \text{if } n<0\end{cases}\\
&= 2\begin{cases} -\int_t^\infty \widehat{f}_n(s)e^{n(t-s)}ds & \text{if } n\geq 0\\ \int_T^t \widehat{f}_n(s)e^{n(t-s)}ds & \text{if } n<0\end{cases}
\end{aligned}$$

Using the obvious isometry for each $k\in\mathbb{Z}$

$$\phi_\delta\colon L^2_{k,\delta}\cong L^2_k$$

given by

$$\phi_\delta(f) = fe^{\frac{t\delta}{2}},$$

we see that $Q$ corresponds to the parametrix for $\phi_\delta\frac{\partial}{\partial\overline{z}}\phi_\delta^{-1}$ acting on the cylinder $C_T$ with APS boundary conditions constructed in Proposition 2.5 of [2]. It follows from this proposition, then, that $Q$ induces a bounded map $Q\colon L^2_{k,\delta}(C_T)\to L^2_{k+1,\delta}(C_T)$ whenever $\frac{\delta}{2}\notin\mathbb{Z}$.

Since

$$\int\mu_{D_r}(w)\int\mu_{C'_T}(z)|Q_zf(w,z)|^2e^{t\delta}\leq C\int\mu_{D_r}(w)\int\mu_{C'_T}(z)|f(w,z)|^2e^{t\delta},$$

we have by Fubini's theorem that $Q_z$ maps $L^2_{0,\delta}(D_r\times C_T)$ to $L^2_{0,\delta}(D_r\times C'_T)$. Similarly, since $Q_z$ commutes with differentiation, we have that $Q$ extends to a map from $L^2_{k,\delta}$ to itself for each $k$, and hence a map from $L^2_{\infty,\delta}$ to itself. □

We now prove Proposition 9.1.2

**Proof.** Clearly, $A^{0,p}_\delta$ is a fine sheaf. We must show that given $U\subset R$, $\omega\in\Omega^{0,p}_{\infty,\delta}(U,E|_U)$ with $\overline{\partial}_A\omega=0$, and $x\in U$, there is an open neighborhood $V\subset U$ of $x$ and a form $\eta_V\in\Omega^{0,p-1}_{\infty,\delta}(V;E|_V)$ such that

$$\overline{\partial}_A\eta = \omega|_V.$$

When $x\notin\Sigma_*$, this follows from the $\overline{\partial}$-Poincaré lemma. When $x\in\Sigma_*$, we imitate the usual proof of the this lemma.

Suppose that $x\in\Sigma_*$, but $x$ is not a singular point. Then, there is a neighborhood of $x\in R$ of the form $D_r\times D_{e^{-T}}$ with respect to which $x$ corresponds to the origin, $\Sigma_*$ corresponds to $D_r\times\{0\}$, holomorphic sections of $E$ correspond to holomorphic functions (under a holomorphic trivialization of $\mathcal{E}|_{D_r\times C_T}$), and $\omega|_{D_r\times D_{e^{-T}-0}}$ corresponds to a form in $L^2_{\infty,\delta}(D_r \times C_T; \Lambda^{0,p})$ satisfying $\overline{\partial}\omega = 0$. We would like to find a form $\eta \in L^2_{\infty,\delta}(D_r \times C_T)$ with $\overline{\partial}\eta = \omega$.

When $p = 2$, write $\omega = f d\overline{w} \wedge d\overline{z}$, and let $\eta = -Qzf d\overline{w}$. The fact that $\eta$ lies in the requisite Sobolev spaces and satisfies $\overline{\partial}\eta = \omega$ follows from Lemma 9.1.7.

When $p = 1$, write $\omega = g d\overline{w} + h d\overline{z}$. We reduce to the case where $h = 0$ by subtracting $\overline{\partial} Q_z h$.

When $p = 1$, and $\omega$ is of the form

$$\omega = g d\overline{w},$$

let

$$Q_w g(w, z) = \frac{1}{2\pi i} \int_{D_r} \frac{g(\xi, z)}{\xi - w} d\xi \wedge d\overline{\xi}.$$

By Cauchy's integral formula, $Q_w g$ is a smooth function with $\frac{\partial}{\partial \overline{w}} Q_w g = g$. The condition that $\overline{\partial}\omega = 0$ is equivalent to the condition that $g$ is holomorphic in the $z$ variable, a property which is clearly preserved by the operator $Q_w$; thus, $\overline{\partial} Q_w g = \omega$.

To verify the integrability condition on $Q_w g$, recall Young's convolution inequality

$$\|f * g\|_{L^2} \leq \|f\|_{L^1} \|g\|_{L^2},$$

which gives us that

$$\int_{D_r} |Q_w g(w, z)|^2 dw \wedge d\overline{w} \leq 2r^2 \int_{D_r} |g(w, z)|^2 dw \wedge d\overline{w}.$$

Thus, by Fubini's theorem, $Q_w g \in L^2_{0,\delta}$. Using elliptic regularity for $\overline{\partial}$ gives us then that $\eta \in L^2_{\infty,\delta}(D_{r'} \times C_T)$, for any $r' < r$.

Thus, we have verified the exactness of the sequence of sheaves

$$A^{0,0}_\delta(E) \xrightarrow{\overline{\partial}_A} A^{0,1}_\delta(E) \xrightarrow{\overline{\partial}_A} A^{0,2}_\delta(E) \longrightarrow 0 \tag{82}$$

over the smooth locus of $R$.

When $x$ is a singular point, we can identify a neighborhood $U$ of $x$ with the quotient of $D_r \times D_{e^{-T}}$ by the $\mathbb{Z}/p\mathbb{Z}$ action

$$\zeta \times (w, z) = (\zeta w, \zeta^p z).$$

For each $j = 0, ..., p - 1$, we can identify

$$L^2_{\infty,\delta}(U - \Sigma_* \cap U, E|_{U - \Sigma_* \cap U}) \cong L^2_{\infty,\delta}(\frac{D_r \times C_T}{\mathbb{Z}/p\mathbb{Z}}, \mathcal{O}_j)$$

with the space of $\mathbb{Z}/p\mathbb{Z}$-equivariant sections of

$$L^2_{\infty,\delta}(D_r \times C_T),$$

i.e. those functions which satisfy

$$\zeta^j f(w, z) = f(\zeta w, z - q \log \zeta).$$

This is true because the $\mathbb{Z}/p\mathbb{Z}$ acts by isometries which do not change the $t$ coordinate. Letting $\zeta^*$ denote pull-back (of forms) under the map which defines the $\mathbb{Z}/p\mathbb{Z}$ action, and letting $w, z$ denote the $D_r$ and $C_T$ coordinates on $D_r \times C_T$ respectively, note that

$$\zeta^* d\overline{w} = \zeta^{-1} d\overline{w}$$

and

$$\zeta^* d\overline{z} = d\overline{z}.$$

Exactness of the Diagram 82 then follows once we have demonstrated that the parametrices $Q_w$ and $Q_z$ constructed earlier commute with the $\mathbb{Z}/p\mathbb{Z}$ action:

$$\begin{aligned}
\zeta^* Q_w(f d\overline{w})(w,z) &= \frac{1}{2\pi i}\int_{D_r} \frac{f(\xi, z - q\log(\zeta))d\xi \wedge d\overline{\xi}}{\zeta w - \xi} \\
&= \frac{1}{2\pi i\zeta}\int_{D_r} \frac{f(\zeta\xi, z - q\log(\zeta))d\xi \wedge d\overline{\xi}}{w - \xi} \\
&= (Q_w \zeta^* f)\zeta^*(d\overline{w}) \\
\zeta^* Q_z(f d\overline{z})(w,z) &= \frac{1}{2\pi i}\int_{C_T} \frac{f(\zeta w, \eta - q\log(\zeta))d\eta \wedge d\overline{\eta}}{1 - e^{z-\eta}} \\
&= (Q_z \zeta^* f)\zeta^*(d\overline{z})
\end{aligned}$$

The fact that Sequence 82 is a resolution of the sheaf $\mathcal{E} \otimes \mathcal{I}_{\Sigma_*}$ when $0 < \frac{\delta}{2} < 1$ is an immediate consequence of Lemma 9.1.5. □

**Remark 9.1.8.** *More generally, for $\frac{\delta}{2} \notin \mathbb{Z}$, the above sequence provides a resolution of $\mathcal{E} \otimes \mathcal{I}_{\Sigma_*}^{\lceil \frac{\delta}{2} \rceil}$, the subsheaf of $\mathcal{E}$ generated by sections which vanish to order $\lceil \frac{\delta}{2} \rceil$ along $\mathcal{I}_{\Sigma_*}$.*

**Remark 9.1.9.** *Of course, the arguments given above generalize in an obvious way to any dimension, replacing the curve $\Sigma_*$ by a (complex) codimension-one sub-orbifold.*

**Example 9.1.10.** *If $X$ is a compact, complex curve, and $P \subset X$ is a finite collection of points, then we can give $X - P$ a complete metric with finitely many cylindrical ends. The above shows that the $L^2_\delta$-index of the $\overline{\partial}$ operator acting on $\mathcal{O}_{X-P}$ is the holomorphic Euler characteristic $\chi(X, \mathcal{I}_P)$, which is easily computed via the short exact sequence*

$$0 \longrightarrow \mathcal{I}_P \longrightarrow \mathcal{O}_X \longrightarrow \oplus_{p\in P}\mathbb{C}_p \longrightarrow 0$$

*to be $1 - g - |P|$. This agrees with the computation appearing in [2] of the extended $L^2$-index of the $\overline{\partial}$ operator.*

In Section 10, we need a version of Proposition 9.1.2 in the case where $(E, \overline{\partial}_A)$ does not quite extend over $R$. Assume rather that we are given $(E, A)$, a Hermitian line bundle with connection over $R^o$, which naturally extends over $R - \Sigma_+$, together with an isomorphism

$$j\colon \pi^*(E_+)|_{R^o} \to E|_{R^o}$$

such that

$$j^*(A) = A_+ + i\delta_0\eta,$$

where $\delta_0 \in (0,1)$, $E_+$ is an orbifold bundle over $\Sigma$, and $A_+$ is a connection which extends over $R - \Sigma_-$. Then, we can form the line bundle $\mathcal{E}_0$ over $R$ given by patching $A$ to $A_+$ using the complex gauge transformation $t \mapsto e^{-\delta_0 t}$; i.e. the $\overline{\partial}_{A_+}$-holomorphic section $\Phi$ of $\pi^*(E_+)$ is identified with the $\overline{\partial}_A$-holomorphic section $e^{-t\delta_0}j^*(\Phi)$ of $E|_{[T,\infty)}$.

**Corollary 9.1.11.** *For any real $\delta \in (0, \delta_0)$, the complex $(A^{0,p}_\delta(E), \overline{\partial}_A)$ forms a fine resolution of the sheaf $\mathcal{E}_0 \otimes \mathcal{I}_{\Sigma_-}$.*

**Proof.** We must prove the following analogue of Lemma 9.1.5: Consider the complex structure $\overline{\partial}_A$ on the trivial bundle $D_r \times C_T$ induced by the connection $d - i\delta_0 d\theta$. Then, for any real $\delta \in (0, \delta_0)$, the map which takes a function $\phi$ over $D_r \times D_{e^{-T}}$ to the function on $D_r \times C_T$ defined by

$$(w, z) \mapsto e^{-t\delta_0}\phi(w, e^{-z})$$

induces an identification between the space of holomorphic functions on $D_r \times D_{e^{-T}}$ and the space of $\overline{\partial}_A$-holomorphic function on $D_r \times C_T$ which lie in $L^2_\delta$. This implies that the complex $(A^{0,*}_\delta(E), \overline{\partial}_A)$, which was shown to be a resolution, resolves $\mathcal{E}_0 \otimes \mathcal{I}_{\Sigma_-}$. ∎

9.2. **Comparison of Deformation Theories.** With this background in place, we turn to the identification of the deformation theories.

**Definition 9.2.1.** *Let $U_1, U_2 \subset \mathcal{C}^\infty(Y)$ be a pair of subsets, $\delta \in \mathbb{R}$, $k \in \mathbb{Z} \cup \infty$, $k \geq 2$. The space of $U_1, U_2$-*extended $L^2_{k,\delta}$ functions*, denoted $L^2_{k,\delta,U_1,U_2}$, is the space of functions which can be written as*

$$f_0 + f_1 + f_2,$$

*where $f_0 \in L^2_{k,\delta}$, $f_1(t, y) = \lambda(t)g_1(y)$, $f_2(t, y) = \lambda(-t)g_2(y)$, for $g_i \in U_i$, and $\lambda$ a cut-off function with $\lambda \equiv 0$ for $t \leq 0$, $\lambda \equiv 1$ for $t \geq 1$.*

**Example 9.2.2.** *The global sections of $A^{p,q}_{\delta,\Sigma_*}(E)$ are the $L^2_{\infty,\delta}$ sections of the bundle $\Lambda^{p,q}$ extended by the smooth forms in $\pi^*(\Omega^{p,q}(E_-)), \pi^*(\Omega^{p,q}(E_+))$.*

Let $\mathcal{C}_1$ and $\mathcal{C}_2$ be a pair of non-degenerate, irreducible critical manifolds (Definition 5.9.3) for cs and let $[A,\Phi]\in\mathcal{M}(\mathcal{C}_1,\mathcal{C}_2)$ correspond to a gradient flow line between them; i.e.

$$\lim_{t\mapsto\pm\infty}[A,\Phi]|_{\{t\}\times Y}=[B_\pm,\Psi_\pm]$$

with $[B_-,\Psi_-]\in\mathcal{C}_1$, $[B_+,\Psi_+]\in\mathcal{C}_2$ (we return to the case of flows to the trivial connection in Subsection 10.1). A neighborhood of $[A,\Phi]$ in $\mathcal{M}(\mathcal{C}_1,\mathcal{C}_2)$ can be described as follows. Choose spaces $U_1,U_2$ of $(B_+,\Psi_+)$ $(B_-,\Psi_-)$ in $\mathcal{C}(Y)$ which project to open neighborhoods of $[B_+,\Psi_+]$, $[B_-,\Psi_-]$ in $\mathcal{C}_1,\mathcal{C}_2\subset\mathcal{B}(Y)$ respectively, and consider the configuration space $\mathcal{B}$ of $L^2_{k,\delta,U_1,U_2}$ perturbations of $(A,\Phi)$, modulo the action of the gauge group $L^2_{k+1,\delta}$-gauge transformations which are 1 at infinity. The Seiberg-Witten equations, viewed as a section of a bundle over $\mathcal{B}$, cut out a neighborhood of $[A,\Phi]\in\mathcal{M}(\mathcal{C}_1,\mathcal{C}_2)$. The tangent space of $\mathcal{B}$ is the quotient of the space of $T_{(B_-,\psi_-)}U_1$, $T_{(B_-,\psi_-)}U_2$-extended $L^2_{k,\delta}$ sections of $\Lambda^1(i\mathbb{R})\times W$ by the image of $L^2_{k+1,\delta}$ under the derivative of the gauge group action. The linearization

$$D_{[A,\Phi]}\mathrm{sw}\colon \mathrm{T}_{[\mathrm{A},\Phi]}\mathcal{B}\to\Omega^+(R^o,i\mathbb{R})\oplus\Gamma(R^o,\mathrm{W}^-),$$

given by

$$(a,\phi)\mapsto(d^++\widehat{\sigma}(\phi,\Phi),\mathbf{D}_A\phi+a\cdot\Phi),$$

where the range is given the $L^2_{k-1,\delta}$ topology, is Fredholm for all sufficiently small $\delta>0$ (and even when $\delta=0$, if the two critical manifolds are non-degenerate, isolated points). Here, $\widehat{\sigma}$ denotes the bilinear form corresponding to the quadratic function $\sigma$ appearing in the four-dimensional Seiberg-Witten Equaitons (Equation 11). Thus, the deformation complex, $(C^*_{sw},d^*_{sw})$, for $\mathcal{M}(\mathcal{C}_1,\mathcal{C}_2)$ is

$$\Omega^0(R^o,i\mathbb{R})\longrightarrow\Omega^1(R^o,i\mathbb{R})\oplus\Gamma_\delta(R^o,W^+)\longrightarrow\Omega^+(R^o,i\mathbb{R})\oplus\Gamma_\delta(R^o,W^-),\tag{83}$$

This complex is graded so that the middle term in the above sequence has degree zero. The terms in degree $i\neq0$ are required to lie in the Sobolev space $L^2_{\delta,k-i}$, while the term in degree zero lies in $L^2_{\delta,k}$ extended by $T_{(B_1,\psi_1)}U_1\oplus T_{(B_2,\psi_2)}U_1$.

The Morse theory of [27] adapted to the Seiberg-Witten context reads as follows:

**Theorem 9.2.3.** *Let $\mathcal{C}_1,\mathcal{C}_2$ be a pair of non-degenerate critical manifolds for* cs*. Given any $(A,\Phi)\in\mathcal{M}(C_1,C_2)$, with*

$$\lim_{t\mapsto\pm\infty}[A,\Phi]_{\{t\}\times Y}=\pi^*[B_\pm,\Psi_\pm],$$

*there is a Kuranishi map*

$$\psi\colon H^0(C^*_{sw},d^*_{sw})\to H^1(C^*_{sw},d^*_{sw})$$

*describing a neighborhood of $[A,\Phi]$ in $\mathcal{M}(\mathcal{C}_1,\mathcal{C}_2)$.*

**Definition 9.2.4.** *Let $\mathcal{C}_1, \mathcal{C}_2$ be a pair of non-degenerate critical manifolds for* cs. *The* formal dimension *of the space of flows from $\mathcal{C}_1$ to $\mathcal{C}_2$, written e-dim$(\mathcal{M}(\mathcal{C}_1,\mathcal{C}_2))$ is the quantity defined by the formula:*

$$\text{e-}dim(\mathcal{M}(\mathcal{C}_1,\mathcal{C}_2)) := \chi(H^*(C^*_{sw}, d^*_{sw})).$$

When $\mathcal{C}_i$ $i = 1, 2$ are both non-degenerate, isolated points, the above quantity is the relative degree appearing in the Floer complex (see Section 13).

We have given an identification between the space of flows $\mathcal{M}(\mathcal{C}^+(E_1), \mathcal{C}^+(E_2))$ and an open subset of the space of divisors $\mathcal{D}(E_1, E_2)$ interpolating between $E_1$ and $E_2$. In this section, we compare Kuranishi descriptions of these spaces.

**Theorem 9.2.5.** *Suppose that $E_i$ $i = 1, 2$ correspond to irreducible critical manifolds $\mathcal{C}^+(E_i)$. Then, like the moduli space of divisors, the moduli space $\mathcal{M}(\mathcal{C}^+(E_1), \mathcal{C}^+(E_2))$ is locally modeled about any point $(A, \Phi) \in \mathcal{M}(\mathcal{C}^+(E_1), \mathcal{C}^+(E_2))$ on the zeros of a map*

$$H^0(D, \mathcal{E}|_D) \to H^1(D, \mathcal{E}|_D),$$

*where $\mathcal{E}$ is the holomorphic structure induced by $(A, \Phi)$ on the line bundle interpolating $E_1$ and $E_2$, and $D = \alpha^{-1}(0)$ is the zero locus of the associated section. In particular,*

$$\text{e-}dim\mathcal{M}(\mathcal{C}^+(E_1), \mathcal{C}^+(E_2)) = \text{e-}dim\mathcal{D}(E_1, E_2). \tag{84}$$

**Remark 9.2.6.** *Given a section $\alpha$ over an orbifold line bundle $\mathcal{E}$, the sheaf $\mathcal{E}/\alpha\mathcal{O}_R$ is supported along $D = \alpha^{-1}(0)$. Indeed, the pair $(D, \mathcal{E}/\alpha\mathcal{O}_R)$ is naturally a (possibly singular) orbifold curve equipped with an orbifold line bundle, which we write as $(D, \mathcal{E}|_D)$.*

Using the identifications of Isomorphism 47, together with the identifications ([6])

$$\begin{array}{ccc} \Omega^1(i\mathbb{R}) & \xrightarrow{d^+} & \Omega^+(i\mathbb{R}) \\ \cong\downarrow & & \cong\downarrow \\ \Omega^{0,1} & \xrightarrow{a\mapsto(i\mathrm{Im}\Lambda\partial a,\overline{\partial}a)} & \Omega^0 \oplus \Omega^{0,2}, \end{array}$$

we get a natural map of complexes from the Seiberg-Witten deformation complex (Diagram 83) to the complex $(C^*_r, d^*_r)$ defined by:

$$\Omega^{0,0} \longrightarrow \Omega^{0,1} \oplus (\Omega^{0,0} \oplus \Omega^{0,2})(E) \longrightarrow \Omega^{0,2} \oplus \Omega^{0,1}(E), \tag{85}$$

where once again, the terms in degree $i \neq 0$ lie in $L^2_{k-i,\delta}$, while the term in degree 0 lies in $L^2_{k-i,\delta}$ is extended by $T_{(B_1,\psi_1)}U_1 \oplus T_{(B_1,\psi_1)}U_1$.

To identify the cohomology groups of $(C^*_r, d^*_r)$ and $(C^*_{sw}, d^*_{sw})$, we appeal to the following standard fact from homological algebra.

**Lemma 9.2.7.** *Let $(C^*, d^*)$ be a cochain complex, $i \in \mathbb{Z}$, and $A, B$ be Abelian groups together with maps*

$$f\colon A \to C^i \quad \text{and} \quad g\colon C^i \to B$$

*with the properties that*

$$d^i \circ f = g \circ d^{i-1} = 0$$

*and*

$$f \circ g \colon A \to B$$

*is an isomorphism. Then the natural map of chain complexes*

$$\begin{array}{ccccccccc} \dots & \longrightarrow & C^{i-1} & \xrightarrow{d^{i-1}} & C^i & \xrightarrow{g\oplus d^i} & B \oplus C^{i+1} & \longrightarrow & \dots \\ & & \downarrow & & \downarrow & & \downarrow & & \\ \dots & \longrightarrow & A \oplus C^{i-1} & \xrightarrow{f+d^{i-1}} & C^i & \xrightarrow{d^i} & C^{i+1} & \longrightarrow & \dots \end{array}$$

*induces an isomorphism in cohomology.*

**Lemma 9.2.8.** *The natural map from the Seiberg-Witten deformation complex $(C^*_{sw}, d^*_{sw})$ to the complex $(C^*_r, d^*_r)$ induces an isomorphism in cohomology.*

**Proof.** We appeal to Lemma 9.2.7, with $A = ReL^2_{k+1}\Omega^{0,0}_\delta$ and $B = L^2_{k-1,\delta}\Omega^0\omega \subset \Omega^+_\delta(R^o, i\mathbb{R})$. We must verify that the map induced from $A$ to $B$ by $\Lambda d^1_{sw} \circ d^0_r$ is an isomorphism. But

$$\begin{aligned} \Lambda d^1_{sw} d^0_r f &= \Lambda d^1_{sw}(2\overline{\partial} f, f \otimes \Phi) \\ &= 2(iIm\Lambda\partial\overline{\partial} + |\Phi|^2)f, \end{aligned}$$

which is (half) the differential of the operator $\mathcal{N}$ encountered in Section 8. This was shown to be an isomorphism in Lemma 8.0.40. □

**Lemma 9.2.9.** *The natural map from the complex appearing in Diagram 85 to the complex $(C^*_{\mathbb{C}}, d^*_{\mathbb{C}})$ given by*

$$\Omega^{0,0} \longrightarrow \Omega^{0,1} \oplus \Omega^{0,0}(E) \longrightarrow \Omega^{0,2} \oplus \Omega^{0,1}(E) \longrightarrow \Omega^{0,2}(E), \tag{86}$$

*with the conventions for Sobolev completion the same as for the complex $(C^*_{sw}, d^*_{sw})$, induces an isomorphism in cohomology.*

**Proof.** (Of Theorem 9.2.5) Once again, we appeal to Lemma 9.2.7, this time with $A = L^2_k\Omega^{0,2}_\delta(E)$ and $B = L^2_{k-2}\Omega^{0,2}_\delta(E)$. Note that

$$d^2_{\mathbb{C}} \circ d^1_r|_A = d^2_{\mathbb{C}} \circ (d^2_{\mathbb{C}})^* \colon L^2_{k,\delta}\Omega^{0,2}(E) \to L^2_{k,\delta}\Omega^{0,2}(E),$$

where

$$(d^2_{\mathbb{C}})^* \colon L^2_{k,\delta}\Omega^{0,2}(E) \to C^1_{\mathbb{C}}$$

denotes the formal ($L^2$) adjoint of the differential operator defining $d^2_{\mathbb{C}}$.

To see injectivity of $d^2_{\mathbb{C}} \circ (d^2_{\mathbb{C}})^*$, use integration by parts to identify its kernel with

$$\mathrm{Ker}(d^2_{\mathbb{C}})^* = \mathrm{Ker} D_{[A,\Phi]}\mathrm{sw}|_{S^-};$$

but we have shown in Section 6 that the $S^-$-component of any element in the kernel of $D_{[A,\Phi]}\mathrm{sw}$ vanishes, when $\Phi$ is a non-vanishing section of $S^+$. To see surjectivity, note that as a map from $L^2_{k,0}$ to $L^2_{k-2,0}$, $d^2_{\mathbb{C}} \circ (d^2_{\mathbb{C}})^*$ is Fredholm, as it asymptotically approaches

$$(-\frac{d^2}{dt^2} + (D_{(B_\pm,\Psi_\pm)}\mathrm{sw}_3)^2)|_{S^-},$$

(where $(B_\pm, \Psi_\pm)$ are the limiting values at $\pm\infty$ of the pair $(A, \Phi)$) and the operator $D_{(B_\pm,\Psi_\pm)}\mathrm{sw}_3$ has no kernel in $S^-$, since $\Psi_\pm$ are nonvanishing sections of $S^+$. Moreover, when $\delta = 0$, $d^2_{\mathbb{C}} \circ (d^2_{\mathbb{C}})^*$ is formally self-adjoint, hence it has index zero. Thus, for small $\delta > 0$, the operator remains Fredholm of index zero. □

**Proof.** A subcomplex of a complex is said to *carry the cohomology* if the inclusion map induces an isomorphism in cohomology. Consider the subcomplex $(C^*_\delta, d^*_\delta) = (C^*_{\mathbb{C}}, d^*_{\mathbb{C}}) \cap L^2_\delta$, the subcomplex consisting of un-extended Sobolev sections. For small $\delta$, this is a Fredholm complex, so its cohomology is carried by harmonic representatives. In particular, elliptic regularity tells us that the cohomology is carried by the subcomplex of sections which lie in $L^2_{\infty,\delta}$. Moreover, we have a short exact sequence of complexes

$$0 \longrightarrow (C^*_\delta, d^*_\delta) \longrightarrow (C^*_{\mathbb{C}}, d^*_{\mathbb{C}}) \longrightarrow T_{(B_-,\psi_-)}U_1 \oplus T_{(B_+,\psi_+)}U_2 \longrightarrow 0,$$

where the tangent space above is to be thought of as a complex concentrated in degree zero. Thus, by an easy application of the five-lemma, the cohomology of the complex $(C^*_{\mathbb{C}}, d^*_{\mathbb{C}})$ is carried by $L^2_{\infty,\delta}$ forms (appropriately extended in degree zero). We have seen (Subsection 5.6 together with Subsection 5.7) that the tangent space to the critical manifolds carries the cohomology of the deformation complex for divisors in $\Sigma_*$, $(C^*_\Sigma, d^*_\Sigma)$:

$$\Omega^{0,0}(\Sigma_*) \longrightarrow \Omega^{0,1}(\Sigma_*) \oplus \Omega^{0,1}(\Sigma_*, E_*) \longrightarrow \Omega^{0,1}(\Sigma_*, E_*).$$

This tells us that the extended $L^2_{\infty,\delta}$ subcomplex of $(C^*_{sw}, d^*_{sw})$, thought of as a subcomplex of

$$\Omega^{0,0} \longrightarrow \Omega^{0,1} \oplus \Omega^{0,0}(E) \longrightarrow \Omega^{0,2} \oplus \Omega^{0,1}(E) \longrightarrow \Omega^{0,2}(E),$$

consisting of all $\Omega^{0,*}(\Sigma_*)$-extended $L^2_{\infty,\delta}$ forms, carries the cohomology of this complex.

But this latter complex is the mapping cone for the map

$$\underline{\Phi}\colon \Omega^{0,*} \to \Omega^{0,*}(E)$$

induced by

$$\omega \mapsto \omega \otimes \Phi,$$

viewed as a map between complexes of $\Omega^{0,*}(\Sigma_*)$-extended $L^2_{\infty,\delta}$ forms. Thus, by the isomorphisms in Lemmas 9.2.8 and 9.2.9, together with Corollary 9.1.4, we have that

$$H^*(C^*_{sw}, d^*_{sw}) \cong H^*(C^*_{\mathbb{C}}, d^*_{\mathbb{C}}) \cong H^*\left(\frac{\Omega^{0,*}(E)}{\Phi\Omega^{0,*}}\right) \cong H^*(\mathcal{E}/\alpha\mathcal{O}_R).$$

□

## 10. Flows to the Reducible Locus

In this section, we will investigate the analogue of Theorem 7.0.18 for flows to the reducible locus.

Throughout this section we make the following assumptions. Let $Y = S(N)$ be the unit sphere bundle over an orbifold bundle over $\Sigma$ with

$$\deg(N) < 0.$$

Fix a line bundle $E$ over $Y$ which is the pull-back of an orbifold line bundle over $\Sigma$, i.e. $c_1(E)$ is a torsion class. Moreover, assume that the reducible locus $\mathfrak{J}(E)$ in the corresponding $\mathrm{Spin}_c(3)$ structure is non-degenerate. (See Proposition 5.8.4 and Corollary 5.8.5.)

**Definition 10.0.10.** *Let $\lfloor\frac{K}{2}\rfloor$ be the orbifold line bundle over $\Sigma$ of maximal degree among all orbifold line bundles whose degree is less than $\frac{1}{2}\deg(K_\Sigma)$ and whose pull-back is isomorphic to $E$.*

*Moreover, let*

$$\delta_0 = \frac{\deg(\lfloor\frac{K}{2}\rfloor) - \frac{1}{2}\deg(K_\Sigma)}{\deg(N)}.$$

**Remark 10.0.11.** *Given $E \cong \pi^*(E_0)$ for some orbifold line bundle $E_0$ over $\Sigma$, there is a unique line bundle $\lfloor\frac{K}{2}\rfloor$ which satisfies the above properties. This is clear, since the set of all orbifold line bundles over $\Sigma$ whose pull-back is isomorphic to $E$ is the set of all bundles of the form $E_0 \otimes N^k$, for $k \in \mathbb{Z}$; and*

$$\deg(E_0 \otimes N^k) = \deg(E_0) + k\deg(N).$$

*Note that the definition of $\lfloor\frac{K}{2}\rfloor$ depends on the choice of $E$.*

Fix an orbifold line bundle $E_0$ with $\pi^*(E_0) \cong E$, and $\deg(E_0) < \frac{1}{2}\deg(K_\Sigma)$. In fact, fix isomorphisms

$$j_-\colon \pi^*(E_0) \to E$$

and

$$j_+\colon \pi^*(\lfloor\tfrac{K}{2}\rfloor) \to E.$$

Throughout this section, $\widehat{E}$ will denote the bundle over $R$ obtained by gluing $\pi^*(E_0)$ to $\pi^*(\lfloor\frac{K}{2}\rfloor)$ (viewed as bundles over $R - \Sigma_+$ and $R - \Sigma_-$ respectively) using the transition function $e^{\delta_0 t} j_+^{-1} \circ j_-$ over the overlap $R^o \subset R$. Note that by construction the restrictions $\widehat{E}|_{\Sigma_-}$ and $\widehat{E}|_{\Sigma_+}$ are canonically identified with the bundles $E_0$ and $\lfloor\frac{K}{2}\rfloor$ respectively.

**Definition 10.0.12.** *Consider the space of generalized vortex pairs $(A,\alpha)$ which connect vortices in $E_0$ with reducibles in $\pi^*(\lfloor\frac{K}{2}\rfloor)$ via the given isomorphisms $j_\pm$. The* moduli space of generalized vortices which connect $E_0$ to $\mathfrak{J}(E)$, *denoted $\mathcal{M}(E_0,\mathfrak{J}(E))$, is the quotient of this space of pairs by the gauge automorphisms whose derivatives exponentially decay to zero. Dividing the space of generalized vortex pairs by the smaller goup of gauge transformations which exponentially decay to* 1, *we obtain a larger space, the* based moduli space, *which is denoted $\mathcal{M}^\circ(E_0,\mathfrak{J}(E))$. This space is a circle bundle over the (unbased) moduli space.*

There is a corresponding notion for the moduli space of divisors.

**Definition 10.0.13.** *Pick a reference point $x\in R$. The* based moduli space of divisors connecting $E_0$ to $\lfloor\frac{K}{2}\rfloor$, *denoted $\mathcal{D}^\circ(E_0,\lfloor\frac{K}{2}\rfloor)$, is the quotient of the space of holomorphic pairs in $\widehat{E}$ by the subgroup of the complex gauge group which fixes the fiber $\widehat{E}|_x$. The $\mathbb{C}^*$ action on this fiber induces a $\mathbb{C}^*$ action on $\mathcal{D}^\circ(E_0,\lfloor\frac{K}{2}\rfloor)$, whose quotient is $\mathcal{D}(E_0,\lfloor\frac{K}{2}\rfloor)$. Over the points represented by sections which do not vanish identically, the $\mathbb{C}^*$ action on $\mathcal{D}^\circ(E_0,\lfloor\frac{K}{2}\rfloor)$ is free.*

**Remark 10.0.14.** *The space $\mathcal{D}^\circ(E_0,\lfloor\frac{K}{2}\rfloor)$ inherits a deformation theory from $\mathcal{D}(E_0,\lfloor\frac{K}{2}\rfloor)$ (see Remark 7.0.16).*

A gradient flow line in the based moduli space from the critical manifold corresponding to $E_0$, $\mathcal{C}^+(E_0)$, to the reducible locus $\mathfrak{J}(E)$ induces, as in Section 7, a holomorphic pair over $\widehat{E}|_{R-\Sigma_+}$. Indeed, the following analogue of Theorem 7.0.18 tells us that this data naturally extends over all of $\widehat{E}$ over $R$:

**Theorem 10.0.15.** *If $\mathcal{C}^+(E_0)$ is a critical manifold in the* $\mathrm{Spin}_c(3)$ *structure determined by $E$, whose reducible critical manifold $\mathfrak{J}(E)$ is non-degenerate, then there is a natural identification between the based moduli space of flows connecting $\mathcal{C}^+(E_0)$ to $\mathfrak{J}(E)$ with the open subset of based divisors in $\mathcal{D}^\circ(E_0,\lfloor\frac{K}{2}\rfloor)$ which do not contain $\Sigma_-$.*

**Remark 10.0.16.** *In particular, it is an easy consequence of this theorem that this space of flows is non-empty. Given a divisor $D_0\subset\Sigma$ corresponding to a critical point in $\mathcal{C}^+(E_0)$, the divisor*

$$\pi^*(D_0)+\big(\deg(D_0)-\deg(\lfloor\tfrac{K}{2}\rfloor)\big)\,[\Sigma_+]$$

*corresponds to a divisor in $\mathcal{D}(E_0,\lfloor\frac{K}{2}\rfloor)$.*

**Remark 10.0.17.** *The space of flows from $\mathfrak{J}(E)$ to $\mathcal{C}^+(E_0)$ is empty, while all flows from $\mathfrak{J}(E)$ to $\mathfrak{J}(E)$ are stationary. This statement follows easily from Proposition 5.9.5.*

The proof of this theorem follows the same lines as Theorem 7.0.18 and has two parts corresponding to the discussions in Sections 7 and 8 respectively.

To exhibit a map from the $\mathcal{M}^\circ$ to $\mathcal{D}^\circ$ we will use an alternative description of $\mathcal{D}^\circ$. Fix a vortex $(A_-, \alpha_-)$ in $E_0$ and a constant curvature connection $B_0$ in $\lfloor\frac{K}{2}\rfloor$. Consider the space $\mathcal{S}$ of holomorphic pairs $(\widehat{A}, \widehat{\alpha})$ in $\widehat{E}$ so that:

- $(\widehat{A}, \widehat{\alpha})|_{\Sigma_-} = (A_-, \alpha_-)$
- $\widehat{A}|_{\Sigma_+} = B_0$.

The quotient, $\mathcal{Q}$, of this space by the action of the group of complex gauge transformations which are the identity over $\Sigma_\pm$ is naturally identified with the based moduli space of divisors with $[\widehat{A}, \widehat{\alpha}]|_{\Sigma_-} = [A_-, \alpha_-]$ and $[\widehat{A}]|_{\Sigma_+} = [B_0]$ (as complex gauge equivalence classes). We will now obtain an element of $\mathcal{Q}$ from a based flow line connecting the critical point in $\mathcal{C}^+(E_0)$ corresponding to $[A_-, \alpha_-]$ to the reducible corresponding to $\pi^*(B_0)$.

**Lemma 10.0.18.** *The real number $\delta_0 \geq 0$ defined above is the smallest eigenvalue of the operator $i\nabla^B_{\frac{\partial}{\partial\varphi}}$ acting on $\Gamma(Y, S^+)$. Let $\delta_1$ be the smallest non-negative eigenvalue of $\mathbf{D}_B$ belonging to any eigenvector in $\Gamma(Y, S^+) \subset \Gamma(Y, W)$, then $\delta_0 \leq \delta_1$.*

**Proof.** Given $\delta \in \mathbb{R}$, the $\delta$-eigenspace of $i\nabla_{\frac{\partial}{\partial\varphi}}$ for the is the kernel of $i\nabla^{B+i\delta\eta}_{\frac{\partial}{\partial\varphi}}$. Thus, pushing $B+i\delta\eta$ forward (according to Proposition 5.1.3), we identify the $\delta$-eigenspace of $i\nabla^B_{\frac{\partial}{\partial\varphi}}$ with sections of a bundle over $\Sigma$ whose curvature form is

$$F_B + \delta F_N.$$

This proves the first claim. The second follows from the fact that $i\nabla^B_{\frac{\partial}{\partial\varphi}}$ is the $S^+$-component of the Dirac operator, acting on $S^+$. □

**Proposition 10.0.19.** *Let $(A, \alpha)$ be a flow on $\mathbb{R} \times Y$ which converges exponentially to a reducible solution $[B, 0]$ as $t \mapsto +\infty$, then the pair*

$$(e^{t\delta_0}\overline{\partial}_A e^{-t\delta_0}, e^{t\delta_0}\alpha)$$

*converges exponentially to a the pull-back of a holomorphic pair $(B_0, \alpha_0)$ in $\lfloor\frac{K}{2}\rfloor$, where $B_0$ has constant curvature.*

**Proof.** There are constants $\mu_1 > 0$, $\delta_2 > \delta_1 > 0$ such that for any gradient flow line $(A, \alpha)$ as above, there is an asymptotic expansion (for $t \gg 0$) of the form

$$(A, \alpha) = (B + O(e^{-\mu_1 t}), e^{-\delta_1 t}\alpha_1 + O(e^{-\delta_2 t})).$$

Here, $\delta_1$ is the constant from Lemma 10.0.18, so that $\delta_0 \leq \delta_1$. If $\delta_0 < \delta_1$, let $\alpha_0 = 0$, otherwise, let $\alpha_0 = \alpha_1$. In either case, we see that $e^{t\delta_0}\alpha$ converges exponentially to $\alpha_0$.

The section $e^{t\delta_0}\alpha$ is harmonic for the Dirac operator

$$e^{t\delta_0}\mathbf{D}_A e^{-t\delta_0} = \mathbf{D}_A - \delta_0\rho(dt),$$

which is the same as $\mathbf{D}_{A+i\delta_0\eta}$, when restricted to $S^+$. Moreover, the connection $A + i\delta_0\eta$ converges as $t \mapsto +\infty$ to $B + i\delta_0\eta = \pi^*(B_0)$, so we see that $\alpha_0$ corresponds to a holomorphic section of $\lfloor \frac{K}{2} \rfloor$. □

As an immediate consequence, we see that the sections $e^{t\delta_0}\alpha$ over $\pi^*(E_0)$, and $\alpha_0$ over $\Sigma_+$ glue together to give a continuous section $\widehat{\alpha}$ of $\widehat{E}$. We can apply Theorem 7.0.19 to extend the holomorphic structure $\overline{\partial}_{A+i\delta_0}$ over $\pi^*(\lfloor \frac{K}{2} \rfloor)$, giving us a holomorphic structure on all of $\widehat{E}$ with respect to which $\widehat{\alpha}$ is holomorphic. Indeed the complex gauge transformation obtained in Theorem 7.0.19 is the idenity over $\Sigma_\pm$ and hence we obtain an element of $\mathcal{Q}$.

Thus, we have constructed a map from the moduli space of flows to the space of divisors. To invert this map, as before (Section 8), we have to solve a Kazdan-Warner equation for a complex gauge transformation. First, we solve the problem "at infinity":

**Lemma 10.0.20.** *Let $\widehat{E}$ be the line bundle over $R$ defined as above. Suppose that $(\widehat{A}, \widehat{\alpha})$ is a holomorphic pair in $\mathcal{S}$. Then, there is a complex gauge transformation $e^u$ on $E$ which takes the restriction of $(\widehat{A}, \widehat{\alpha})|_{R^o}$ to a pair $(A, \alpha)$ such that*

1. *$(A, \alpha)$ converges exponentially as $t \mapsto -\infty$ to the vortex in $E_0$ with section $\widehat{\alpha}|_{\Sigma_-}$, and*
2. *$(A, \alpha)$ converges exponentially as $t \mapsto +\infty$ to a reducible solution $[B, 0]$.*

**Proof.** The hypothesis that $(\widehat{A}, \widehat{\alpha}) \in \mathcal{S}$ gives us that

$$(2\Lambda F_{\widehat{A}} - \Lambda F_{K_\Sigma} - i|\widehat{\alpha}|^2)\big|_{\Sigma_-} \equiv 0$$

and

$$\left( i(2\Lambda F_{\widehat{A}} - \Lambda F_{K_\Sigma}) - \frac{2\pi(2\deg(\lfloor \frac{K}{2} \rfloor) - \deg(K_\Sigma))}{\mathrm{Vol}(\Sigma)} \right) \big|_{\Sigma_+} \equiv 0.$$

Suppose that $\deg(\lfloor \frac{K}{2} \rfloor) = \frac{1}{2}\deg(K_\Sigma)$, then $\alpha|_{\Sigma_+}$ is a holomorphic section of an orbifold bundle whose degree is $\frac{1}{2}\deg(K_\Sigma)$ and whose pull-back to $Y$ is $E$. Thus, $\alpha|_{\Sigma_+} = 0$, by our assumption that $\mathfrak{J}(E)$ is non-degenerate (see Proposition 5.8.4). Now, the restriction $(A, \alpha)$ of $(\widehat{A}, \widehat{\alpha})$ to $R^o$ satisfies the requirements (by Lemma 8.0.29).

Suppose that $\deg(\lfloor\frac{K}{2}\rfloor)\neq\frac{1}{2}\deg(K_\Sigma)$, so that $\delta_0>0$. Recall that a complex gauge transformation $e^u$ acting on a connection $A$ on $E$ over $\mathbb{R}\times Y$ (i.e. the transformation which takes sections $\alpha\in\Gamma(E,R^o)$ to $e^u\alpha$) transforms the curvature according to

$$i\Lambda F_{e^uA}\mapsto 2i\Lambda\overline{\partial}\partial u+i\Lambda F_A=-e^{-2\xi t}\frac{\partial}{\partial t}e^{2\xi t}\frac{\partial}{\partial t}u+\Delta_Y u+i\Lambda F_A,$$

following Lemma 8.0.34. Thus, if $e^u$ is a complex gauge transformation with

$$u(t)=\begin{cases}0 & \text{if } t\ll 0\\ -\delta_0 t & \text{if } t\gg 0\end{cases},$$

this complex gauge transformation takes the restriction of $(\widehat{A},\widehat{\alpha})$ to a pair $(A,\alpha)$ which exponentially decays as $t\mapsto-\infty$ a vortex in $E_0$, and as $t\mapsto\infty$, it decays to a reducible solution $[B,0]$. □

Thanks to this lemma, finding a flow-line then becomes equivalent to solving the Kazdan-Warner equation

$$\Delta_{cyl}(u)+h(e^u-1)=k \tag{87}$$

for an exponentially decaying function $u$, where $h=|\alpha|^2$ satisfies

1. $h\geq 0$,
2. $\lim_{t\mapsto-\infty}h=h_-$ exponentially, for some function $h_-\not\equiv 0$
3. $\lim_{t\mapsto+\infty}h=0$ exponentially.

Under these hypotheses, Equation 87 has a unique exponentially decaying solution $u$, according to Proposition 8.0.30.

Suppose that $(\widehat{A},\widehat{\alpha})$ is a holomorphic pair as in Lemma 10.0.20. Then, a gauge transformation over $R$ which is the identity over $\Sigma_\pm$ carries the pair $(A,\alpha)$ provided by that lemma to another pair which is equivalent to $(A,\alpha)$ by a gauge transformation which exponentially converges to 1. Thus, thanks to the uniqueness statement in Proposition 8.0.30, our construction maps $\mathcal{Q}$ to $\mathcal{M}^\circ(E_0,\mathfrak{J}(E))$ as promised.

10.1. **Deformation Theory of Flows to the Reducible Locus.** We can slightly modify the discussion in Section 9 to obtain the following analogue of Theorem 9.2.5:

**Theorem 10.1.1.** *Let $E$, $E_0$ be as in Theorem 10.0.15. Then, the based moduli space of flows to the reducible locus $\mathcal{M}^\circ(E_0,\mathfrak{J}(E))$ is locally modeled about any point $(A,\Phi)\in\mathcal{M}^\circ(E_0,\mathfrak{J}(E))$ on the zeros of a map*

$$\mathbb{C}\oplus H^0(D,\mathcal{E}|_D)\to H^1(D,\mathcal{E}|_D),$$

*where $\mathcal{E}$ is the holomorphic structure induced by $(A,\Phi)$ on the line bundle interpolating $E_0$ and $\lfloor\frac{K}{2}\rfloor$ from Theorem 10.0.15, and $D=\alpha^{-1}(0)$ is the zero-set of the associated section. In particular,*

$$\textit{e-dim}\mathcal{M}(\mathcal{C}^+(E_0),\mathfrak{J}(E))=\textit{e-dim}(\mathcal{D}(E_0,\lfloor\tfrac{K}{2}\rfloor))+1. \tag{88}$$

**Proof.** Recall that the tangent space to the reducible locus is given by

$$T_{[B_+,0]}U_2 = \text{Coker}\left(\overline{\partial}\colon \Omega^{0,0}(\Sigma_+) \to \Omega^{0,1}(\Sigma_+)\right).$$

Bearing this in mind, the proof of Theorem 9.2.5 shows that the cohomology groups of the deformation complex are isomorphic to a mapping cone for

$$\underline{\Phi}\colon \Omega^{0,*} \to \Omega^{0,*}(E),$$

where forms in $\Omega^{0,*}$ are extended by $\pi^*(\Omega^{0,*}(\Sigma_-))$ at $-\infty$ but only by $T_{[B_+,0]}U_2$ at $+\infty$; while forms in $\Omega^{0,*}(E)$ are extended by $\pi^*(\Omega^{0,*}(\Sigma_-, E_0))$ at $-\infty$ and by nothing at $+\infty$.

With the above conventions on $\Omega^{0,*}$, the cohomology groups

$$H^*(\Omega^{0,*}, \overline{\partial}) \cong \begin{cases} 0 & \text{if } * = 0 \\ H^*(R, \mathcal{O}_R) & \text{otherwise} \end{cases},$$

which we can see by mapping this complex to the corresponding complex of $L^2_{\infty,\delta}$ forms which are extended by all smooth forms pulling back from $\Sigma_\pm$ at the ends $\pm\infty$, noting that the quotient complex is merely a $\mathbb{C}$ in degree zero, and then applying Corollary 9.1.4.

We wish now to identify the cohomology groups of $(\Omega^{0,*}(E), \overline{\partial}_A)$, given the above decay conventions, with the sheaf cohomology groups $H^*(R, \mathcal{E})$. There are two cases. When $\delta_0 \neq 0$, the connection $A$ used over $R^o$ does not extend over $\Sigma_+$, but rather it satisfies the hypotheses of Corollary 9.1.11, so this corollary provides the requisite identification. When $\delta_0 = 0$, $\Omega^{0,*}(E)$ maps to the complex extended by $\Omega^{0,*}(\Sigma_-, E_0), \Omega^{0,*}(\Sigma_+, \lfloor\frac{K}{2}\rfloor)$, with quotient the complex

$$\Omega^{0,0}(\Sigma_+, \lfloor\tfrac{K}{2}\rfloor) \xrightarrow{\overline{\partial}_B} \Omega^{0,1}(\Sigma_+, \lfloor\tfrac{K}{2}\rfloor).$$

But the cohomology groups of this complex must vanish by the non-degeneracy hypothesis on $\mathfrak{J}(E)$ (see Proposition 5.8.4).

Putting together the above facts, we see that in either case, the cohomology groups

$$H^*(C^*_{sw}, d^*_{sw}) \cong \begin{cases} \mathbb{C} \oplus H^0(\mathcal{E}/\alpha\mathcal{O}_R) & * = 0 \\ H^1(\mathcal{E}/\alpha\mathcal{O}_R) & * = 1 \\ 0 & \text{otherwise} \end{cases}$$

The theorem then follows. □

## 11. Resolving Sheaves over Orbifold Singularities

Given any two relatively prime integers with $0 < q < p$, let $a \in \mathbb{Z}/p\mathbb{Z}$ act on $\mathbb{C}^2$ by

$$a \times (w, z) = (\zeta^a w, \zeta^{qa} z),$$

where $\zeta$ is a primitive $p^{th}$ root of unity. Topologically, the quotient

$$C_{p,q} = \mathbb{C}^2/(\mathbb{Z}/p\mathbb{Z})$$

is a cone on a Lens space of type $(p,q)$.

This quotient has the structure of a complex two-dimensional orbifold. Its (algebraic) coordinate ring $\mathbb{A}[C_{p,q}]$ is given by the subring

$$\mathbb{C}[w,z]^{\mathbb{Z}/p\mathbb{Z}} \subset \mathbb{C}[w,z]$$

of algebraic functions on $\mathbb{C}^2$ which are invariant under the action of $\mathbb{Z}/p\mathbb{Z}$; i.e. those functions $f(w,z)$ for which

$$f(w,z) = f(\zeta w, \zeta^q z).$$

The variety $C_{p,q}$ inherits $p$ distinct orbifold sheaves $\{\mathcal{O}_j\}_{j=0}^{p-1}$, consisting of those functions $f(w,z)$ on $\mathbb{C}^2$ for which

$$f(\zeta w, \zeta^q z) = \zeta^j f(w,z).$$

Equivalently, the representations $\{\rho_j\}_{j=0}^{p-1}$ of $\mathbb{Z}/p\mathbb{Z}$ on the one-dimensional complex vector space $\mathbb{C}$ given by

$$\rho_j(a)(v) = \zeta^{aj} v,$$

along with the above action of $\mathbb{Z}/p\mathbb{Z}$ on $\mathbb{C}^2$, give us $p$ distinct $\mathbb{Z}/p\mathbb{Z}$-equivariant line bundles over $\mathbb{C}^2$, which we will denote by $\{\mathcal{L}_j\}_{j=0}^{p-1}$. Now, the sections of $\mathcal{O}_j$ can be thought of as the $\mathbb{Z}/p\mathbb{Z}$-equivariant sections of the line bundle $\mathcal{L}_j$.

Since the action of $\mathbb{Z}/p\mathbb{Z}$ on $\mathbb{C}^2$ is free away from the origin, the equivariant line bundle $\mathcal{L}_j$ descends to give an honest line bundle on $C_{p,q} - [0,0]$ so that the sheaves $\mathcal{O}_j|_{C_{p,q}-[0,0]}$ are all locally free. Indeed, these $p$ line bundles are all distinct topological bundles over $C_{p,q} - [0,0]$, representing the $p$ elements in

$$H^2(C_{p,q} - [0,0], \mathbb{Z}) \cong H^2(L_{p,q}, \mathbb{Z}) \cong \mathbb{Z}/p\mathbb{Z}.$$

The sheaves $\mathcal{O}_j$ (for $j \neq 0$) do not extend over $[0,0]$ as locally free sheaves (as their first Chern classes do not extend). However, there is a canonical way to extend these sheaves over the resolution $\widehat{C}_{p,q}$. Describing this procedure is the goal of this section.

First, recall that the minimal resolution

$$r\colon \widehat{C}_{p,q} \to C_{p,q}$$

is the map which contracts the chain of two-spheres in $\widehat{C}_{p,q}$, labeled $\{S_i\}_{i=1}^m$. Here, $m$ is the number of terms in the Hirzebruch-Jung continued fraction expansion of $p/q$:

$$\frac{p}{q} = a_1 - \cfrac{1}{a_2 - \cfrac{1}{\ddots - \frac{1}{a_m}}} = \langle a_1, ..., a_m \rangle; \tag{89}$$

and $S_i$ has self-intersection number $-a_i$ (see Section 3). Pulling back the orbifold sheaf $\mathcal{O}_j$ via $r$ (as a coherent sheaf), we obtain a coherent sheaf $r^*(\mathcal{O}_j)$ over $\widehat{C}_{p,q}$.

In order to give a more concrete description of this sheaf, one which is useful in calculations, we must first introduce some terminology. Suppose $\mathbf{w}$ is in a monoid generated by elements $\mathbf{v}_1, ..., \mathbf{v}_m$. Then the $m$-tuple of non-negative integers $(x_1, ..., x_m)$ is called a *maximal decomposition of* $\mathbf{w}$ *in terms of* $(\mathbf{v}_1, ..., \mathbf{v}_m)$ if it is maximal, in the dictionary ordering, among all $m$-tuples $(y_1, ..., y_m)$ of non-negative integers which satisfy

$$\mathbf{w} = \sum_i y_i \mathbf{v}_i.$$

In particular, if $(d_1, ..., d_m)$ is an ordered $m$-tuple of natural numbers, and $j$ is a non-negative integer, then the minimal decomposition of $j$ with respect to $(d_1, ..., d_m)$ (if it exists) can be found inductively by a Euclidean algorithm letting

$$x_1 = \lfloor \frac{j}{d_1} \rfloor,$$

and $(x_2, ..., x_m)$ be the maximal decomposition of $j - x_1 d_1$ in terms of $(d_2, ..., d_m)$.

**Proposition 11.0.2.** *For each* $j = 0, ..., p-1$*, the sheaf* $r^*(\mathcal{O}_j)$ *is an invertible sheaf over* $\widehat{C}_{p,q}$ *for which the natural map*

$$\mathcal{O}_j \to r_* r^*(\mathcal{O}_j).$$

*is an isomorphism. Let* $(x_1, ..., x_m)$ *be a maximal decomposition of* $j$ *in terms of* $(d_1, ..., d_m)$*, where* $d_i$ *is the the denominator of the Hirzebruch-Jung continued fraction* $\langle a_i, a_{i+1}, ..., a_m \rangle$*. Then, the first Chern class of the line bundle corresponding to* $r^*(\mathcal{O}_j)$ *is determined by*

$$\langle c_1(r^*(\mathcal{O}_j)), [S_i] \rangle = x_i,$$

*where* $[S_i]$ *is the homology class of the* $i^{th}$ *sphere in* $\widehat{C}_{p,q}$*.*

Thus, this proposition gives a natural correspondence between orbifold sections of the sheaf $\mathcal{O}_j$ over $C_{p,q}$ and sections of a line bundle (depending on $j$) over $\widehat{C}_{p,q}$. The result should be compared with the following theorem of Gonzalez-Sprinberg and Verdier [12] (see also [19]):

**Theorem 6** (Gonzalez-Sprinberg, Verdier)**.** *Let* $G$ *be a finite subgroup of* $SU(2)$*. The resolution* $\mathbb{C}^2/G$ *is obtained by inserting a configuration of* $2-$*spheres which are in a one-to-one correspondence with the non-trivial irreducible representations of* $G$*. Given a representation*

$$\rho\colon G \to GL(V),$$

*letting* $\mathcal{O}_\rho$ *denote the sheaf of orbifold sections of of the orbifold bundle*

$$\mathbb{C}^2 \times_{\{G,\rho\}} V$$

*over* $\mathbb{C}^2/G$*, we can form*

$$r(\mathcal{O}_\rho) = r^*(\mathcal{O}_\rho)/\mathrm{Tors},$$

*the quotient of the pull-back sheaf by its nilpotent elements. Then,* $r(\mathcal{O}_\rho)$ *is a locally free sheaf for which the natural map*

$$\mathcal{O}_\rho \to r_*(r(\mathcal{O}_\rho))$$

*is an isomorphism and*

$$\langle c_1(r(\rho)), [S_\sigma]\rangle = \begin{cases} 1 & \text{if } \sigma = \rho \\ 0 & \text{otherwise} \end{cases},$$

*where $S_\sigma$ denotes the sphere corresponding to the representation $\sigma$.*

The hypotheses of our Proposition are consistent with the hypotheses of the Gonzalez-Sprinberg and Verdier theorem exactly when $q = p - 1$, i.e. when $G \cong \mathbb{Z}/p\mathbb{Z}$ is given by matrices of the form

$$\begin{pmatrix} \zeta^a & 0 \\ 0 & \zeta^{-a} \end{pmatrix},$$

with $a \in \mathbb{Z}/p\mathbb{Z}$. There are $p-1$ terms in the continued fraction expansion of $p/p-1$, and the denominators are $(p-1, p-2, ..., 1)$, so that the maximal expansion of any $0 < j < p$ is given by coefficients $x_i = \delta_{i,p-j}$ (Kronecker delta). Thus, in this special case, both Theorem 6 and Proposition 11.0.2 are in agreement.

The proof of Proposition 11.0.2 requires a concrete description of $\widehat{C}_{p,q}$, which is based on elements of the theory of toric varieties (see [9]). We set up the notation presently.

The map

$$\chi\colon \mathbb{Z}^2 \mapsto \mathbb{C}[w, z, w^{-1}, z^{-1}]$$

from the standard lattice to the Laurent polynomials

$$\chi(i, j) = w^i z^j$$

sets up a correspondence between additive sub-monoids of $\mathbb{Z}^2$ and subrings of $\mathbb{C}[w, z, w^{-1}, z^{-1}]$. Under this correspondence, the lattice $\Lambda \subset \mathbb{Z}^2$ given by

$$i + qj \equiv 0 \pmod p$$

corresponds to the subring of $\mathbb{Z}/p\mathbb{Z}$-invariant Laurent polynomials. Moreover, the upper right quadrant of $\Lambda$,

$$\Lambda_+ = \{v \in \Lambda \big| \Pi_1(v) \geq 0, \Pi_2(v) \geq 0\},$$

corresponds to the coordinate ring $\mathbb{A}[C_{p,q}]$. In the above expression,

$$\Pi_\alpha\colon \mathbb{Z}^2 \to \mathbb{Z}$$

denotes the projection onto the $\alpha^{th}$ coordinate.

As in [9], the monoid $\Lambda_-$

$$\Lambda_- = \{v \in \Lambda \big| \Pi_1(v) \leq 0, \Pi_2(v) \geq 0\}$$

has a canonical generating set $(\mathbf{v}_0, ..., \mathbf{v}_{m+1})$, arranged so that

$$-p = \Pi_1\mathbf{v}_0 < \Pi_1\mathbf{v}_1 < ... < \Pi_1\mathbf{v}_m < \Pi_1\mathbf{v}_{m+1} = 0. \tag{90}$$

Consecutive pairs of vectors $\mathbf{v}_i$, $\mathbf{v}_{i+1}$ generate the monoid of lattice points in the sector of $\Lambda$ lying between those vectors. The vectors $\mathbf{v}_0 = (-p, 0)$, $\mathbf{v}_{m+1} = (0, p)$, and the vectors $\{\mathbf{v}_1, ..., \mathbf{v}_m\}$ are the lattice points appearing on the boundary of the convex hull of the set $\Lambda_{_} - \{(0,0)\}$, off the coordinate axes. This convexity forces

$$0 = \Pi_2\mathbf{v}_0 < \Pi_2\mathbf{v}_1 < ... < \Pi_2\mathbf{v}_m < \Pi_2\mathbf{v}_{m+1} = p. \tag{91}$$

These vectors satisfy the relations

$$\mathbf{v}_{i-1} + \mathbf{v}_{i+1} = a_i\mathbf{v}_i, \tag{92}$$

where the $a_i$ appear in the continued fraction expansion in Equation 89. Moreover,

$$\Pi_1\mathbf{v}_i = -d_i, \tag{93}$$

where the $d_0 = p$, $d_{m+1} = 0$, and $d_i$ for $i = 1, ..., m$ are as in the statement of the Proposition 11.0.2. (Similarly, for $i = 1, ..., m$, $\Pi_2\mathbf{v}_i$ is the denominator of the Hirzebruch-Jung continued fraction expansion of $\langle a_i, ..., a_1\rangle$.) Equations 92, and 89 (along with the statement about $\Pi_2\mathbf{v}_i$) are neatly summarized by the statement that, for each $i = 1, ..., m$, there is a matrix $A_i \in \mathrm{Sl}_2(\mathbf{Z})$ with $A_i\mathbf{v}_{i-1} = \mathbf{v}_i$, and $A_i\mathbf{v}_i = \mathbf{v}_{i+1}$.

The ordering properties expressed in Inequalities 90 and 91 give us a constructive technique for giving the maximal decomposition of a vector $\mathbf{w} \in \Lambda_-$.

**Lemma 11.0.3.** *Suppose $\mathbf{w} \in \Lambda_-$, and $(x_0, ..., x_m)$ is the maximal decomposition of $\Pi_1\mathbf{w}$ with respect to $(d_0, ..., d_m)$. Then, letting*

$$x_{m+1} = \frac{\Pi_2(\mathbf{w} - \sum_{i=0}^m x_i\mathbf{v}_i)}{p},$$

*the $(x_0, ..., x_{m+1})$ form the maximal decomposition of the $\mathbf{w}$ with respect to the generators $(\mathbf{v}_0, ..., \mathbf{v}_{m+1})$. In particular, we have for the maximal decomposition*

$$\mathbf{v} = \sum_i x_i\mathbf{v}_i$$

*that*

$$d_i > \sum_{j>i} x_j d_j. \tag{94}$$

**Proof.** Letting $\mathbf{u} = \sum_{i=0}^m x_i\mathbf{v}_i$, $\mathbf{u} \in \Lambda_-$ has $\Pi_1\mathbf{u} = \Pi_1\mathbf{w}$ by Equation 93. The maximality of the $x_i$, combined with the ordering property in Inequality 91 ensure that $\mathbf{u}$ has minimal $\Pi_2$ among all vectors in $\Lambda \cap {\Pi_1}^{-1}\Pi_1\mathbf{w}$. This implies the first statement.

The inequality 94 follows immediately from the construction of maximal decompositions of natural numbers. □

If $\mathbf{v}, \mathbf{w} \in \Lambda_-$, we say that $\mathbf{v}$ precedes $\mathbf{w}$, denoted $\mathbf{v} \prec \mathbf{w}$, if the maximal decomposition of $\mathbf{v}$ with respect to the generators $(\mathbf{v}_0, ..., \mathbf{v}_m)$ precedes that of $\mathbf{w}$ in the dictionary ordering.

**Corollary 11.0.4.** *If $\mathbf{v}, \mathbf{w} \in \Lambda_-$ and $\Pi_1 \mathbf{w} < \Pi_1 \mathbf{v}$ then there is a $\lambda \in \Lambda_+$ and $\mathbf{u} \in \Lambda_-$ with $\mathbf{u} \preceq \mathbf{w}$ such that*

$$\mathbf{v} = \mathbf{u} + \lambda.$$

**Proof.** Let $\mathbf{v} = \sum b_i \mathbf{v}_i$, $\mathbf{w} = \sum c_i \mathbf{v}_i$ be the maximal decompositions. We might as well assume $\mathbf{v} \not\preceq \mathbf{w}$ (otherwise the conclusion is satisfied for $\lambda = 0$). Thus, there must be an integer $k$ with $b_k > c_k$. Let $\ell$ be the integer with the property that for all $j < \ell$, $b_j = c_j$, and $b_\ell \neq c_\ell$.

The condition that $\Pi_1 \mathbf{w} < \Pi_1 \mathbf{v}$, along with Inequality 94 implies that $b_\ell < c_\ell$. Let

$$\mathbf{u} = \sum_{j<\ell} c_j \mathbf{v}_j + (b_\ell + 1)\mathbf{v}_\ell.$$

Clearly, $\mathbf{u} \preceq \mathbf{w}$ as required. Moreover, by Inequality 94 once again, we have that $\Pi_1 \mathbf{u} < \Pi_1 \mathbf{v}$. The existence of some $k$ with $b_k > c_k$, along with the ordering property Inequality 91 forces $\Pi_2(\mathbf{u}) < \Pi_2(\mathbf{v})$. Thus, $\mathbf{v} - \mathbf{u} \in \Lambda_+$. □

The space $\widehat{C}_{p,q}$ is covered by $m + 1$ affine coordinate patches $\{U_i\}_{i=0}^m$, where

$$U_i = \text{Spec } \mathbb{C}[\chi(-\mathbf{v}_i), \chi(\mathbf{v}_{i+1})]$$

corresponds to the submonoid $\Lambda_i = \langle -\mathbf{v}_i, \mathbf{v}_{i+1}\rangle$ of $\Lambda$ generated by $-\mathbf{v}_i$ and $\mathbf{v}_{i+1}$. Moreover, the set

$$U_{i-1} \cup U_i$$

is the total space of a line bundle over

$$S_i = \text{Spec } \mathbb{C}[\chi(-\mathbf{v}_i)] \cup \text{Spec } \mathbb{C}[\chi(\mathbf{v}_i)],$$

which is a projective line. Its transition function can be read off from the relation in Equation 92, so we see that $c_1$ of the line bundle is $-a_i$.

The ring of global functions on this variety is given by

$$\bigcap_i \mathbb{C}[\chi(-\mathbf{v}_i), \chi(\mathbf{v}_{i+1})] = \mathbb{C}[\Lambda_+],$$

inducing the map

$$r \colon \widehat{C}_{p,q} \to C_{p,q}.$$

A module over the monoid $\Lambda_+$ naturally induces a module over the ring $\mathbb{C}[\Lambda_+]$. We claim that for all $j$, the modules $\Gamma(C_{p,q}, \mathcal{O}_j)$ over the coordinate ring $\mathbb{C}[\Lambda_+]$ arise in this way.

**Lemma 11.0.5.** *The module*

$$M_j = [-j, +\infty) \times [0, \infty) \cap \Lambda \subset \Lambda$$

*over* $\Lambda_+$ *induces a module over* $\mathbb{C}[\Lambda_+]$ *which is isomorphic to* $\Gamma(C_{p,q}, \mathcal{O}_j)$.

**Proof.** From its description, $\Gamma(C_{p,q}, \mathcal{O}_j)$ is an $\mathbb{A}[C_{p,q}]$-submodule of of the module $\mathbb{C}[w, z]$. Multiplication by $z^{-j}$ is an $\mathbb{A}[C_{p,q}]$ automorphism of $\mathbb{C}[w, z, w^{-1}, z^{-1}]$, taking the module $\Gamma(C_{p,q}, \mathcal{O}_j)$ to the $\mathbb{C}[\Lambda_+]$-module $M_j$. □

Let $\mathbf{w}_j$ be the unique vector in $\Lambda$ with

$$\Pi_1 \mathbf{w}_j = -j \text{ and } 0 \le \Pi_2 \mathbf{w}_j \le p.$$

**Lemma 11.0.6.** *The module* $M_j$ *over the monoid* $\Lambda_+$ *is generated by all vectors* $\mathbf{u} \in \Lambda_-$ *with* $\mathbf{u} \preceq \mathbf{w}_j$.

**Proof.** By Inequality 90, any $\mathbf{u}$ with $\mathbf{u} \preceq \mathbf{w}_j$ must be contained in $M_j$. This fact, together with Corollary 11.0.4 implies the lemma. □

Since $\Lambda_+$ is a submonoid of $\Lambda_i = \langle -v_i, v_{i+1} \rangle$, we can consider the monoid $\Lambda_i M_j$, the $\Lambda_i$-submodule of $\Lambda$ generated by $M_j$.

**Lemma 11.0.7.** *The* $\Lambda_i$ *module* $\Lambda_i M_j \subset \Lambda$ *is free and generated by the single element* $\sum_{k \le i} x_k \mathbf{v}_k$. *Equivalently, the* $\mathbb{C}[\Lambda_i]$*-module* $M_j \otimes_{\mathbb{C}[\Lambda_+]} \mathbb{C}[\Lambda_i]$ *is free and generated by the single element* $\chi(\sum_{k \le i} x_k \mathbf{v}_k)$.

**Proof.** We must verify that the map

$$\Lambda_i \to \Lambda_i M_j$$

taking

$$\lambda \mapsto \lambda + \sum_{k \le i} x_k \mathbf{v}_k$$

is surjective (as it is clearly injective). By the ordering of the $\mathbf{v}_k$ (Inequality 90, 91), we have that

$$-\mathbf{v}_k \in \Lambda_i \quad \text{if} \quad k \le i,$$
$$\mathbf{v}_k \in \Lambda_i \quad \text{if} \quad k > i.$$

By Lemma 11.0.6, it suffices to show that any vector $\mathbf{u} \preceq \mathbf{w_j}$ lies in the image. Writing

$$\mathbf{u} = \sum_k y_k \mathbf{v}_k,$$

where

$$y_k \le x_k,$$

we can write

$$\mathbf{u} = \sum_{k\leq i} x_k \mathbf{v}_k + \sum_{k>i} y_k \mathbf{v}_k + \sum_{\ell\leq i} (x_k - y_k)(-\mathbf{v}_k),$$

where the second two terms are in $\Lambda_i$, by the ordering property. □

**Remark 11.0.8.** *Recall (see for example* [16]*) that under the composition*

$$U_i \subset \widehat{C}_{p,q},$$

*the restriction of* $r^*(\mathcal{O}_j)$ *corresponds to the module* $M_j \otimes_{\mathbb{C}[\Lambda_+]} \mathbb{C}[\Lambda_i]$.

**Lemma 11.0.9.** *The module* $M_j$ *is precisely the intersection of the* $\Lambda_i M_j$*;*

$$M_j = \bigcap_{i=0}^{m} \Lambda_i M_j.$$

**Proof.** Clearly,

$$M_j \subset \bigcap_{i=0}^{m} \Lambda_i M_j.$$

Given $\mathbf{v} \in \bigcap_{i=0}^{m} \Lambda_i M_j$, as $\Pi_2 \lambda \geq 0$ for any $\lambda \in \Lambda_0$, we must have

$$\Pi_2 \mathbf{v} \geq 0.$$

Similarly, since $\Pi_1 \lambda \geq 0$ for any $\lambda \in \Lambda_m$, we have that

$$\Pi_1 \mathbf{v} \geq -j.$$

These conditions force $\mathbf{v} \in M_j$. □

Proposition 11.0.2 follows rather easily from these lemmas.

**Proof.** [of Proposition 11.0.2] We have exhibited an affine cover, the

$$\{\text{Spec } \mathbb{C}[\Lambda_i]\}_{i=0}^{r},$$

of $\widehat{C}_{p,q}$ over which $r^*(\mathcal{O}_j)$ is free by Lemma 11.0.7. The statement that $\mathcal{O}_i \mapsto r_* r^*(\mathcal{O}_i)$ follows from Lemma 11.0.9.

Restricting to $S_i$, we have that the generator of $r^*(\mathcal{O}_j)$ from Lemma 11.0.7 over $U_i \cap S_i$ differs from the generator of $r^*(\mathcal{O}_j)$ over $U_{i+1} \cap S_i$ by multiplication by $\chi(\mathbf{v}_i)^{x_i}$. Thus, it follows immediately that

$$\langle c_1(r^*(\mathcal{O}_j)), S_i \rangle = x_i.$$

□

Suppose now that $R$ is a complex orbifold with finitely many $\mathbb{C}^2/(\mathbb{Z}/p\mathbb{Z})$ singular points, and $\mathcal{F}$ is a reflexive sheaf over $R$ which is locally free away from the singular points, and is isomorphic to $\mathcal{O}_{j(x)}$ near the singular point $x \in R$. Consider the minimal resolution

$$r\colon \widehat{R} \to R,$$

obtained by locally performing the procedure discussed above. We obtain the following precise statement of Theorem 3 from Section 1:

**Theorem 7.** *The sheaf $r^*(\mathcal{F})$ is a locally free sheaf over $\widehat{R}$, with Chern classes determined as in Proposition 11.0.2. Moreover, $r$ induces an isomorphism between the moduli space of divisors in $r^*(\mathcal{F})$ over $\widehat{R}$ with the moduli space of divisors in $\mathcal{F}$ over $R$.*

**Proof.** The fact that $r^*(\mathcal{F})$ is locally free follows directly from Proposition 11.0.2, after passing to germs near the singularity. Similarly, we have that the natural map

$$\mathcal{F} \to r_* r^* \mathcal{F} \tag{95}$$

is an isomorphism. This gives us the fact that the map
between moduli spaces is, as a map of sets, bijective.

To check that the local models agree, take any nonzero section

$$s \in H^0(\widehat{R}, r^*\mathcal{F}) - 0.$$

This section induces a short exact sequence of sheaves

$$0 \longrightarrow \mathcal{O}_{\widehat{R}} \xrightarrow{\ s\ } r^*\mathcal{F} \longrightarrow r^*\mathcal{F}_{s^{-1}(0)} \longrightarrow 0,$$

whose direct image (via the rationality of the singularities, and the fact that Natural Map 95 is an isomorphism) is the short exact sequence

$$0 \longrightarrow \mathcal{O}_R \xrightarrow{\ r(s)\ } \mathcal{F} \longrightarrow \mathcal{F}_{rs^{-1}(0)} \longrightarrow 0.$$

Rationality of the singularities implies that all

$$\mathrm{R}^j r_*(\mathcal{O}_{\widehat{R}}) = 0$$

for $j > 0$, so by the Leray spectral sequence

$$H^i(R, \mathrm{R}^j r_*(\mathcal{O}_{\widehat{R}})) \Rightarrow H^{i+j}(\widehat{R}, \mathcal{O}_{\widehat{R}})$$

(see [11]), we have natural isomorphisms

$$H^i(\widehat{R}, \mathcal{O}_{\widehat{R}}) \cong H^i(R, \mathcal{O}_R).$$

Similarly, we have natural identifications between

$$H^i(\widehat{R}, r^*\mathcal{F}) \cong H^i(R, \mathcal{F}),$$

because of the vanishing

$$\mathrm{R}^j r_*(\mathcal{O}_{\widehat{R}}) = 0$$

for $i > 0$, which we see as follows. The fibers of $r$ have (complex) dimension at most 1, so the $\mathrm{R}^j r_*$ must vanish for all $j > 1$, by a standard result of sheaf theory (see [11]). The fact that

$$\mathrm{R}^1 r_*(r^*(\mathcal{F})) = 0$$

follows from Grauert's comparison theorem (see for example [13]), together with the fact that

$$H^1(\widehat{r^{-1}(x)}, r^*(\mathcal{F})|_{r^{-1}(x)}) = 0$$

for any $x \in R$. This is true because the fiber is zero-dimensional over each non-singular point, and over each singular point, it is a wedge of projective lines, and the sheaf restricts over each line is a locally free sheaf with non-negative first Chern class (as in the Proposition 11.0.2), so its $H^1$ must vanish. Alternatively, the vanishing of $\mathrm{R}^1 r_*$ follows from a more general argument in [12] (see also [1]). □

## 12. The Dimension Formula

In Section 11, we gave a correspondence between divisors in the ruled surface $R$ with divisors in its desingularization $\widehat{R}$. The purpose of this section is to use this correspondence, together with the Riemann-Roch formula, to derive the explicit formula from Section 1 (Corollary 1.0.4) for the expected dimension of the space of divisors in $R$ in terms of Seifert data.

Let $E_1, E_2$ be a pair of orbifold line bundles over $\Sigma$. Theorem 7 gives an identification between the moduli space of effective divisors connecting $E_1$ with $E_2$, $\mathcal{D}(E_1, E_2)$, with the moduli space of divisors in a line bundle we will denote by $\widehat{E}_{1,2}$ over $\widehat{R}$. Letting $\mathcal{D}(\widehat{E}_{1,2})$ denote this latter moduli space, the ordinary Riemann-Roch theorem gives us that

$$\text{e-dim}\mathcal{D}(E_1, E_2) = \text{e-dim}\mathcal{D}(\widehat{E}_{1,2}) = \langle c_1(\widehat{E}_{1,2})^2 + c_1(K_{\widehat{R}})c_1(\widehat{E}_{1,2}), [\widehat{R}]\rangle.$$

Since $\deg(Y) \neq 0$, the intersection form of $\widehat{R}$ splits (over $\mathbb{Q}$) into the summands corresponding to the resolutions of the line bundles $Y \times_{S^1} \mathbb{C}$ and $Y^* \times_{S^1} \mathbb{C}$, so we can naturally decompose

$$\text{e-dim}\mathcal{D}(E_1, E_2) = \dim_Y(E_1) + \dim_{Y^*}(E_2),$$

where $\dim_Y(E)$ for any Seifert fibered space $Y$ and orbifold line bundle $E$ over $\Sigma$ is given by the formula

$$\dim_Y(E) = \langle c_1(\widehat{E})^2 + c_1(K_X)c_1(\widehat{E}), [\widehat{X}]\rangle, \tag{96}$$

where $\widehat{X}$ is the smooth, non-compact complex surface obtained as the minimal resolution of $X = Y \times_{S^1} \mathbb{C}$, and $\widehat{E}$ is the desingularization of $E$, as given by the procedure of Section 11.

We will now derive a concrete formula for $\dim_Y(E)$ in terms of Seifert data. Let $g$ denote the genus of $\Sigma$, and suppose that $\Sigma$ has $n$ singularities with multiplicities $(\alpha_1, ..., \alpha_n)$. Let $(b, \beta_1, ..., \beta_n)$ be the Seifert invariants of $Y$. Let $\widehat{\Sigma} \subset \widehat{X}$ denote the proper transform of the zero-section $\Sigma \subset X$, and for $i = 1, ..., n$, $j = 1, ..., m_i$, let $S_i^j$ denote the $j^{th}$ sphere in the configuration of spheres lying over the $i^{th}$ singular point in $\Sigma \subset X$. The homology classes induced by the curves

$$\widehat{\Sigma}, S_1^1, ...S_{m_1}^1, S_1^2, ..., S_{m_2}^2, ..., S_1^n, ..., S_{m_n}^n$$

form a basis $\mathfrak{B}$ for $H_2(\widehat{X})$. With respect to the Poincaré dual basis $\mathfrak{B}^*$, the intersection form for $\widehat{X}$ can be written as a sparse matrix $M$ (with zeros in all the entries not marked):

$$\begin{pmatrix}
b & 1 & & & & 1 & & & & \cdots & 1 & & & \\
1 & -a_1^1 & 1 & & & & & & & & & & & \\
 & 1 & -a_2^1 & \ddots & & & & & & & & & & \\
 & & \ddots & \ddots & 1 & & & & & & & & & \\
 & & & 1 & -a_{m_1}^1 & & & & & & & & & \\
1 & & & & & -a_1^2 & 1 & & & & & & & \\
 & & & & & 1 & -a_2^2 & \ddots & & & & & & \\
 & & & & & & \ddots & \ddots & 1 & & & & & \\
 & & & & & & & 1 & -a_{m_2}^2 & & & & & \\
\vdots & & & & & & & & & \ddots & & & & \\
1 & & & & & & & & & & -a_1^n & 1 & & \\
 & & & & & & & & & & 1 & -a_2^n & \ddots & \\
 & & & & & & & & & & & \ddots & \ddots & 1 \\
 & & & & & & & & & & & & 1 & -a_{m_n}^n
\end{pmatrix},$$

where the $a_1^j, ..., a_{m_j}^j$ are the coefficients appearing in the Hirzebruch-Jung continued fraction expansion for $\alpha_j/\beta_j$.

Let $d_i^j$ be given by the denominator Hirzebruch-Jung continued fraction

$$\frac{d_i^j}{d_{i+1}^j} = \langle a_i^j, ..., a_{m_j}^j \rangle,$$

so that

$$a^j_i = \frac{d^j_{i-1}}{d^j_i} + \frac{d^j_{i+1}}{d^j_i}. \tag{97}$$

Let $E$ be an orbifold bundle over $\Sigma$ with Seifert invariants $(e_0, \epsilon_1, ..., \epsilon_n)$, and let $\xi^j_i$ be the coefficients appearing in the minimal expansion of $\epsilon_j$ with respect to $(d^j_1, ..., d^j_{m_j})$,

$$\epsilon_j = \sum_{i=1}^{m_j} \xi^j_i d^j_i. \tag{98}$$

It follows from Section 11, that the first Chern class of $\widehat{E}$ over $\widehat{X}$, the resolution of $\pi^*(E)$, is uniquely characterized by its evaluations

$$\langle c_1(\widehat{E}), [\widehat{\Sigma}]\rangle = e_0,$$

$$\langle c_1(\widehat{E}), [S^j_i]\rangle = \xi^j_i.$$

Thus, if we let $\Xi$ be the vector formed by concatenating $e_0$, $\xi^1_1, ..., \xi^1_{m_1}$, $\xi^2_1, ..., \xi^2_{m_1}$, ..., $\xi^n_1, ..., \xi^n_{m_n}$, we can write $c_1(\widehat{E})$, in the basis $\mathfrak{B}^*$ as the vector

$$\mathbf{x} = M^{-1}\Xi.$$

We will write this vector more explicitly, in the following Lemma.

**Lemma 12.0.10.** *Let*

$$\mathbf{x} = (x^0, x^1_1, ..., x^1_{m_1}, x^2_1, ..., x^2_{m_1}, ..., x^n_1, ..., x^n_{m_n})$$

*denote the coefficients of* $c_1(\widehat{E})$ *in terms of the basis* $\mathfrak{B}^*$ *for* $H^2(\widehat{X})$. *Then,*

$$x^0 = \frac{\deg(E)}{\deg(Y)}. \tag{99}$$

*and*

$$x^j_\ell = d^j_\ell \left( \frac{x^0}{d^j_0} - \sum_{i=1}^{\ell} \frac{1}{d^j_{i-1} d^j_i} \sum_{k=i}^{m_j} d^j_k \xi^j_k \right). \tag{100}$$

**Remark 12.0.11.** *In the above formulas, we have let* $d^j_0 = \alpha_j$.

**Proof.** The equation $M\mathbf{x} = \Xi$ can be written as a system of equations:

$$bx^0 + x^1_1 + x^2_1 + ...x^n_1 = e_0,$$

and, for $j = 1, ..., n, \ell = 1, ..., m_j - 1$

$$x^j_{\ell-1} - a^j_\ell x^j_\ell + x^j_{\ell+1} = \xi^j_\ell,$$

$$x^j_{m_j-1} - a^j_{m_j} x^j_{m_j} = \xi^j_{m_j},$$

where we declare $x_0^j = x^0$ for $j = 1, ..., n$. This is equivalent to the system of equations:

$$\begin{aligned} (b + \frac{d_1^1}{d_0^1} + ... + \frac{d_1^n}{d_0^n})x^0 &= e_0 + \sum_{k=1}^{m_1} d_k^1 \xi_k^1 + ... + \sum_{k=1}^{m_n} d_k^n \xi_k^n \\ x_{\ell-1}^j - \frac{d_{\ell-1}^j}{d_\ell^j} x_\ell^j &= \frac{1}{d_\ell^j} \sum_{k=\ell}^{m_j} d_k^j \xi_k^j. \end{aligned}$$

By Equation 98, the first equation is equivalent to Equation 99. The other can be reexpressed as Equation 100. □

The canonical class can be described (via the adjunction formula) as the cohomology class with

$$\langle c_1(K), [\widehat{\Sigma}] \rangle = -b + 2g - 2,$$
$$\langle c_1(K), [S_i^j] \rangle = a_i^j - 2,$$

so, letting $\kappa$ denote the vector obtained by concatenating the above values, the vector $M^{-1}\kappa$ represents the coefficients of $c_1(K)$ in our basis. Now, we can rewrite Equation 96 as:

$$(\Xi - \kappa)^\dagger M^{-1} \Xi.$$

We can combine this formula with Equation 97 and the results of Lemma 12.0.10 to get that the dimension is given by

$$\begin{aligned} \dim_Y(E) &= \sum_{j=1}^{n} \sum_{\ell=1}^{m_j} (\frac{x^0}{d^0} - \sum_{i=1}^{\ell} \frac{1}{d_{i-1}^j d_i^j} \sum_{k=i}^{m_j} d_k^j \xi_k^j)(d_\ell^j \xi_\ell^j - d_{\ell-1}^j - d_{\ell+1}^j + 2d_\ell^j) \\ &\quad + x^0(e_0 + b + 2 - 2g). \end{aligned} \tag{101}$$

We will rewrite this formula after a few observations.

First, notice that for a fixed $j$,

$$\begin{aligned} \sum_{\ell=1}^{m_j} (d_\ell^j \xi_\ell^j - d_{\ell-1}^j - d_{\ell+1}^j + 2d_\ell^j) &= \left( \sum_{\ell=1}^{m_j} d_\ell^j \xi_\ell^j \right) + (d_1^j - d_0^j) - (d_{m_j}^j) \\ &= \epsilon_j + d_1^j - d_0^j + 1. \end{aligned}$$

Thus, summing this over $j$, we see that

$$\begin{aligned} &\sum_{j=1}^{n} \sum_{\ell=1}^{m_j} (d_\ell^j \xi_\ell^j - d_{\ell-1}^j - d_{\ell+1}^j + 2d_\ell^j) + (e_0 + b + 2 - 2g) \\ &= \sum_{j=1}^{n} (\frac{\epsilon_j}{d_0^j} + \frac{d_1^j}{d_0^j} - 1 + \frac{1}{d_0^j}) + (e_0 + b + 2 - 2g) \\ &= \deg(E) + \deg(Y) - \deg(K_\Sigma) \end{aligned} \tag{102}$$

Moreover, by the discrete analogue of integration by parts, applied twice,

$$\begin{aligned}
&\sum_{\ell=1}^{m_j}(\sum_{i=1}^{\ell}\frac{1}{d^j_{i-1}d^j_i}\sum_{k=i}^{m_j}d^j_k\xi^j_k)((d^j_{\ell+1}-d^j_\ell)-(d^j_\ell-d^j_{\ell-1}))\\
&= (d^j_{m+1}-d^j_m)(\sum_{i=1}^{m_j}\frac{1}{d^j_{i-1}d^j_i}\sum_{k=i}^{m_j}d^j_k\xi^j_k)+\sum_{\ell=1}^{m_j}(\frac{1}{d^j_\ell}-\frac{1}{d^j_{\ell-1}})\sum_{k=\ell}^{m_j}d^j_k\xi^j_k\\
&= (\sum_{i=1}^{m_j}\frac{1}{d^j_{i-1}d^j_i}\sum_{k=i}^{m_j}d^j_k\xi^j_k)-\frac{1}{d^j_0}\sum_{\ell=1}^{m_j}d^j_\ell\xi^j_k+\sum_{\ell=1}^{m_j}\xi^j_\ell\\
&= (\sum_{i=1}^{m_j}\frac{1}{d^j_{i-1}d^j_i}\sum_{k=i}^{m_j}d^j_k\xi^j_k)-\frac{\epsilon_j}{d^j_0}+\sum_{\ell=1}^{m_j}\xi^j_\ell \qquad (103)
\end{aligned}$$

Substituting Equations 102, 103, and 99 back into Equation 101, we get:

$$\begin{aligned}
\dim_Y(E) &= \sum_{j=1}^{n}(\sum_{\ell,k=1}^{m_j}d^j_k\xi^j_kd^j_\ell\xi^j_\ell\sum_{i=1}^{\min(k,\ell)}\frac{1}{d^j_{i-1}d^j_i}+\sum_{\ell=1}^{m_j}d^j_\ell\xi^j_\ell\sum_{i=1}^{\ell}\frac{1}{d^j_{i-1}d^j_i}-\sum_{\ell=1}^{m_j}\xi^j_\ell)\\
&\quad +e_0+\frac{\deg(E)}{\deg(Y)}(\deg(E)-\deg(K_\Sigma)). \qquad (104)
\end{aligned}$$

## 13. Examples

In this section, we give some examples of the theory. We restrict for simplicity to the case of rational homology spheres, i.e. those Seifert fibrations

$$\pi\colon Y\to\Sigma$$

which have non-zero degree, and whose base orbifold has genus zero. These are the manifolds for which the reducible locus consists of isolated points.

If the base orbifold has less than three singular points (so that the total space $Y$ is a lens space), then the Euler characteristic of the base orbifold is positive, so that there are no irreducible critical points. This agrees with the calculations which can be made for the critical points of the equations using the Levi-Civita connection for a metric of positive curvature (see [34]).

**Definition 13.0.12.** *An orbifold is called* simply-connected *if its genus is zero and its multiplicities are pairwise relatively prime.*

**Remark 13.0.13.** *A more natural definition is to say that an orbifold is simply-connected if it admits no non-trivial, connected orbifold covering spaces. This condition is equivalent to the above definition, according to the general theory of orbifolds (see* [10], [30]*).*

For simply-connected orbifolds, the isomorphism class of an orbifold line bundle is uniquely specified by its degree. According to Corollary 5.8.5, if

$$\pi\colon Y \to \Sigma$$

is a Seifert fibered space over a simply-connected orbifold $\Sigma$, the the reducible locus for any $\mathrm{Spin}_c(3)$ structure is smooth. According to Theorem 2.0.19, these simply-connected orbifolds are precisely those which arise as the base spaces of integral homology Seifert fibered spaces. In fact, given a simply-connected orbifold $\Sigma$, there are exactly two Seifert fibered spaces over $\Sigma$ whose underlying three-manifold is an integral homology sphere. The two fibrations correspond to the two generators of the topological Picard group of $\Sigma$, and the two underlying three-manifolds are naturally (orientation-reversing) diffeomorphic. Since $\Sigma$ is uniquely specified by its (relatively prime, positive, integral) multiplicities $(\alpha_1, ..., \alpha_n)$ this three-manifold can be unambiguously be described by the same data. It is typically denoted $\Sigma(\alpha_1, ..., \alpha_n)$. The convention that it correspond to a line bundle with negative degree (unambiguouisly) gives $\Sigma(\alpha_1, ..., \alpha_n)$ the additional structure of a Seifert fibration.

The quantity

$$\lfloor -\frac{\chi(\Sigma)}{2} \rfloor = \lfloor g - 1 + \frac{1}{2}\sum(\frac{1}{\alpha_i} - 1)\rfloor$$

gives the maximum dimension of any irreducible critical manifold. In particular, since for any homology sphere each $\alpha_i \geq 2$, the above result says that a Seifert fibered homology sphere with fewer than five singular five singular fibers has a discrete critical manifold. (When we allow five singular fibers, some homology spheres, e.g. $\Sigma(2,3,5,7,11)$, still have discrete critical manifolds while others, e.g. $\Sigma(3,5,7,11,13)$ do not.)

Note that these critical manifolds are always rational; indeed they are projective spaces. (It is interesting to compare this with the analogous statement for $SU(2)$ representation varieties, see [10].)

When all the critical manifolds are isolated points we can proceed as in the classical Floer picture (see for example [8]), to form a relatively graded chain complex, the *Seiberg-Witten Floer complex* which depends on the Seifert fibration

$$\pi\colon Y \to \Sigma$$

and the $\mathrm{Spin}_c(3)$ structure on $Y$, $(\mathrm{CF}_*(Y, W), \partial)$. The underlying group is a free Abelian group generated by the irreducible critical points $\{\mathfrak{b}\}$ of cs. The relative grading is defined by

$$\mathrm{gr}(\mathfrak{a}) - \mathrm{gr}(\mathfrak{b}) = \text{e-dim}\mathcal{M}(\mathfrak{a}, \mathfrak{b}),$$

the expected dimension of the moduli space of parameterized flow lines from $\mathfrak{a}$ to $\mathfrak{b}$. The moduli spaces of flows always have even expected dimension, since they are identified with moduli spaces of divisors. This says that any two critical points have even relative grading, so the boundary map must be identically zero, so that the homology groups of this complex the *irreducible Seiberg-Witten Floer homology*, denoted $\mathrm{HF}^{irr}_*(Y,W)$ is, as a graded group, isomorphic to the Seiberg-Witten Floer chain complex. In general, these homology groups depend the auxiliary data – the metric on $Y$ and the choice of connection used on $TY$ (choices which we have made in this paper); in particular they are not topological invariants of the three-manifold underlying $Y$. However, these groups enjoy certain functorial properties (product formulas) useful for calculating Seiberg-Witten invariants for four-manifolds glued along Seifert fibered spaces. According to recent work ([21] and [25]), they can also be applied to the problem of studying the fillable contact structures on these three-manifolds.

Some sample irreducible Seiberg-Witten Floer Homology are displayed, for the Seifert fibered homology spheres chosen in [7]. (The relative grading is normalized so that the solutions with nowhere vanishing spinor lie in degree zero.)

$$\mathrm{HF}^{irr}_{2i}(\Sigma(2,3,6k-1)) = \begin{cases} \mathbb{Z}^{2\lfloor \frac{k}{2}\rfloor} & \text{if } i=0, \\ 0 & \text{otherwise;} \end{cases}$$

$$\mathrm{HF}^{irr}_{2i}(\Sigma(2,3,6k+1)) = \begin{cases} \mathbb{Z}^{2\lfloor \frac{k+1}{2}\rfloor} & \text{if } i=0, \\ 0 & \text{otherwise;} \end{cases}$$

$$\mathrm{HF}^{irr}_{2i}(\Sigma(2,5,10k-1) = \begin{cases} \mathbb{Z}^{2} & \text{if } 0\le i<k-1, \\ \mathbb{Z}^{2(\lfloor \frac{k}{2}\rfloor+1)} & \text{if } i=k-1, \\ 0 & \text{otherwise} \end{cases}$$

$$\mathrm{HF}^{irr}_{2i}(\Sigma(2,5,10k+1) = \begin{cases} \mathbb{Z}^{2} & \text{if } 0\le i<k, \\ \mathbb{Z}^{2\lfloor \frac{k+1}{2}\rfloor} & \text{if } i=k, \\ 0 & \text{otherwise} \end{cases}$$

$$\mathrm{HF}^{irr}_{2i}(\Sigma(2,5,10k\pm 3)) = \mathrm{HF}^{irr}_{2i}(\Sigma(2,5,10k\pm 1))$$

$$\mathrm{HF}^{irr}_{2i}(\Sigma(2,7,14k-1)) = \begin{cases} \mathbb{Z}^2 & \text{if } 2\big| i \text{ and } 0 \leq i \leq 2k-2, \\ \mathbb{Z}^2 & \text{if } 2k-2 \leq i \leq 3k-3, \\ \mathbb{Z}^{2\lfloor \frac{k}{2}+1 \rfloor} & \text{if } i = 3k-2, \\ \mathbb{Z}^{2\lfloor \frac{k}{2} \rfloor} & \text{if } i = 3k-1, \\ 0 & \text{otherwise;} \end{cases}$$

$$\mathrm{HF}^{irr}_{2i}(\Sigma(2,7,14k+1)) = \begin{cases} \mathbb{Z}^2 & \text{if } 2\big| i \text{ and } 0 \leq i \leq 2k-1, \\ \mathbb{Z}^2 & \text{if } 2k-2 \leq i \leq 3k-2, \\ \mathbb{Z}^{2\lfloor \frac{k+3}{2} \rfloor} & \text{if } i = 3k-1, \\ \mathbb{Z}^{2\lfloor \frac{k+1}{2} \rfloor} & \text{if } i = 3k-1, \\ 0 & \text{otherwise;} \end{cases}$$

$$\mathrm{HF}^{irr}_{2i}(\Sigma(2,7,14k-3)) = \begin{cases} \mathbb{Z}^2 & \text{if } 2\big| i \text{ and } 0 \leq i \leq 2k-2, \\ \mathbb{Z}^2 & \text{if } 2k-2 \leq i \leq 3k-3, \\ \mathbb{Z}^{2k} & \text{if } i = 3k-2, \\ 0 & \text{otherwise;} \end{cases}$$

$$\mathrm{HF}^{irr}_{2i}(\Sigma(2,7,14k+3)) = \begin{cases} \mathbb{Z}^2 & \text{if } 2\big| i \text{ and } 0 \leq i \leq 2k, \\ \mathbb{Z}^2 & \text{if } 2k \leq i \leq 3k-1, \\ \mathbb{Z}^{2(k+1)} & \text{if } i = 3k, \\ 0 & \text{otherwise;} \end{cases}$$

$$\mathrm{HF}^{irr}_{2i}(\Sigma(2,7,14k\pm 5)) = \mathrm{HF}^{irr}_{2i}(\Sigma(2,7,14k\pm 3))$$

$$\mathrm{HF}^{irr}_{2i}(\Sigma(3,4,12k-1)) = \begin{cases} \mathbb{Z}^2 & \text{if } 2\big| i \text{ and } 0 \leq i \leq 2k-2, \\ \mathbb{Z}^2 & \text{if } 2k-2 \leq i \leq 3k-3, \\ \mathbb{Z}^{2\lfloor \frac{k+2}{2} \rfloor} & \text{if } k = 3k-2, \\ 0 & \text{otherwise;} \end{cases}$$

$$\mathrm{HF}^{irr}_{2i}(\Sigma(3,4,12k+1)) = \begin{cases} \mathbb{Z}^2 & \text{if } 2\big| i \text{ and } 0 \leq i < 2k, \\ \mathbb{Z}^2 & \text{if } 2k \leq i \leq 3k, \\ \mathbb{Z}^{2\lfloor \frac{k+1}{2} \rfloor} & \text{if } k = 3k, \\ 0 & \text{otherwise;} \end{cases}$$

$$\mathrm{HF}^{irr}_{2i}(\Sigma(3,4,12k-5)) = \begin{cases} \mathbb{Z}^2 & \text{if } 2\big| i \text{ and } 0 \leq i \leq 2k-2, \\ \mathbb{Z}^2 & \text{if } 2k-2 \leq i \leq 3k-3, \\ \mathbb{Z}^{2k} & \text{if } k = 3k-2, \\ 0 & \text{otherwise;} \end{cases}$$

$$\mathrm{HF}^{irr}_{2i}(\Sigma(3,5,15k-2)) = \begin{cases} \mathbb{Z}^2 & \text{if } 3\big| i \text{ and } 0 \le i \le 3k-2, \\ \mathbb{Z}^2 & \text{if } 2\big| i \text{ and } 3k-4 \le i \le 5k-3, \\ \mathbb{Z}^2 & \text{if } 5k-3 \le i \le 6k-4, \\ \mathbb{Z}^{2\lfloor \frac{k+2}{2} \rfloor} & \text{if } i = 6k-3, \\ \mathbb{Z}^{2\lfloor \frac{k}{2} \rfloor} & \text{if } i = 6k-2, \\ 0 & \text{otherwise;} \end{cases}$$

$$\mathrm{HF}^{irr}_{2i}(\Sigma(3,5,15k+2)) = \begin{cases} \mathbb{Z}^2 & \text{if } 3\big| i \text{ and } 0 \le i \le 3k+1, \\ \mathbb{Z}^2 & \text{if } 2\big| i \text{ and } 3k-1 \le i \le 5k, \\ \mathbb{Z}^2 & \text{if } 5k \le i \le 6k-2, \\ \mathbb{Z}^{2\lfloor \frac{k+3}{2} \rfloor} & \text{if } i = 6k-1, \\ \mathbb{Z}^{2\lfloor \frac{k+1}{2} \rfloor} & \text{if } i = 6k, \\ 0 & \text{otherwise.} \end{cases}$$

## References

[1] M. Artin and J.-L. Verdier. Reflexive modules over rational double points. *Math. Ann.*, 270:79–82, 1985.

[2] M. F. Atiyah, V. K. Patodi, and I. M. Singer. Spectral asymmetry and Riemannian geometry, I. *Math. Proc. Camb. Phil. Soc.*, 77:43–69, 1975.

[3] S. Bando. Einstein-Hermitian metrics on non-compact Kähler manifolds. Preprint.

[4] N. Berline, E. Getzler, and M. Vergne. *Heat Kernels and Dirac Operators*. Number 298 in Grundlehren der mathematischen Wissenschaften. Springer-Verlag, 1991.

[5] S. B. Bradlow. Special metrics and stability for holomorphic bundles with global sections. *J. Differential Geometry*, 33:169–213, 1991.

[6] S. K. Donaldson and P. B. Kronheimer. *The Geometry of Four-Manifolds*. Oxford Mathematical Monographs. Oxford University Press, 1990.

[7] R. Fintushel and R. J. Stern. Seifert fibered homology three-spheres. *Proc. of the London Math. Soc.*, 61:109–137, 1990.

[8] A. Floer. An instanton-invariant for 3-manifolds. *Comm. Math. Phys.*, 119:215–240, 1988.

[9] W. Fulton. *Introduction to Toric Varieties*. Number 131 in Annals of Mathematics Studies. Princeton University Press, 1993.

[10] M. Furuta and B. Steer. Seifert fibered homology 3-spheres and the Yang-Mills equations on Riemann surfaces with marked points. *Advances in Mathematics*, pages 38–102, 1992.

[11] R. Godement. *Topologie Algébrique et Théorie des Faisceaux*. Number 1252 in Actualités Scientifiques et Industrielles. Hermann, 1964.

[12] G. Gonzalez-Sprinberg and J.-L. Verdier. Construction géométrique de la correspondance de McKay. *Ann. Scient. Éc. Norm. Sup.*, pages 409–449, 1983.

[13] H. Grauert, Th. Peternell, and R. Remmert, editors. *Several Complex Variables VII: Sheaf-Theoretical Methods in Complex Analysis*, volume 74 of *Encyclopaedia of Mathematical Sciences*. Springer-Verlag, 1991.

[14] R. C. Gunning. *Introduction to Holomorphic Functions of Several Variables*, volume III. Brooks/Cole Publishing Co., 1990.

[15] G.-Y. Guo. Yang-mills fields on cylindrical manifolds and holomorphic bundles i. *Comm. in Math. Phys.*, 179(3):737–776, 1996.

[16] R. Hartshorne. *Algebraic Geometry.* Number 52 in Graduate Texts in Mathematics. Springer-Verlag, 1977.

[17] A. Jaffe and C. H. Taubes. *Vortices and Monopoles.* Number 2 in Progess in Physics. Birkhäuser, 1980.

[18] J. Kazdan and F. W. Warner. Curvature functions for compact 2-manifolds. *Annals of Mathematics*, 99:14–47, 1974.

[19] H. Knörrer. Group representations and the resolution of rational double points. In J. McKay, editor, *Finite Groups – Coming of Age*, number 45 in Contemporary Mathematics. American Mathematical Society, 1985.

[20] D. Kotschick. The Seiberg-Witten equations on symplectic four-manifolds [after C. H. Taubes]. exposé 812, Séminar N. Bourbaki, 1995-1996.

[21] P. B. Kronheimer and T. S. Mrowka. Monopoles and contact structures. To appear in *Inventiones Math.*

[22] P. B. Kronheimer and T. S. Mrowka. The genus of embedded surfaces in the projective plane. *Math. Research Letters*, 1:797–808, 1994.

[23] P. B. Kronheimer and T. S. Mrowka. Embedded surfaces and the structure of Donaldson's polynomial invariants. *J. Differential Geometry*, pages 573–734, 1995.

[24] H. B. Lawson and M.-L. Michelsohn. *Spin Geometry.* Number 38 in Princeton Mathematics Series. Princeton University Press, 1989.

[25] P. Lisca and G. Matić. Tight contact structures and the Seiberg-Witten invariants. To appear in *Inventiones Math*, 1996.

[26] R. Lockhart and R. McOwen. Elliptic differential operators on non-compact manifolds. *Annali di Scuola Norm. Sup. de Pisa*, IV-12:409–448, 1985.

[27] J. W. Morgan, T. S. Mrowka, and D. Ruberman. *The $L^2$-Moduli Space and a Vanishing Theorem for Donaldson Polynomial Invariants.* Number II in Monographs in Geometry and Topology. International Press, 1994.

[28] J. W. Morgan, Z. Szabó, and C. H. Taubes. A product formula for Seiberg-Witten invariants and the generalized Thom conjecture. Preprint.

[29] T. S. Mrowka. Seiberg-Witten invariants for four-manifolds. Lectures at Harvard University.

[30] P. Scott. The geometries of 3-manifolds. *Bull. London Math. Soc.*, 15:401–487, 1983.

[31] H. Seifert. Topologie dreidimensionaler gefaserter Räume. *Acta Math.*, 60:147–238, 1932.

[32] C. H. Taubes. *$L^2$-Moduli spaces on open 4-manfolds.* Number I in Monographs in Geometry and Topology. International Press, 1993.

[33] C. H. Taubes. The Seiberg-Witten invariants and symplectic forms. *Math. Research Letters*, 1:809–822, 1994.

[34] E. Witten. Monopoles and four-manifolds. *Math. Research Letters*, 1:769–796, 1994.

Department of Mathematics, Massachusetts Institute of Technology, Cambridge MA 02139

The Mathematical Sciences Research Institute, 1000 Centennial Dr., Berkeley, Ca 94720